\documentclass[12pt]{article}
\usepackage{amsmath,amsfonts,amssymb,mathabx,txfonts}
\usepackage{mathrsfs,upgreek}
\usepackage{dsfont}
\usepackage[utf8]{inputenc}
\usepackage[T1]{fontenc}
\usepackage{lmodern} \usepackage{ucs}
\usepackage{fullpage}

\usepackage{graphicx}
\usepackage{caption}
\usepackage{subcaption}

\usepackage{array}
\usepackage{multicol}
\usepackage{multirow}
\usepackage{color}
\usepackage{authblk}
\usepackage{esint}
\usepackage{enumitem,hyperref}

\usepackage{accents}

\usepackage{tikz}
\usetikzlibrary{decorations}
\usetikzlibrary{decorations.pathreplacing}
\usetikzlibrary{arrows}

\newtheorem{lemma}{Lemma}[section]
\newtheorem{theo}[lemma]{Theorem}
\newtheorem{rmk}[lemma]{Remark}
\newtheorem{proposition}[lemma]{Proposition}

\newtheorem{coro}[lemma]{Corollary}
\newtheorem{example}[lemma]{Example}

\newcommand*{\overbigdot}[1]{%
   \accentset{\mbox{\large\bfseries .}}{#1}}

\makeatletter
\def\namedlabel#1#2{\begingroup
    #2%
    \def\@currentlabel{#2}%
    \phantomsection\label{#1}\endgroup
}
\makeatother

\makeatletter
\renewcommand*{\eqref}[1]{%
  \hyperref[{#1}]{\textup{\tagform@{\ref*{#1}}}}%
}
\makeatother

\newenvironment{Proof}{\noindent
    \abovedisplayskip = 0.5\abovedisplayskip
    \belowdisplayskip=\abovedisplayskip{\bfseries Proof. }}{\QED\medskip}

    \newenvironment{ProofOf}[1]{\noindent
    \abovedisplayskip = 0.5\abovedisplayskip
    \belowdisplayskip=\abovedisplayskip{\bfseries Proof of  #1. }}{\QED\medskip}

\newcommand{\QED}{\mbox{}\hfill \raisebox{-0.2pt}{\rule{5.6pt}{6pt}\rule{0pt}{0pt}} \medskip\par}
\newcommand{\N}{\mathbb{N}}
\newcommand{\Z}{\mathbb{Z}}
\newcommand{\R}{\mathbb{R}}
\newcommand{\C}{\mathbb{C}}
\newcommand{\T}{\mathbb{T}}

\newcommand{\Ker}{\mathrm{Ker}}
\newcommand{\ds}{\displaystyle}
\newcommand{\ud}{\, {\mathrm{d}}}
\newcommand{\avex}[1]{\big\langle{#1} \big\rangle_{\mathbb T^d}}

\newcommand{\Td}{\mathbb T^d}

\title{Plane wave stability analysis of Hartree and quantum dissipative  systems}

\author{Thierry~Goudon\thanks{ {Corresponding author, \tt thierry.goudon@inria.fr}}}
\author{Simona Rota Nodari\thanks{ {\tt simona.rotanodari@univ-cotedazur.fr}}}
\date{}

\affil[1]{\small Universit\'e C\^ote d'Azur, Inria,  CNRS, LJAD,

Parc Valrose, F-06108 Nice, France}

\begin{document}
\maketitle

\begin{abstract}
We investigate the stability of plane wave solutions 
of equations describing quantum particles interacting with a complex environment.
The models take the form of PDE systems with a non local (in space or in space and time)
self-consistent potential; such a coupling lead to challenging issues 
compared to the usual non linear Schr\"odinger equations. 
The analysis relies on the identification of  suitable Hamiltonian structures
and Lyapounov functionals.
We point out analogies and differences 
between the original model, involving a coupling with a wave equation, 
and its asymptotic counterpart obtained in the large wave speed regime.
In particular, while the analogies provide interesting intuitions, our analysis shows that it is illusory to 
obtain 
results on the former based on  a perturbative analysis from the latter. 
\end{abstract}

\vspace*{.5cm}
{\small
Open quantum systems. Particles interacting with a vibrational field. Schr\"odinger-Wave equation. Plane wave. Orbital stability.
\\[.4cm]

\noindent{\bf Math.~Subject Classification.} 
35Q40 %PDEs in connection with quantum mechanics 
35Q51 % Soliton-like equations
35Q55 % NLS-like equations (nonlinear Schrödinger) }

\section{Introduction}

This work is concerned with the stability analysis of certain solutions of the following 
Hartree-type equation
\begin{subequations}
    \begin{alignat}{1}
        \label{Hartree}&\ds i\partial_{t}U+\frac{1}{2}\Delta_{x}U =\gamma\left(\sigma_{1}\star_{x}\int_{\mathbb R^n}\sigma_{2} \Psi\ud z\right)U,\\
        \label{Poissonz_Hartree}&\ds -\Delta_{z}\Psi=-\gamma\sigma_{2}(z)\left(\sigma_{1}\star_{x}|U|^{2}\right)(x)    \end{alignat}
\end{subequations}
endowed with the initial condition
\begin{equation}\label{Ci_Ha}
U\big|_{t=0}=U^{\mathrm{Init}},
\end{equation}
and of the following
Schr\"odinger-Wave system:
\begin{subequations}
  \begin{alignat}{1}
\label{Schro-s}
 i\partial_t U +\ds\frac{ 1}{2} \Delta_x U=\gamma \Phi U,
\\
%[.3cm]
\label{Schro-w}\ds\frac1{c^2}\partial_{tt}^2\Psi-\Delta _z \Psi=-\gamma\sigma_2(z)\sigma_1\star  |U|^2(t,x),
\\
%[.3cm]
\label{Schro-p}\Phi(t,x)=\ds\iint_{\mathbb T^d\times\mathbb R^n}\sigma_1(x-y)\sigma_2(z)\Psi(t,y,z)\ud z\ud y,
  \end{alignat}
\end{subequations}
where $\gamma,c>0$ are given positive parameters,
completed with 
\begin{equation}\label{Ci_SW}
U\big|_{t=0}=U^{\mathrm{Init}}, \qquad \Psi\big|_{t=0}=\Psi^{\mathrm{Init}},\qquad
\partial_t\Psi\big|_{t=0}=\Pi^{\mathrm{Init}}.
\end{equation}
The variable $x$ lies in the torus $\mathbb T^d$, meaning that the equations are understood with $(2\pi)-$periodicity in all directions. In \eqref{Schro-w}, the additional variable $z$ lies in $\mathbb R^n$ and, as explained below, it is crucial to assume $n\geq 3$.
For reader's convenience, the scaling of the equation is fully detailed in Appendix~\ref{adim}; for our purposes
the God-given form functions $\sigma_1,\sigma_2$ are fixed once for all and the features of the coupling 
are embodied in the parameters $\gamma, c$.
The system \eqref{Hartree}-\eqref{Poissonz_Hartree} can be obtained,
at least formally, from \eqref{Schro-s}-\eqref{Schro-p} by  letting the parameter $c$
run to $+\infty$, while $\gamma$ is kept fixed.
By the way, system  \eqref{Hartree}-\eqref{Poissonz_Hartree} can be cast in the more usual form 
\begin{equation}\label{hartree}
    i\partial_{t}U+\frac{1}{2}\Delta_{x}U=-\gamma^2\kappa\left(\Sigma\star_{x}|U|^{2}\right)U,\hspace{1cm}t\in\R,\;x\in\R^{d}.
\end{equation}
where\footnote{The Fourier transform of an integrable  function $\varphi: \mathbb R^n\to \mathbb C$ is defined by 
$\widehat
\varphi(\xi)=\int_{\mathbb R^n} \varphi(z)e^{-i\xi\cdot z} \ud z$.}
\begin{equation}\label{kappa1}
\kappa=\ds\int_{\mathbb R^n} \sigma_2(z)(-\Delta_z)^{-1} \sigma_2(z)\ud z
=\ds\int_{\mathbb R^n} \ds\frac{|\widehat\sigma_2(\xi)|^2}{|\xi|^2}\ds\frac{\ud \xi}{(2\uppi)^n}>0 \text{ and } \Sigma=\sigma_1\star \sigma_1.\end{equation}
Letting now $\Sigma$ resemble the delta-Dirac mass, the asymptotic leads to 
the standard cubic non linear Schr\"odinger equation
\begin{equation}\label{nls}
    i\partial_{t}U+\frac{1}{2}\Delta_{x}U=-\gamma^2\kappa |U|^{2} U.
\end{equation}
 in the \emph{focusing} case.
These asymptotic  connections  can be expected to shed some light on the dynamics of \eqref{Schro-s}-\eqref{Schro-p} 
and to be helpful to guide the intuition about the behavior of the solutions, see \cite{Vi4,Vi3}.  
\\

The motivation for investigating these systems  takes its roots in the general landscape of the  
analysis of ``open systems'', describing the dynamics of particles  driven by 
 momentum and energy exchanges with a complex environment.
 Such problems are modeled as Hamiltonian systems, 
 and it is expected that the interaction mechanisms ultimately produce 
  the dissipation of the particles' energy, an  idea which  dates back to A.~O.~Caldeira and A.~J.~Leggett \cite{CL}. 
  These issues have been investigated 
  for various classical and quantum couplings, and 
  with many different mathematical viewpoints, see e.~g.~\cite{bach,Fau,JP1,JP2,KKS,KKSb,bDPL}. 
The case  in which the environment
is described as a vibrational field, like in the definition of the potential by 
\eqref{Schro-w}-\eqref{Schro-p}, is particularly appealing.
In fact, \eqref{Schro-s}-\eqref{Schro-p}
is a quantum version 
of a model introduced by S.~De Bi\`evre and L.~Bruneau, dealing with a single classical particle \cite{BdB}.
Intuitively, the model of \cite{BdB} can be thought of  
as if in each space position $x\in\R^{d}$ 
 there is a membrane  oscillating in a  direction $z\in\R^{n}$, transverse to the motion of the particles. When a particle hits a membrane,  its kinetic energy  activates
 vibrations and the energy is evacuated at infinity in the $z-$direction. 
These energy transfer mechanisms eventually act as a sort of friction force
  on the particle, an intuition rigorously justified in 
   \cite[Theorem~2 and Theorem~4]{BdB}. 
   We refer the reader to 
   \cite{AdBLP,dBP,dBPS,bDPL,dBS} for further theoretical and numerical insight about this model.
    The model of \cite{BdB} has been revisited by considering  \emph{many} interacting particles, which leads to 
    Vlasov-type equations, still coupled to 
    a wave equation for defining the potential \cite{GV}.  
    Unexpectedly, asymptotic arguments indicate a connection with the \emph{attractive} Vlasov-Poisson dynamic 
\cite{BGV}. In turn, the particles-environment interaction 
can be interpreted in terms of Landau damping \cite{Vi1,Vi2}.
The quantum version \eqref{Schro-s}-\eqref{Schro-p} of the De Bi\`evre-Bruneau model has been discussed in \cite{Vi4,Vi3}, 
with a connection to the kinetic model by means of a semi-classical analysis inspired from \cite{Lions_Paul_Wigner}.
Note that in \eqref{Schro-s}-\eqref{Schro-p}, the vibrational field remains of classical nature; 
a fully quantum 
framework is dealt with in \cite{Fau,DuSh} for instance.
\\

A remarkable feature of these systems 
is the presence of conserved quantities,
here inherited from the framework designed in  \cite{BdB} for a classical particle, 
and the study of these models brings out the critical role of the wave speed $c>0$ and the   
dimension $n$ of the space for the wave equation (we can already notice that  $n\geq 3$ is necessary for  \eqref{kappa1} to be meaningful), see \cite{BdB,Vi1,Vi2,Vi3}.
For the Schr\"odinger-Wave system \eqref{Schro-s}-\eqref{Schro-p} 
the energy
\begin{equation}\label{eq:hamiltonian}
H_{SW}(U,\Psi,\Pi)=
\frac{ 1}{4}\int_{\Td}|\nabla  U |^2 \ud x+ \iint_{\Td\times\mathbb R^n} \left(\ds{\color{black} c^2}\Pi^2+{\color{black}\frac{1}{4}}|\nabla_z \Psi|^2\right)\ud x\ud z+\frac{\gamma}{2}\int_{\Td}\Phi |U|^2\ud x,
\end{equation}
is conserved since we can readily check that 
\[
\ds\frac{\ud }{\ud t} H_{SW}(U,\Psi,{\color{black}-\frac{1}{2c^2}}\partial_t \Psi)=0.
\]
Similarly, for the Hartree system \eqref{Hartree}-\eqref{Poissonz_Hartree}, we get 
\[\ds\frac{\ud }{\ud t} H_{Ha}(U)=0\]
where we have set \[
H_{Ha}(U)=\frac{ 1}{4}\int_{\Td}|\nabla  U |^2 \ud x
-\gamma^2\frac{\kappa}{4}\int_{\Td}\Sigma(x-y) |U(t,x)|^2|U(t,y)|^2\ud y\ud x
.\]
Furthermore, for both model,  the $L^2$ norm is conserved.
Of course, these conservation properties play a central role for the analysis of the equations.
However, \eqref{Hartree}-\eqref{Poissonz_Hartree} has further fundamental properties  
which occur only for the asymptotic model:
 firstly, \eqref{Hartree}-\eqref{Poissonz_Hartree} is 
 Galilean invariant, which means that, given a solution $(t,x)\mapsto u(t,x)$ and for any $p_0\in \mathbb T^d$,
 the function $(t,x)\mapsto u(t,x-tp_0)e^{i(x-tp_0/2)}$ is a solution too; secondly, the momentum
  $p(t)=\mathrm{Im}\int\bar u(t,x) \nabla_x u(t,x)\ud x$ is conserved and, accordingly, 
  the center of mass follows a straight
line at constant speed.
That these properties are not satisfied by the more complex system
 \eqref{Schro-s}-\eqref{Schro-p} makes its analysis more challenging.
Finally, we point out that, in contrast to the usual 
nonlinear Schr\"odinger equation or Hartree-Newton system, where $\Sigma$ is the Newtonian potential,  the equations \eqref{Hartree}-\eqref{Poissonz_Hartree} or 
\eqref{Schro-s}-\eqref{Schro-p}
do not fulfil a scale invariance property. This also leads to specific mathematical difficulties: despite the possible regularity
of $\Sigma$, many results  and approaches 
of the Newton case do not extend to a general kernel, due to the lack of scale invariance.

When the problem is set on the whole space $\mathbb R^d$, 
one is interested in the stability of solitary waves, which are  solutions
of the equation with the specific  form $u(t,x)=e^{i\omega t} Q(x)$, and, for \eqref{Schro-s}-\eqref{Schro-p}, $\psi(t,x,z)=\Psi(x,z)$.
The details of the solitary 
wave are embodied into the Choquard equation, satisfied by the profile $Q$,
\cite{Lieb_Choquard, Lions_Choquard}.
It turns out that the Choquard equation have  infinitely many solutions; 
among these solutions,
it is relevant to select the solitary wave
which minimizes
 the energy functional under a mass constraint,  \cite{Lieb_Choquard,MaZ}
and to study the orbital stability 
of this  minimal energy state.
This program has been investigated  for \eqref{nls} and \eqref{Hartree}-\eqref{Poissonz_Hartree}
 in the specific case where $\Sigma(x)=\frac{1}{|x|}$ in dimension $d=3$, by various approaches \cite{Cazenave_Lions_orbital_stab, Lenzmann, PLLcc1, PLLcc2,MaMe,Weinstein1,Weinstein2}.
 Quite surprisingly, the specific form of the 
  potential plays a critical role in the analysis (either through explicit formula or through scale invariance properties),
  and dealing with a general convolution kernel, as smooth as it is, leads to new difficulties, that can be treated 
  by a perturbative argument, see  \cite{jmaa,Zhang} for the case of the Yukawa potential, and \cite{Vi3} for \eqref{Hartree}-\eqref{Poissonz_Hartree} and
\eqref{Schro-s}-\eqref{Schro-p}.
\\

Here, we adopt a different viewpoint. We consider the case where the problem holds on the torus $\mathbb T^d$, 
and we are specifically interested in the stability of \emph{plane wave solutions}
of \eqref{Schro-s}-\eqref{Schro-p} and \eqref{Hartree}-\eqref{Poissonz_Hartree}. We refer the reader to 
\cite{SRN1,SRN2,Faou, Keller} for results on the nonlinear Schr\"odinger equation \eqref{nls} in this framework. 
The discussion on the stability of these plane wave solutions will make the following smallness condition
\begin{equation}\label{small}
4 \gamma^2 \kappa 
 \|\sigma_{1} \|_{L^1}^2<1
\end{equation} 
(assuming the plane wave has an amplitude unity) appear. 
Despite its restriction to the periodic framework, the interest of this study is two-fold:
on the one hand, it points out some difficulties specific to the coupling and provides useful hints for future works;
on the other hand, it clarify the role of the parameters, by making   
stability conditions \emph{explicit}.
\\

The paper is organized as follows. In Section~\ref{SetUp}, we clarify 
the positioning of the paper. To this end, we further discuss some mathematical features of the model.
We also introduce the main assumptions on the parameters  that will be used throughout the paper and 
we provide an overview of the results. Section~\ref{S:Hartree} is concerned with the stability  analysis 
of the Hartree equation \eqref{Hartree}-\eqref{Poissonz_Hartree}.
Section~\ref{S:SW0} deals with the Schr\"odinger-Wave system
at the price of restricting to the case where the wave vector of the plane wave solution vanishes: $k=0$.
For reasons explained in details below, the general case is much more difficult.
Section~\ref{S:kn0} justifies that in general the mode $k\neq 0$ is linearly 
and orbitally
unstable. The proof splits into two steps. The former is concerned by the spectral instability; it relies on a suitable reformulation of the linearized operator, which allows us
to count indirectly the eigenvalues. The latter step proves instability by using a contradiction argument and estimates established through the Duhamel formula.
Finally, in Appendix \ref{adim}, we provide a physical interpretation of the parameters involved,  and for the sake of completeness, in Appendices~\ref{AppB} and \ref{AppC}, we discuss the well-posedness of the Schr\"odinger-Wave system \eqref{Schro-s}-\eqref{Schro-p} and its link with the Hartree equation \eqref{Hartree}-\eqref{Poissonz_Hartree} in the regime of large $c$'s.

\section{Set up of the framework}
\label{SetUp}

\subsection{Plane wave solutions and dispersion relation}
\label{PWsol}

For any $k\in \Z^d$, we start by seeking  solutions to \eqref{Schro-s}-\eqref{Schro-p} of the form
\begin{equation}\label{eqrelequiSW}
 U(t,x)=U_k(t,x):=\exp\big(i(\omega t+k\cdot x)\big),\quad  
  \Psi(t,x,z)=\Psi_*(z), \quad \partial_t\Psi(t,x,z)={\color{black} -2c^2}\Pi_* (z)=0,
\end{equation}
with $\omega\in \mathbb R$. Note that the $L^2$ norm of $U_k$ is $(2\uppi)^{d/2}$ and $\Psi_*$ actually does not depend on the time variable, nor on $x$.
Since $|U_k(t,x)|=1$ 
is constant, 
the wave equation simplifies to 
$$ {\frac{1}{c^2}}\partial_{tt}^2\Psi-\Delta _z \Psi=-\gamma \sigma_2(z) 
\avex{ \sigma_1},
$$
where $\avex{\cdot}$ stands for the average over $\mathbb T^d$: $\avex{f}=\int_{\mathbb T^d} f(x)\ud x$. As a consequence, 
$z\mapsto \Psi_*(z)$ is a solution to \eqref{Schro-w} if 
\begin{equation*}
  \Psi_*(z)=-\gamma 
  \Gamma(z) \avex{\sigma_1},
\end{equation*}
with $\Gamma$ the solution of 
\[-\Delta_z \Gamma(z) =\sigma_2(z).\]
This auxiliary function $\Gamma$ is thus defined by the convolution of $\sigma_2$ with the elementary solution of the Laplace operator in dimension $n$,
or equivalently by means of Fourier transform: 
\begin{equation}\label{formGamma}
\Gamma(z)=\int_{\mathbb R^n}
\frac{C_n}{|z-z'|^{n-2}}\sigma_2(z')\ud z'=
\mathscr F^{-1}_{\xi\to z}\Big(\ds\frac{\widehat \sigma_2(\xi)}{|\xi|^2}\Big) .
\end{equation}
The corresponding potential \eqref{Schro-p}
 is actually a constant which 
reads
\[
 -\gamma 
\ds\iint_{\mathbb T^d\times\mathbb R^n} \sigma_1(x-y)
\sigma_2(z)\Gamma(z)
 \avex{\sigma_1}\ud z\ud y
=-\kappa\gamma 
\avex{\sigma_1}^2
\]
with $$\kappa=\ds\int_{\mathbb R^n} \sigma_2(z)\Gamma(z)\ud z=\ds\int_{\mathbb R^n} |\nabla_z\Gamma(z)|^2\ud z>0$$
(we remind the reader  that this formula coincides with \eqref{kappa1} and makes sense only when $n\geq 3$).
It remains to identify the condition on the coefficients so that $U_k$ satisfies the Schr\"odinger equation \eqref{Schro-s}: this leads to the following dispersion relation
 \begin{equation}\label{eq:dispersion}
\omega+\frac{k^2 }2-
\Upsilon_*=0,\qquad 
\Upsilon_*=\gamma^2 \kappa\avex{\sigma_1}^2
>0
\end{equation}
with $k^2=\sum_{j=1}^d k_j^2$.
We can compute explicitly the associated energy:
\[
H_{SW}(U_k,\Psi_*,\Pi_*)
=\frac{(2\uppi)^d}{2}
\left(
\ds\frac{k^2}{2}
-\ds\frac{ \gamma^2 \kappa}{2} 
\avex{\sigma_1}^2\right)
=\ds\frac{(2\uppi)^d }{4}(
k^2-\Upsilon_* )
.\]
Of course, among these solutions, the constant mode $U_0(t,x)=e^{i\omega t}{\mathbf 1(x)}$ has minimal energy.
\\

It turns out that the plane wave
$
U_k(t,x)=e^{i\omega t}e^{ik\cdot x}$
equally satisfies \eqref{Hartree}-\eqref{Poissonz_Hartree} provided the dispersion relation \eqref{eq:dispersion}
holds. Incidentally, we can  check that 
\[
H_{Ha}(U_k)
=\frac{(2\uppi)^d}{2}
\left(
\ds\frac{k^2}{2}
-\ds\frac{ \gamma^2 \kappa}{2}  
\avex{\Sigma}\right)
=\ds\frac{(2\uppi)^d  
}{4}(
k^2-\Upsilon_*  
)
\]
is made minimal 
when $k=0$.

\subsection{Hamiltonian structure and symmetries of the problem}

The conservation properties 
play a central role in the stability analysis, for instance 
in the reasonings that use
concentration-compactness arguments \cite{Cazenave_Lions_orbital_stab}.
Based on the conserved quantities, one can try to construct 
a Lyapounov functional, intended to evaluate how far a solution is from an equilibrium state.
Then the stability analysis relies on the ability to prove a coercivity estimate
on the variations 
of the Lyapounov functional, see \cite{Tao,Weinstein1,Weinstein2}.
This viewpoint can 
be further extended by identifying analogies with finite dimensional Hamiltonian systems with symmetries, 
which has permitted to set up a quite general framework \cite{GSS,GSS2}, revisited recently in \cite{SRN1}.
The strategy relies on the ability in 
 exhibiting 
a  Hamiltonian formulation of the problem
\[
\partial_t X=\mathbb J\partial_X \mathscr H(X),\]
where  the symplectic structure is given by the skew-symmetric operator  $\mathbb J$.
 As a consequence of Noether’s Theorem, this formulation encodes the conservation properties of the system. In particular, it implies that $t\mapsto \mathscr H(X(t))$ is a conserved quantity.
For the problem under consideration, as it will be detailed below, $X$ is a vectorial unknown
with components possibly
depending on different variables ($x\in \mathbb T^d$ and $z\in \mathbb R^n$).
This induces specific difficulties, in particular because the nature of the coupling is non local  
and delicate spectral issues arise related to the essential spectrum of the wave equation in $\mathbb R^n$.
Next, we can easily observe that the systems \eqref{Hartree}-\eqref{Poissonz_Hartree} and
\eqref{Schro-s}-\eqref{Schro-p} are invariant under multiplications by a phase factor of $U$, the 
``Sch\"odinger unknown'', and under translations in the $x$ variable. This leads to the conservation of the $L^2$ norm of $U$ and of the total momentum. However,
the systems \eqref{Hartree}-\eqref{Poissonz_Hartree} and
\eqref{Schro-s}-\eqref{Schro-p} cannot be handled by a direct application of the results in \cite{SRN1,GSS,GSS2}: the basic assumptions are simply not satisfied. Nevertheless, our approach is strongly inspired from  \cite{SRN1,GSS,GSS2}.
As we will see later, for the Hartree system, a decisive advantage comes from the 
conservation of the total momentum and the Galilean invariance of the problem. 
For the Schr\"odinger-Wave problem, since the expression of the total momentum mixes up contribution from the 
``Schr\"odinger unknown'' $U$  and the ``wave unknown'' $\Psi$, the information on its conservation does not seem readily useful. \footnote{For the problem set on $\mathbb R^d$, it is still possible, in the spirit of results obtained in \cite{Faou} for NLS, to justify that orbital stability holds on a finite time interval: the solution remains at a distance $\epsilon$ from the orbit of  the ground state over time interval of order $\mathscr O(1/\sqrt\epsilon)$, see \cite[Theorem~4.2.11 \& Section~4.6]{LV_PhD}. The argument relies on the dispersive properties of the wave equation through Strichartz' estimates.}

In what follows, we find advantages in  changing the unknown by writing 
$U(t,x)=e^{ik\cdot x} u(t,x)$; in turn the Schr\"odinger equation 
$i\partial_ t U+\frac12\Delta U=\Phi U$ becomes
\[i\partial_ t u+\frac12\Delta u- \frac{k^2}{2}u + ik\cdot \nabla u =\Phi u.\]
Accordingly, the parameter $k$ will appear in the definition 
the 
energy functional $\mathscr H$.
This explains a major difference between  \eqref{Hartree}-\eqref{Poissonz_Hartree} and
\eqref{Schro-s}-\eqref{Schro-p}: for 
the former, a coercivity estimate can be obtained for   the energy functional $\mathscr H$, for the latter, when $k\neq 0$
there are terms which cannot be controlled easily.
This is reminiscent of 
the momentum conservation in \eqref{Hartree}-\eqref{Poissonz_Hartree} and 
the lack of Galilean invariance for \eqref{Schro-s}-\eqref{Schro-p}.
The detailed analysis of the linearized 
 operators sheds more light on the different behaviors of the systems   \eqref{Hartree}-\eqref{Poissonz_Hartree} 
 and  \eqref{Schro-s}-\eqref{Schro-p}.

\subsection{Outline of the main results}

Let us collect the assumptions on the form functions $\sigma_1$ and $\sigma_2$ 
that govern the coupling:
\begin{itemize}
\item [\namedlabel{H1}{\textbf{(H1)}}] $\sigma_1:\mathbb T^d\rightarrow [0,\infty)$ is $C^\infty$ smooth, radially symmetric; $\avex{\sigma_1}\neq 0$;
\item [\namedlabel{H2}{\textbf{(H2)}}] $\sigma_2:\mathbb R^n\rightarrow [0,\infty)$ is $C^\infty$ smooth, radially symmetric and compactly supported;
 \item [\namedlabel{H3}{\textbf{(H3)}}] 
 $(-\Delta)^{-1/2}\sigma_2\in L^2(\mathbb R^n)$;
\item [\namedlabel{H4}{\textbf{(H4)}}] for any $\xi\in \mathbb R^n$, $\widehat \sigma_2(\xi)\neq 0$.
\end{itemize}
Assumptions \ref{H1}-\ref{H2} are natural in the framework introduced in \cite{BdB}.
Hypothesis \ref{H3} 
can equivalently  be rephrased as
$(-\Delta)^{-1}\sigma_2\in\overbigdot{H}^1(\mathbb R^n)$;
 it appears in many places of the analysis of such coupled systems and,  at least, 
it makes the constant 
$\kappa$ in \eqref{kappa1} meaningful. This constant is a component of the stability constraint \eqref{small}.
Hypothesis 
\ref{H4} equally appeared in \cite[Eq.~(W)]{BdB} when discussing large time asymptotic issues.  
Assumptions \ref{H1}-\ref{H4} are assumed throughout the paper. \\

Our results can be summarized as follows. 
We assume \eqref{small} and consider $k\in \mathbb Z^d$ and $\omega\in \mathbb R$ satisfying \eqref{eq:dispersion}.
For the Hartree equation, the analysis is quite complete:
\begin{itemize}
\item the plane wave $e^{i(\omega t+k\cdot x)}$ is spectrally stable (Theorem~\ref{spectralLkH});
\item for any initial perturbation with zero mean, the solutions of the linearized Hartree equation are $L^2$-bounded, uniformly over $t\geq 0$ (Theorem~\ref{Prop:Hlstab});
\item the plane wave $e^{i(\omega t+k\cdot x)}$ is orbitally stable (Theorem~\ref{prop:orbital_Hartree_k}).
\end{itemize}
For the Schr\"odinger-Wave system, %only 
the case $k=0$ is fully addressed as follows:
%{\color{blue} Modifier}
\begin{itemize}
\item  the plane wave $(e^{i\omega t}\mathbf 1(x), -\gamma\Gamma(z)\avex{\sigma_1}, 0)$ is spectrally stable (Corollary \ref{spectralstabSWk0});
\item for any initial perturbation of $(e^{i\omega t}\mathbf 1(x), -\gamma\Gamma(z)\avex{\sigma_1}, 0)$ with zero mean, the solutions of the linearized Schr\"odinger-Wave system are $L^2$-bounded, uniformly over $t\geq 0$ (Theorem~\ref{LinStabk0SW});
\item the plane wave $(e^{i\omega t}\mathbf 1(x), -\gamma\Gamma(z)\avex{\sigma_1}, 0)$ is orbitally stable (Theorem~\ref{prop:orbital_SW_k}).
\end{itemize}
When $k\neq 0$, 
the situation is much more involved; at least we prove that
in general the plane wave solution $(e^{i(\omega t+k\cdot x)},-\gamma \Gamma(z)\avex{\sigma_1},0)$ is spectrally unstable, see Section~\ref{S:kn0}  and Corollary \ref{spectralinstabSWk},  and orbitally unstable, see Theorem~\ref{Th:unsta}. 

Finally, let us mention that the approach presented here has been developed on an even simpler model, where the Schr\"odinger equation is replaced by a mere finite dimensional differential system \cite{SRN3}. 

 \section{Stability analysis of the Hartree system \eqref{Hartree}-\eqref{Poissonz_Hartree}}
\label{S:Hartree}

To study the stability of the plane wave solutions of the Hartree system, it is useful to write the solutions of \eqref{Hartree}-\eqref{Poissonz_Hartree} in the form
\begin{equation*}
  U(t,x)=e^{i k\cdot x}u(t,x)
\end{equation*}
with $u(t,x)$ solution to 
\begin{equation}
  \label{Hartree_k}
  i\partial_t u+\ds\frac12\Delta u-\ds\frac{k^2}{2}u+ik\cdot \nabla u=-\gamma^2 \kappa(\Sigma\star |u|^2)u.
\end{equation}
If $k\in \Z^d$ and $\omega\in \mathbb R$ satisfy the dispersion relation \eqref{eq:dispersion}, $u_\omega(t,x)=e^{i\omega t}\mathbf 1(x)$ is a solution to \eqref{Hartree_k} with initial condition $u_\omega(0,t)=\mathbf 1(x)$.
Therefore, studying the stability properties of $U_k(t,x)=e^{i\omega t}e^{ik\cdot x}$ as a solution to \eqref{Hartree}-\eqref{Poissonz_Hartree} amounts to studying the stability of $u_\omega(t,x)=e^{i\omega t}\mathbf 1(x)$ as a solution to \eqref{Hartree_k}.

The problem \eqref{Hartree_k} has an Hamiltonian symplectic structure  when considered on the \emph{real} Banach space $H^1(\T^d;\R)\times H^1(\T^d;\R)$.
Indeed, if we write $u=q+ip$, with $p,q$ real-valued, we obtain 
\[\partial_t 
\begin{pmatrix}
q
\\
p\end{pmatrix}
=\mathbb  J \nabla_{(q,p)}\mathscr H(q,p)
\]
with 
\[
\mathbb J=\begin{pmatrix}
0 & 1 
\\
-1 & 0
\end{pmatrix}\]
and
\begin{align*}
  \mathscr H(q,p)=&\,
\frac12\left(\ds\frac12\ds\int_{\mathbb T^d}  
|\nabla q|^2+|\nabla p|^2 \ud x+ \frac{k^2}{2}\int_{\mathbb T^d}(p^2+q^2)
\ud x-\int_{\T^d}pk\cdot \nabla q\ud x+\int_{\T^d}qk\cdot \nabla p\ud x\right)\\
&-\frac{\gamma^2 \kappa}{4} \ds\int_{\mathbb T^d}
\Sigma\star (p^2+q^2)(p^2+q^2)\ud x.
\end{align*}
Coming back to $u=q+ip$, we can write
\begin{align}\label{ham_Hartree_k}
  \mathscr H(u)=&\,
  \frac12\left(\ds\frac12\ds\int_{\mathbb T^d}  
  |\nabla u|^2\ud x+ \frac{k^2}{2}\int_{\mathbb T^d}|u(x)|^2
  \ud x+\int_{\T^d}k\cdot(-i \nabla u)\overline{u}\ud x\right)\nonumber\\
  &-\frac{\gamma^2 \kappa}{4} \ds\int_{\mathbb T^d}
  (\Sigma\star |u|^2)(x)|u(x)|^2\ud x.
\end{align}
As observed above, $\mathscr H$ is a constant of the motion.

Moreover, it is clear that \eqref{Hartree_k} is invariant under multiplications by a phase factor so that $F(u)=\frac{1}{2}\|u\|^2_{L^2}$ is conserved by the dynamics.
The quantities 
\begin{equation*}
  G_{j}(u)=\frac12 \int_{\T^d}\left(\frac{1}{i}\partial_{x_j}u\right)\overline{u}\ud x
\end{equation*}
are constants of the motion too, that correspond to the invariance under translations.
Indeed, a direct verification leads to
\[
\ds\frac{\ud}{\ud t}G_j(u)(t)
=
\frac{\kappa\gamma^2}2 \int_{\T^d} \int_{\T^d} \partial_{x_j}\Sigma (x-y)\star |u|^2(t,y) |u|^2(t,x) \ud y \ud x=0
.
\] 

Finally, we shall endow the Banach space
$H^1(\mathbb T^d;\mathbb R)\times H^1(\mathbb T^d;\mathbb R)$ with the inner 
product
\[
\left\langle \begin{pmatrix}
q
\\
p\end{pmatrix} \Big|  \begin{pmatrix}
q'
\\
p'\end{pmatrix}\right\rangle=\ds\int_{\mathbb T^d} \big( pp'+qq')\ud x.\]
that can be also interpreted as an inner product for complex-valued functions:
\begin{equation}
  \label{eq:innerprodC_Hartree}
  \langle u | u'\rangle =\mathrm{Re} \ds\int_{\mathbb T^d} u\overline{u'}\ud x.
\end{equation}

\subsection{Linearized problem and spectral stability}

Let us expand the solution of 
\eqref{Hartree_k} around $u_\omega$
as 
$u(t,x)=u_\omega(t,x)(1+w(t,x))$.
The linearized 
equation for the fluctuation reads
\begin{equation}\label{lin_hartree}
i\partial_ t w +\ds\frac12\Delta _x w+ik\cdot\nabla_x w= -2\gamma^2\kappa
(\Sigma\star \mathrm{Re} (w)).
\end{equation}
We split $w=q+ip$, $q=\mathrm{Re}(w)$, $p=\mathrm{Im}(w)$ so that
\eqref{lin_hartree} recasts as 
\begin{equation}\label{lin_hartree_qp}
\partial_t\begin{pmatrix}
q\\p\end{pmatrix}=\mathbb L_k\begin{pmatrix}
q\\p\end{pmatrix}
\end{equation}
with the linear operator
\begin{equation}\label{lin_hartree_qp_op}
\mathbb L_k:
\begin{pmatrix}
q\\p\end{pmatrix}\longmapsto 
\begin{pmatrix}-k\cdot\nabla_x q 
-\ds\frac12\Delta _x p
\\
\ds\frac12\Delta _xq + 2\gamma^2\kappa  
\Sigma\star q-  k\cdot\nabla_x p
\end{pmatrix}.
\end{equation}
{\color{black}
From now on, while $(q,p)$ has been  introduced as a pair of real-valued functions, we consider
$\mathbb L_k$ as acting on the $\mathbb C$-vector space of complex-valued functions $L^{2}(\T^d;\C)\times L^{2}(\T^d;\C)$, 
and we study its spectrum.}

\begin{theo}[Spectral stability for the Hartree equation]\label{spectralLkH} Let $k\in \Z^d$ and $\omega\in \mathbb R$ such that the dispersion relation \eqref{eq:dispersion} is satisfied. Suppose \eqref{small} holds. Then the spectrum of $\mathbb L_k$, the linearization of \eqref{Hartree_k} around the plane wave $u_\omega(t,x)=e^{i\omega t}\mathbf 1(x)$, in $L^{2}(\T^d;\C)\times L^{2}(\T^d;\C)$ is contained in $i\R$. Consequently, this wave is spectrally stable in $L^{2}(\T^d)$.
\end{theo}

\begin{Proof}
To prove Theorem \ref{spectralLkH}, we expand $q$, $p$ and $\sigma_1$ by means of their Fourier series
\[\begin{array}{ll}
q(t,x)=\ds\sum_{m\in \mathbb Z^d} Q_{m}(t)e^{im\cdot x},\qquad &
Q_{m}(t)=\ds\frac{1}{(2\uppi)^d}\ds\int_{\mathbb T^d} q(t,x)e^{-im\cdot x} \ud x,
\\[.3cm]
p(t,x)=\ds\sum_{m\in \mathbb Z^d} P_{m}(t)e^{im\cdot x},\qquad &
P_{m}(t)=\ds\frac{1}{(2\uppi)^d}\ds\int_{\mathbb T^d} p(t,x)e^{-im\cdot x} \ud x,
\\
\sigma_1(x)=\ds\sum_{m\in \mathbb Z^d} \sigma_{1,m}e^{im\cdot x},\qquad &
\sigma_{1,m}(t)=\ds\frac{1}{(2\uppi)^d}\ds\int_{\mathbb T^d} \sigma_1(x)e^{-im\cdot x} \ud x.
\end{array}\]
Note that $\sigma_1$ being real and radially symmetric, we have
\begin{equation}\label{condsig1}
\overline{\sigma_{1,m}}=\sigma_{1,m}=\sigma_{1,-m}
\end{equation}
and, by definition, $\avex{\sigma_1}=(2\uppi)^d \sigma_{1,0}$.
As a consequence, we obtain
\begin{align}\label{lin_hartree_qp_fourier}
  \mathbb L_k
    \begin{pmatrix}
      q\\p
    \end{pmatrix}&=
    \begin{pmatrix}
      \sum_{m\in \Z^d} \left(\ds\frac{m^2}2 P_m-i k\cdot m Q_m\right)e^{im\cdot x}\\
      \sum_{m\in \Z^d} \left(-\ds\frac{m^2}{2} Q_m-ik\cdot m   P_m
      +2(2\uppi)^{2d}\gamma^2\kappa  
      |\sigma_{1,m}|^2  Q_m\right)e^{im\cdot x}
    \end{pmatrix}\nonumber\\
    &= \mathbb L_{k,0}\begin{pmatrix}
      Q_0\\P_0
    \end{pmatrix}+\sum_{m\in \Z^d\smallsetminus\{0\}}\mathbb L_{k,m}\begin{pmatrix}
      Q_m\\P_m
    \end{pmatrix}
    {\color{black} e^{im\cdot x}}
\end{align}
with 
\begin{equation}\label{lin_hartree_matrix}
  \mathbb L_{k,0}=\begin{pmatrix}
    0 & 0\\
    2(2\uppi)^{2d}\gamma^2\kappa  
      |\sigma_{1,0}|^2 & 0
  \end{pmatrix}
  \text{ and }
  \mathbb L_{k,m}=\begin{pmatrix}
    -i k\cdot m & \frac{m^2}2\\
    - \frac{m^2}2+2(2\uppi)^{2d}\gamma^2\kappa  
    |\sigma_{1,m}|^2 & -i k\cdot m
  \end{pmatrix}
\end{equation}
for $m\in \Z^d\smallsetminus\{0\}$.

Note that, since the Fourier modes are uncoupled, $\begin{pmatrix}
  q\\ p
\end{pmatrix}$ is a solution to \eqref{lin_hartree_qp} if and only if the Fourier coefficients $\begin{pmatrix}
  Q_m\\ P_m
\end{pmatrix}$
satisfy
\begin{equation*}
  \partial_t \begin{pmatrix}
    Q_m(t)\\ P_m(t)
  \end{pmatrix}= \mathbb L_{k,m} \begin{pmatrix}
    Q_m(t)\\ P_m(t)
  \end{pmatrix}
\end{equation*}
for any $m\in \Z^d$. Similarly, $\lambda \in \C$ is an eigenvalue of the operator $\mathbb L_k$ if and only if there exists at least one Fourier mode $m\in \Z^d$ such that $\lambda$ is an eigenvalue of the matrix $\mathbb L_{k,m}$, \emph{i.e.} there exists $(q_m, p_m)\neq (0,0)$ such that
\begin{equation}\label{eighartree}
\begin{array}{l}
\lambda q_m - \ds\frac{ m^2}{2} p_m +   ik\cdot m q_m= 0,
\\[.3cm]
\lambda p_m+ \ds\frac{ m^2}{2}q_m+   ik\cdot m p_m=2(2\uppi)^{2d}\gamma^2\kappa  
|\sigma_{1,m}|^2q_m.
\end{array}\end{equation}

A straightforward computation gives that $\lambda_0=0$ is the unique eigenvalue of the matrix $\mathbb L_{k,0}$ with eigenvector $(0,1)$. This means that $\mathrm {Ker}(\mathbb L_k)$ contains at least the vector subspace
spanned by the constant function $x\in \mathbb T^d\mapsto \begin{pmatrix}0\\1\end{pmatrix}$, which corresponds to the constant solution $u(t,x)=i$ of \eqref{lin_hartree}.

Next, if $m\in \Z^d\smallsetminus\{0\}$, $\lambda_m$ is an eigenvalue of $\mathbb L_{k,m}$ if it is a solution to 
\begin{equation*}
  (\lambda+ik\cdot m)^2-\frac{m^2}{2}\left(- \frac{m^2}2+2(2\uppi)^{2d}\gamma^2\kappa  
  |\sigma_{1,m}|^2\right)=0.
\end{equation*}
This is a second order polynomial equation for $\lambda$
and the roots are given by
\[
\lambda_{m,\pm}=-ik\cdot m\pm \ds\frac{|m|}2\sqrt{-m^2+4\gamma^2\kappa (2\uppi)^{2d}|\sigma_{1,m}|^2  
}.
\]
If the smallness condition \eqref{small} holds, the argument of the square root is negative for any $m\in \mathbb Z^d\smallsetminus\{0\}$, and thus 
the roots $\lambda$ are all purely imaginary (and we note that $\overline{\lambda_{-m,\pm}}=\lambda_{m,\mp}$).
More precisely, we have the following statement.

\begin{lemma}[Spectral stability for the Hartree equation]
\label{SpecStabH}
Let $k,m\in \Z^d$ and $\mathbb L_{k,m}$ defined as in \eqref{lin_hartree_matrix}. Then
\begin{enumerate}
  \item $\lambda_0=0$ is the unique eigenvalue of $\mathbb L_{k,0}$ and $\mathrm{Ker}(\mathbb L_{k,0})=\mathrm{span}\left\{\begin{pmatrix}0\\1 \end{pmatrix}
    \right\}$; 
  \item for any $m\in \Z^d\smallsetminus\{0\}$, the eigenvalue of $\mathbb L_{k,m}$ are 
  \[
\lambda_{m,\pm}=-ik\cdot m\pm \ds\frac{|m|}2\sqrt{-m^2+4\gamma^2\kappa (2\uppi)^{2d}|\sigma_{1,m}|^2  
}.
\]
  \begin{enumerate}
    \item if $4\gamma^2\kappa (2\uppi)^{2d}  
    \frac{|\sigma_{1,m}|^2}{m^2}\le 1$, then $\lambda_{m,\pm}\in i\R$;
    \item if $4\gamma^2\kappa (2\uppi)^{2d}  
    \frac{|\sigma_{1,m}|^2}{m^2}> 1$, then $\lambda_{m,\pm}\in \C\smallsetminus i\R$. Moreover, $\mathrm{Re}(\lambda_{m,+})>0$.
  \end{enumerate} 
\end{enumerate}
\end{lemma}

Now, \eqref{small} implies  $4\gamma^2\kappa (2\uppi)^{2d}  
\frac{|\sigma_{1,m}|^2}{m^2}< 1$ for all $m\in \Z^d\smallsetminus\{0\}$, so that 
 $\sigma(\mathbb L_k)\subset i\R$ and $u_\omega(t,x)=e^{i\omega t}\mathbf 1(x)$ is spectrally stable. Conversely, if $\sigma_1$, $\sigma_2$ and $\gamma$ are such that there exists $m_*\in \Z^d\smallsetminus\{0\}$ verifying  $4\gamma^2\kappa (2\uppi)^{2d}  
\frac{|\sigma_{1,m_*}|^2}{m_*^2}> 1$, then the plane wave $u_\omega$
 is spectrally unstable for any $k\in \Z^d$ and $\omega\in \mathbb R$ that satisfy the dispersion relation \eqref{eq:dispersion}. This proves {\color{black}  Theorem~\ref{spectralLkH}.}
\end{Proof}

We observe that this result is consistent with the linear stability analysis of \eqref{nls}, see 
\cite[Theorem~1]{Keller}, when replacing formally $\Sigma$ by the delta-Dirac.
The analogy should be considered with caution, though, 
since the functional 
difficulties are substantially different:
here $u\mapsto -\frac12\Delta_{\mathbb T^d}u- 2\gamma^2\kappa \Sigma\star \mathrm{Re}(u)$ is a compact perturbation  of $-\frac12 \Delta_{\mathbb T^d}$, 
 which has a compact resolvent hence a spectral decomposition.

It is important to remark that the analysis of  eigenproblems for $\mathbb L_{k}$ has consequences on the behavior of solutions to \eqref{lin_hartree_qp} of the particular form 
\begin{equation*}
    Q(t,x)=e^{\lambda t}q(x),\qquad P(t,x)=e^{\lambda t}p(x).
\end{equation*}
We warn the reader that spectral stability excludes the \emph{exponential} growth
of the solutions of the linearized problem
when the smallness condition \eqref{small} holds, but a slower growth is still possible.
This can be seen by direct inspection for the mode $m=0$: we have $\partial_t Q_0=0$, so that $Q_0(t)=Q_0(0)$
 and $\partial_t P_0=2(2\uppi)^{2d}\kappa  
  \avex{\sigma_1}^2Q_0(0)$ which
 shows that the solution can grow linearly in time
 \[
 P_0(t)=P_0(0)+2(2\uppi)^{2d}\gamma^2\kappa  
  \avex{\sigma_1}^2Q_0(0)t.
 \]
 In fact, excluding the mode $m=0$ suffices to guaranty the linearized stability.
 
 \begin{theo}[Linearized stability for the Hartree equation]
Suppose \eqref{small}. 
Let $w$ be the solution of  \eqref{lin_hartree} associated to an initial data $w^{\mathrm{Init}}
\in H^1(\mathbb T^d)$ such that $\int_{\mathbb T^d}w^{\mathrm{Init}}\ud x=0$.
Then, there exists a constant $C>0$ such that $\sup_{t\geq 0}\|w(t,\cdot)\|_{H^1}\leq C$.
\label{Prop:Hlstab}
\end{theo}

\begin{Proof}
Note that if $\int_{\mathbb T^d}w^{\mathrm{Init}}\ud x=0$ then the corresponding Fourier coefficients $Q_0(0)$ and $P_0(0)$ are equal to $0$. As a consequence, $Q_0(t)=P_0(t)=0$ for all $t\ge 0$, so that $\int_{\mathbb T^d}w(t,x)\ud x=0$ for all $t\ge 0$.

The proof follows from energetic consideration.
Indeed,  we observe that, on the one hand,
\[\begin{array}{l}
\ds\frac12\ds\frac{\ud}{\ud t}\ds\int_{\mathbb T^d} |\nabla w|^2\ud x
=-\ds\frac{\gamma^2\kappa  
}{2i} \ds\int_{\mathbb T^d}\Sigma\star (w+\overline w) \Delta (w-\overline w) \ud x,
\end{array}\]
and, on the other hand, 
\[\begin{array}{l}
\ds\frac12\ds\frac{\ud}{\ud t}\ds\int_{\mathbb T^d}\Sigma\star (w+\overline w) (w+\overline w)  \ud x
\\
\qquad
=-\ds\frac{1}{2i} \ds\int_{\mathbb T^d}\Sigma\star (w+\overline w) \Delta (w-\overline w) \ud x
-k\cdot \ds\int_{\mathbb T^d} \nabla(w+\overline w) \Sigma\star (w+\overline w)\ud x,
\end{array}
\]
where we get rid of the last term in the right hand side by assuming $k=0$.
This leads to the following energy conservation property
\[
\ds\frac{\ud}{\ud t}\left\{
\ds\frac12\ds\int_{\mathbb T^d} |\nabla w|^2\ud x
- \ds\frac{\gamma^2\kappa  
}2\ds\int_{\mathbb T^d}\Sigma\star (w+\overline w) (w+\overline w)  \ud x
\right\}=0
\]
which holds for $k=0$. We denote by $E_0$ the energy of the initial data $w^{\mathrm{Init}}$. 
Finally, we can simply estimate
\[\left|
\ds\int_{\mathbb T^d}\Sigma\star (w+\overline w) (w+\overline w)  \ud x\right|
\leq \|\Sigma\star(w+\overline w) \|_{L^2}\|w+\overline w \|_{L^2}
\leq \|\Sigma\|_{L^1}\|w+\overline w \|^2_{L^2}
\leq 4\|\Sigma\|_{L^1}\|w\|^2_{L^2}.
\]
To conclude, we use the Poincar\'e-Wirtinger estimate. Indeed, since we have already remarked that the condition $\int_{\mathbb T^d}w^{\mathrm{Init}}\ud x=0$
implies $\int_{\mathbb T^d}w(t,x)\ud x=0$ for any $t\geq 0$, we can write
\begin{align*}
\left\|w(t,\cdot)\right\|_{L^2}^2&=\Big\|w(t,\cdot)-\frac{1}{(2\uppi)^d}\ds\int_{\mathbb T^d} w(t,y)\ud y\Big\|^2_{L^2}
=(2\uppi)^d\ds\sum_{m\in \Z^d\smallsetminus\{0\}} |c_m(w(t,\cdot))|^2 \\
&\leq \ds(2\uppi)^d\ds\sum_{m\in \Z^d\smallsetminus\{0\}} m^2|c_m(w(t,\cdot))|^2
=  \|\nabla w(t,\cdot)\|^2_{L^2}
\end{align*}
for any $t\ge 0$, where the $c_m(w(t,\cdot))$'s are the Fourier coefficients of the function $x\in \T^d\mapsto w(t,x)$.
Hence, for any solution with zero mean,  we infer, for all $t\ge 0$,
\begin{equation*}
  2E_0= \ds\int_{\mathbb T^d} |\nabla w|^2(t,x)\ud x
  - \ds\gamma^2\kappa  
  \ds\int_{\mathbb T^d}\Sigma\star (w+\overline w) (w+\overline w)(t,x)  \ud x \ge (1-4\gamma^2\kappa\|\Sigma\|_{L^1})\int_{\mathbb T^d} |\nabla w(t,x)|^2\ud x.
\end{equation*}
As a consequence, if \eqref{small} is satisfied, we obtain
\begin{equation*}
  \sup_{t\geq 0}\|w(t,\cdot)\|_{H^1}\leq 2\sqrt{\frac{E_0}{1-4\gamma^2\kappa\|\Sigma\|_{L^1}}}.
\end{equation*}

The stability estimate extends to the situation where $k \neq 0$. Indeed, from the solution $w$ of \eqref{lin_hartree}, we set 
\[v(t,x)=w(t,x+tk).\]
It satisfies $i\partial_ t v +\frac12\Delta_x v=-2\gamma^2\kappa  
\Sigma\star \mathrm{Re}(v)$.
Hence, repeating the previous argument, $\|v(t,\cdot)\|_{H^1}=\|w(t,\cdot)\|_{H^1}$ remains uniformly bounded on $(0,\infty)$.
This step of the proof relies on the Galilean invariance of \eqref{hartree}; it could have been 
used from the beginning, but it  does not apply for the Schr\"odinger-Wave system.
\end{Proof}

\begin{rmk}
The analysis applies mutadis mutandis to any  equation of the form \eqref{Hartree},
with the potential defined by a kernel $\Sigma$  and a strength encoded by the constant $\gamma^2\kappa$.
Then, the stability criterion
is set on the quantity $4\gamma^2\kappa (2\uppi)^d  
\frac{|\widehat\Sigma_m|}{m^2} $
For instance, the elementary solution of 
 $ (a^2-\Delta_x)\Sigma=\delta_{x=0}$ with periodic boundary condition  has its Fourier coefficients given by  
$\widehat \Sigma_m=\frac{1}{(2\uppi)^d(a^2+m^2)}>0$.
Coming back to the physical variable, in the one-dimension case, the function $\Sigma$ reads
\[
\Sigma(x)=\ds\frac{e^{-a|x|}}{2a}+\ds\frac{\cosh(ax)}{a(e^{2a\uppi}-1)}.\]
The linearized stability thus holds provided $4\gamma^2\kappa (2\uppi)^{2d}\frac{1  
}{a^2+1}<1$.
\end{rmk}

\subsection{Orbital stability}\label{sec:orbitalHartree}

In this subsection, we wish to establish the \emph{orbital stability} of the plane wave $u_\omega(t,x)=e^{i\omega t}\mathbf 1(x)$ as a solution to \eqref{Hartree_k} for $k\in \Z^d$ and $\omega\in \mathbb R$ that satisfy the dispersion relation \eqref{eq:dispersion}.
As pointed out before, \eqref{Hartree_k} is invariant under multiplications by a phase factor. This leads to define the corresponding orbit through $u(x)=\mathbf 1(x)$ by
\[
\mathscr O_{\mathbf 1}=\{e^{i\theta},\ \theta \in \mathbb R\}.\]
Intuitively, orbital stability means that the solutions 
of \eqref{Hartree_k}
associated to initial data close enough to  the constant function $x\in \mathbb T^d \mapsto 
1=\mathbf 1(x)$ remain at a close distance to the set $\mathscr O_{\mathbf 1}$. 
Stability analysis then amounts to the construction of a suitable Lyapounov functional
satisfying a coercivity property. This functional should be a constant of the motion and be invariant under the action of the group that generates the orbit $\mathscr O_{\mathbf 1}$. Hence, the construction of such a functional relies on the invariants of the equation. Moreover, the plane wave has to be a critical point on the Lyapounov functional so that the coercivity can be deduced from the properties of its second variation.  
The difficulty here is that, in general, the bilinear symmetric form defining the second variation of the Lyapounov function is not positive on the whole space: according to the strategy designed in \cite{GSS}, see also the review \cite{Tao}, it will be enough to prove the coercivity on an appropriate subspace. 
Here and below, we adopt the framework presented in  \cite{SRN1} (see also \cite{SRN2}).

Inspired by the strategy designed in \cite[Section~8~\&~9]{SRN1}, we introduce, for any $k\in \Z^d$ and 
$\omega\in \mathbb R$ satisfying the dispersion relation \eqref{eq:dispersion}, the set 
\begin{equation*}
  \mathscr S_\omega=\Big\{u\in H^1(\T^d;\C),\ F(u)=F(\mathbf 1)=\ds\frac{(2\uppi)^{d}  
}{2}=(2\uppi)^{d} \ds\frac{k^2/2+ \omega}
{2\gamma^2\kappa \avex{\sigma_1}^2}\Big\};
\end{equation*}
$\mathscr S_\omega$
 is therefore the level set of the solutions of
\eqref{Hartree_k}, 
associated to the 
plane wave $(t,x)\mapsto u_\omega(t,x)=e^{i\omega t}\mathbf 1(x)$.
Next, we introduce the functional 
\begin{equation}
  \label{eq:lyapounov_Hartree}
  \mathscr L_\omega(u)=\mathscr H(u)+ \omega F(u)-\sum_{j=1}^dk_jG_j(u),
\end{equation}
which is conserved by the solutions of \eqref{Hartree_k}.
We have
\[
\begin{array}{lll}
\partial_u\mathscr L_\omega(u)(v)&=&
 \mathrm{Re}\left(\ds\frac12
 \ds\int_{\mathbb T^d} (-\Delta u)\overline v\ud x
 +\ds\frac{k^2}{2}\ds\int_{\mathbb T^d} u\overline v\ud x
 \right.
 \\[.3cm]
 &&\left.
\hspace*{2cm} -\gamma^2 \kappa
\ds\iint_{\mathbb T^d\times\mathbb T^d}
 \Sigma(x-y)|u(y)|^2u(x) \overline{v(x)}\ud y\ud x
 {+} \omega \ds\int_{\mathbb T^d} u \overline v\ud x\right).
\end{array}\]
As a matter of fact,  we observe that
\[\partial_u\mathscr L_\omega(
\mathbf 1)=0\]
owing to the dispersion relation.
Next, we get
\[
\begin{array}{lll}
\partial^2_u\mathscr L_\omega(u)(v,w)&=&
 \mathrm{Re}\left(
\ds\frac12 \ds\int_{\mathbb T^d} (-\Delta +k^2)w \overline v\ud x\right.
 \\[.3cm]
&& \left.-2\gamma^2\kappa
 \ds\iint_{\mathbb T^d\times\mathbb T^d}
 \Sigma(x-y)\mathrm{Re}\big(\overline {u(y)} w(y)\big)
 u(x)\overline{v(x)} \ud y\ud x
 \right.
 \\[.3cm]
 &&\left.
 -\gamma^2 \kappa
\ds\iint_{\mathbb T^d\times\mathbb T^d}
 \Sigma(x-y)|u(y)|^2 w(x) \overline{v(x)} \ud y\ud x
 +\omega \ds\int_{\mathbb T^d} w \overline v\ud x\right).
\end{array}\]
Still by using the dispersion relation, we obtain
\[
\partial^2_u\mathscr L_\omega(
\mathbf 1)(v,w)
=\mathrm{Re}\left(\ds\int_{\mathbb T^d}
\underbrace{\left(-\ds\frac{\Delta w}{2}-2\gamma^2\kappa 
\Sigma\star  \mathrm{Re}(w)
\right)}_{:=\mathbb Sw}
\overline{v(x)}
\ud x\right)=\langle \mathbb S w | v\rangle.\]
$\mathbb S:H^2(\T^d)\subset L^2(\T^d) \to L^2(\T^d)$ is an unbounded linear operator and its spectral properties will play an important role for the orbital stability of $u_\omega$. Note that the operator $\mathbb S$ is the linearized operator \eqref{lin_hartree_qp_op}, up to the advection term $k\cdot \nabla$.
The main result of this subsection is the following.
\begin{theo}[Orbital stability for the Hartree  equation]
  \label{prop:orbital_Hartree_k} Let $k\in \Z^d$ and $\omega\in \mathbb R$ such that the dispersion relation \eqref{eq:dispersion} is satisfied. Suppose \eqref{small} holds. Then the plane wave $u_\omega(t,x)=e^{i\omega t}\mathbf 1(x)$ is orbitally stable, \emph{i.e.}
  \begin{equation}
    \label{eq:orbital_Hartree_k}
    \forall \varepsilon>0, \ \exists \delta>0,\  \forall 
    v^{\mathrm{Init}}\in H^1(\T^d;\C),\
     \|v^{\mathrm{Init}}-\mathbf 1\|_{H^1}<\delta
      \Rightarrow \sup_{t\ge 0} 
      \mathrm{dist}(v(t),\mathscr O_{\mathbf 1})<\varepsilon
  \end{equation}
   where {\color{black} $\mathrm{dist}(v,\mathscr O_{\mathbf 1})=\inf_{{\theta} \in [0,2\uppi[}\|v-e^{i {\theta} }\mathbf 1\|_{H^1}$}
 and $(t,x)\mapsto v(t,x)\in C^0([0,\infty);H^1(\T^d))$ stands for the solution 
  of \eqref{Hartree_k} with Cauchy data $v^{\mathrm{Init}}$.
\end{theo}

{\color{black} The full proof of Theorem \ref{prop:orbital_Hartree_k} 
will be obtained from a series of intermediate steps, that we detail now.}
The key ingredient to prove Theorem \ref{prop:orbital_Hartree_k} is the following coercivity estimate on the Lyapounov functional.

\begin{lemma}\label{lem:coer_Hartree_k} Let $k\in \Z^d$ and $\omega\in \mathbb R$ such that the dispersion relation \eqref{eq:dispersion} is satisfied.
  Suppose that there exist $\eta>0$ and $c>0$ such that 
  \begin{equation}
    \label{eq:coer_Hartree}
    \forall w\in \mathscr S_\omega, d(w,\mathscr O_{\mathbf 1})<\eta\Rightarrow \mathscr L_\omega(w)-\mathscr L_\omega(\mathbf 1)\ge c \
     \mathrm{dist} (w,\mathscr O_{\mathbf 1})^2.
  \end{equation}
  Then  the plane wave $u_\omega(t,x)=e^{i\omega t}\mathbf 1(x)$ is orbitally stable.
\end{lemma}

%{\color{black} Let us postpone for a while the justification of this claim, and  explain how it can be used to prove Theorem \ref{prop:orbital_Hartree_k}.}

%\begin{ProofOf}{Theorem \ref{prop:orbital_Hartree_k}} 
\begin{Proof} Assume that  {\color{black} \eqref{eq:coer_Hartree}} %that Lemma \ref{lem:coer_Hartree_k} 
holds and 
  suppose, by contradiction, that $u_\omega$ is not orbitally stable. Hence, there exists $0<\varepsilon_0< \frac{2}{3}\eta$ such that 
  \begin{equation*}
    \forall n\in {\N\smallsetminus\{0\},\ } \exists u_n^{\mathrm{Init}}\in H^{1}(\T^d), \|u_n^{\mathrm{Init}}-\mathbf 1\|_{H^1}<\frac{1}{n} \text{ and } \exists t_n\in [0,+\infty[,  \mathrm{dist}(u_n(t_n),\mathscr O_{\mathbf 1})=\varepsilon_0,
  \end{equation*}
 $(t,x)\mapsto u_n(t,x)\in C^0([0,\infty);H^1(\T^d))$ being the solution 
  of \eqref{Hartree_k} with Cauchy data $u_n^{\mathrm{Init}}$.
  To apply the coercivity estimate of Lemma \ref{lem:coer_Hartree_k}, we define $z_n=\left(\frac{F(\mathbf 1)}{F(u_n(t_n))}\right)^{1/2}u_n(t_n)$. It is clear that $z_n\in{\mathscr  S_\omega}$ since $F(z_n)=F(\mathbf 1)$. Moreover, $\big(u_n(t_n)\big)_{n\in\mathbb N\smallsetminus\{0\}}$ is a bounded sequence in $H^1(\T^d)$ and $\lim_{n\to +\infty} F(u_n(t_n))=F(\mathbf 1)$. Indeed, on the one hand, there exists $\gamma\in [0,2\uppi[$ such that
  \begin{equation*}
    \|u_n(t_n)\|_{H^1}\le \|u_n(t_n)-e^{i\theta}\mathbf 1\|_{H^1}+\|e^{i\theta}\mathbf 1\|_{H^1}\le 2d(u_n(t_n),\mathscr O_{\mathbf 1})+\|e^{i\theta}\mathbf 1\|_{H^1}=2\varepsilon_0 + \|\mathbf 1\|_{H^1}
  \end{equation*}
  and, on the other hand,  
  \begin{equation*}
    |F(u_n(t_n))-F(\mathbf 1)|=\frac12|\|u_n(t_n)\|_{L^2}^2-\|\mathbf 1\|^2_{L^2}|\le \|u_n(t_n)-\mathbf 1\|_{L^2}(\varepsilon_0 + \|\mathbf 1\|_{H^1})<\frac{1}{n}(\varepsilon_0 + \|\mathbf 1\|_{H^1}).
  \end{equation*}
  As a consequence, $\lim_{n\to+\infty}\|z_n-u_n(t_n)\|_{H^1}=0$. This implies for $n\in \N$ large enough,
  \begin{equation*}
    \frac{\varepsilon_0}{2}\le d(z_n,\mathscr O_{\mathbf 1})\le \frac{3\varepsilon_0}{2}<\eta.
  \end{equation*}
  Hence, thanks to Lemma \ref{lem:coer_Hartree_k}, we obtain 
  \begin{align*}
    \mathscr{L}_\omega( u_n^{\mathrm{Init}})-\mathscr{L}_\omega(\mathbf 1)&=\mathscr{L}_\omega(u_n(t_n))-\mathscr{L}_\omega(\mathbf 1)=\mathscr{L}_\omega(u_n(t_n))-\mathscr{L}_\omega(z_n)+\mathscr{L}_\omega(z_n)-\mathscr{L}_\omega(\mathbf 1)\\
    &\ge \mathscr{L}_\omega(u_n(t_n))-\mathscr{L}_\omega(z_n)+cd(z_n,{\mathscr O_{\mathbf 1}})^2\ge \mathscr{L}_\omega(u_n(t_n))-\mathscr{L}_\omega(z_n)+\frac{c}{4}\varepsilon_0^2.
  \end{align*}
  Finally, using the fact that $\partial_u \mathscr L_\omega(\mathbf 1)=0$ and $\partial_u^2 \mathscr L_\omega(\mathbf 1)(w,w)\le C \|w\|_{H^1}^2$, we deduce that 
  \begin{align*}
    &\lim_{n\to+\infty}(\mathscr{L}_\omega(u_n^{\mathrm{Init}})-\mathscr{L}_\omega(\mathbf 1))=0,\\
    &\lim_{n\to+\infty}(\mathscr{L}_\omega(u_n(t_n))-\mathscr{L}_\omega(z_n))=0.
  \end{align*}
We are thus led to a contradiction.
%\end{ProofOf}
\end{Proof}

Since $\partial_u \mathscr L_\omega(\mathbf 1)=0$, the coercivity estimate \eqref{eq:coer_Hartree} can be obtained from a similar estimate on the bilinear form {\color{black}  $w\in H^1\mapsto \partial_u^2 \mathscr L_\omega(\mathbf 1)(w,w)$.} %for any $w\in H^1$}.
 As pointed out before, the difficulty lies in the fact that, in general, this bilinear form is not positive on the whole space $H^1$. The following lemma states that it is enough to have a coercivity estimate on $\partial_u^2 \mathscr L_\omega(\mathbf 1)(w,w)$ for any $w\in T_{\mathbf 1}\mathscr S_{\omega}\cap (T_{\mathbf 1}\mathscr O_{\mathbf 1})^{\perp}$. 
Recall that the tangent set to $\mathscr S_\omega$ is given by 
\[T_{ \mathbf 1}\mathscr S_\omega=\{u\in H^1(\T^d;\C), \partial_u F(\mathbf 1)(u)=0\}=
\left\{(q,p)\in H^1(\T^d,\R)\times H^1(\T^d,\R),\ \Big\langle 
\begin{pmatrix} q \\ p\end{pmatrix}\Big| \begin{pmatrix} \mathbf 1 \\ {0}\end{pmatrix}\Big\rangle=0\right\}.
\]
This set is the orthogonal to $\mathbf 1$ with respect to the inner product defined in \eqref{eq:innerprodC_Hartree}. The tangent set to $\mathscr O_{\mathbf 1}$ (which is the orbit generated by the phase multiplication) is
\begin{equation*}
  T_{ \mathbf 1}\mathscr O_{\mathbf 1}=\mathrm{span}_{\R}\{i\mathbf 1\}
\end{equation*}
so that 
\[(T_{ \mathbf 1}\mathscr O_\mathbf 1)^\perp=\{u\in H^{1}(\T^d,\C), \langle u,i\mathbf 1\rangle=0 \}=
\left\{(q,p): \mathbb T^d\to  \mathbb R,\ \Big\langle 
\begin{pmatrix} q \\ p\end{pmatrix}\Big|
\begin{pmatrix} 0 \\ \mathbf  1\end{pmatrix}\Big\rangle
=0 \right\}.\]

\begin{lemma}\label{lem:coer_Hartree_hess} Let $k\in \Z^d$ and $\omega\in \mathbb R$ such that the dispersion relation \eqref{eq:dispersion} is satisfied. Suppose that there exists $\tilde c>0$ 
  \begin{equation}
    \label{eq:coer_Hartree_hess}
    \partial_u^2 \mathscr L_\omega(\mathbf 1)(u,u)\ge \tilde c\|u\|^2_{H^1}
  \end{equation}
  for any $u\in T_{\mathbf 1}\mathscr S_{\mathbf 1}\cap (T_{\mathbf 1}\mathscr O_{\mathbf 1})^{\perp}$.
Then there exist $\eta>0$ and $c>0$ such that \eqref{eq:coer_Hartree} is satisfied.
\end{lemma}

\begin{Proof}
  Let $w\in \mathscr{S}_{\omega}$ such that $  \mathrm{dist}(w,\mathscr O_{\mathbf 1})<\eta$ with $\eta>0$ small enough. By means of an implicit function theorem argument (see \cite[Section 9, Lemma 8]{SRN1}), we obtain that there exists $\theta \in [0,2\uppi[$ and $v\in (T_{\mathbf 1}\mathscr O_{\mathbf 1})^{\perp}$ such that 
  \begin{equation*}
    e^{i{\theta}}w=\mathbf 1+v,\qquad  {\mathrm{dist}}(w,\mathscr O_{\mathbf 1})\le \|v\|_{H^1}\le C  {\mathrm{dist}}(w,\mathscr O_{\mathbf 1})
  \end{equation*}
  for some positive constant $C$. 

  Next, we use the fact that $H^1(\T^d)=T_{   \mathbf 1}\mathscr S_\omega\oplus \mathrm{span}_{\R}\{\mathbf{1}\}$ to write $v=v_1+v_2$ with $v_1\in T_{   \mathbf 1}\mathscr S_\omega\cap (T_{\mathbf 1}\mathscr O_{\mathbf 1})^{\perp}$ and $v_2\in \mathrm{span}_{\R}\{\mathbf{1}\}\cap (T_{\mathbf 1}\mathscr O_{\mathbf 1})^{\perp}$. Since $v=e^{i {\theta}}w-\mathbf 1$ and $F(w)=F(\mathbf 1)$, we obtain 
  \begin{align*}
    0=F(e^{i{\theta}}w)-F(\mathbf 1)=\frac12\int_{\T^d}|v|^2\ud x+\mathrm{Re}\int_{\T^d}(v_1+v_2)\mathbf1\ud x= \frac12\int_{\T^d}|v|^2\ud x+\mathrm{Re}\int_{\T^d}v_2\mathbf1\ud x.
  \end{align*}
  Since $v_2\in \mathrm{span}_{\R}\{\mathbf{1}\}$, it follows that
  \begin{equation*}
    \|v_2\|_{H^1}\le \frac{\|v\|_{H^1}^2}{{\color{black} 2}\|\mathbf 1\|_{L^2}}.
  \end{equation*}
  This implies $$\|v_1\|_{H^1}=\|v-v_2\|_{H^1}\ge \|v\|_{H^1}-{\color{black}\frac{\|v\|_{H^1}^2}{2\|\mathbf 1\|_{L^2}}}\ge\frac12 \|v\|_{H^1}$$
  provided $\|v\|_{H^1}\le \|\mathbf 1\|_{L^2}$.
  As a consequence, if $\|v\|_{H^1}$ is small enough, using that $\partial_u^2 \mathscr L_\omega(\mathbf 1)(w,z)\le C \|w\|_{H^1}\|z\|_{H^1}$, we obtain
  \begin{align*}
    &\partial_u^2\mathscr L_\omega(\mathbf 1)(v_1,v_2)\le C \|v\|^3_{H^1},\\
    &\partial_u^2\mathscr L_\omega(\mathbf 1)(v_2,v_2)\le C \|v\|^4_{H^1}.
  \end{align*}
  This leads to 
  \begin{equation*}
    \partial_u^2\mathscr L_\omega(\mathbf 1)(v,v)=\partial_u^2\mathscr L_\omega(\mathbf 1)(v_1,v_1)+o(\|v\|^2_{H^1}).
  \end{equation*}
  Finally, let $w\in \mathscr{S}_{\omega}$ be such that $d(w,\mathscr O_{\mathbf 1})<\eta$. We have
  \begin{align*}
    &\mathscr L_\omega(w)-\mathscr L_\omega(\mathbf 1)=\mathscr L_\omega(e^{i{\theta}}w)-\mathscr L_\omega(\mathbf 1)=\frac12 \partial_u^2\mathscr L_\omega(\mathbf 1)(v,v)+o(\|v\|^2_{H^1})\\
    &=\frac12 \partial_u^2\mathscr L_\omega(\mathbf 1)(v_1,v_1)+o(\|v\|^2_{H^1})\ge {\color{black}\ds\frac{\tilde c}{2}}\|v_1\|^2_{H^1}+o(\|v\|^2_{H^1})\ge {\color{black} \frac{\tilde c}{4}}\|v\|^2_{H^1}+o(\|v\|^2_{H^1})\\
    &\ge {\color{black} \ds\frac{\tilde c}{8}}  {\mathrm{dist}}(w,\mathscr O_{\mathbf 1})^2
  \end{align*}
  where we use $\partial_u\mathscr L_\omega(\mathbf 1)=0$ and $v_1\in T_{   \mathbf 1}\mathscr S_\omega\cap (T_{\mathbf 1}\mathscr O_{\mathbf 1})^{\perp}$.
\end{Proof}

At the end of the day, to prove the orbital stability of the plane wave $u_\omega(t,x)=e^{i\omega t}\mathbf 1(x)$ it is enough to prove \eqref{eq:coer_Hartree_hess} for any $u\in T_{\mathbf 1}\mathscr S_{\mathbf 1}\cap (T_{\mathbf 1}\mathscr O_{\mathbf 1})^{\perp}$. This can be done by studying the spectral properties of the operator $\mathbb S$. However, in the simpler case of the Hartree equation, the coercivity of $\partial_u^2\mathscr L_\omega(\mathbf 1)$ on $T_{\mathbf 1}\mathscr S_{\mathbf 1}\cap (T_{\mathbf 1}\mathscr O_{\mathbf 1})^{\perp}$ can be also obtained directly from the expression
\begin{equation}\label{eq:Hartree_hess}
  \partial^2_u\mathscr L_\omega( \mathbf 1)(u,u)
=\mathrm{Re}\left(\ds\int_{\mathbb T^d}
\left(-\ds\frac{\Delta u}{2}-2\gamma^2\kappa 
\Sigma\star  \mathrm{Re}(u)
\right)
\overline{u(x)}
\ud x\right)
=\langle \mathbb S u|u\rangle.
\end{equation}
Let $u\in T_{\mathbf 1}\mathscr S_{\mathbf 1}\cap (T_{\mathbf 1}\mathscr O_{\mathbf 1})^{\perp}$ and write $u=q+ip$. This leads to 
\begin{equation*}
  \partial^2_u\mathscr L_\omega( \mathbf 1)(u,u)=\frac12\int_{\T^d}|\nabla q|^2\ud x-2\gamma^2\kappa \int_{\T^d}(\Sigma\star q )q \ud x+ \frac12\int_{\T^d}|\nabla p|^2\ud x.
\end{equation*}
Moreover, since  $u\in T_{\mathbf 1}\mathscr S_{\mathbf 1}\cap (T_{\mathbf 1}\mathscr O_{\mathbf 1})^{\perp}$, we have
\begin{equation*}
  \int_{\T^d}q\ud x=0\text{ and } \int_{\T^d}p\ud x=0.
\end{equation*}
As a consequence, thanks to the Poincaré-Wirtinger inequality, we deduce
\begin{equation}\label{proof:coerc1_hartree}
  \partial^2_u\mathscr L_\omega( \mathbf 1)(u,u)\ge \frac12\int_{\T^d}|\nabla q|^2\ud x-2\gamma^2\kappa \int_{\T^d}(\Sigma\star q )q \ud x+ \frac14\|p\|^2_{H^1}.
\end{equation}
Next, we expand $q$ and $\Sigma$ in Fourier series, \emph{i.e.} $$q(x)=\ds\sum_{m\in \mathbb Z^d} q_{m}e^{im\cdot x} \text{ and } \Sigma(x)=\ds\sum_{m\in \mathbb Z^d} \Sigma_{m}e^{im\cdot x}.$$ Note that, if $\Sigma=\sigma_1\star\sigma_1$, then $\Sigma_m=(2\uppi)^d\sigma_{1,m}^2$. Moreover, $\int_{\T^d}q\ud x=0$ implies $q_0=0$.
Hence, 
\begin{align}\label{proof:coerc2_hartree}
  \frac12\int_{\T^d}|\nabla q|^2\ud x-2\gamma^2\kappa \int_{\T^d}(\Sigma\star q )q \ud x=(2\uppi)^d\sum_{m\in \Z^d\smallsetminus \{0\}}\left(\frac{m^2}{2}-2\gamma^2\kappa(2\uppi)^d\Sigma_m\right)q_m^2\nonumber\\
  =(2\uppi)^d\sum_{m\in \Z^d\smallsetminus \{0\}}\left(1-4\gamma^2\kappa(2\uppi)^d\frac{\Sigma_m}{m^2}\right)\frac{m^2}{2}q_m^2. 
\end{align}
As a consequence, we obtain the following statement.

\begin{proposition}
  \label{prop:coer_Hartree_hess_small} Let $k\in \Z^d$ and $\omega\in\mathbb R$ such that the dispersion relation \eqref{eq:dispersion} is satisfied. Suppose that there exists $\delta\in (0,1)$ such that
  \begin{equation}
    \label{cond_small_Hartree}
    4\gamma^2\kappa(2\uppi)^{2d}\frac{\sigma^2_{1,m}}{m^2}\le \delta
  \end{equation}
  for all $m\in \Z^d\smallsetminus\{0\}$.
  Then, there exists $\tilde c>0$ such that   
  \begin{equation}
    \label{eq:coer_Hartree_hess_bis}
    \partial_u^2 \mathscr L_\omega(\mathbf 1)(u,u)\ge \tilde c\|u\|^2_{H^1}
  \end{equation}
  for any $u\in T_{\mathbf 1}\mathscr S_{\mathbf 1}\cap (T_{\mathbf 1}\mathscr O_{\mathbf 1})^{\perp}$.
\end{proposition}

\begin{Proof}
  If \eqref{cond_small_Hartree} holds, then \eqref{proof:coerc1_hartree}-\eqref{proof:coerc2_hartree} lead to 
  \begin{align*}
    \partial^2_u\mathscr L_\omega( \mathbf 1)(u,u)\ge \frac{1-\delta}{2} (2\uppi)^d\sum_{m\in \Z^d\smallsetminus \{0\}}{m^2}q_m^2+ \frac14\|p\|_{H^1}=\frac{1-\delta}{2}\|\nabla q\|_{L^2}^2+\frac14\|p\|^2_{H^1}\ge \frac{1-\delta}{4}\|u\|_{H^1}^2.
  \end{align*}
  where in  the last inequality we used the Poincaré-Wirtinger inequality together with the fact that $\int_{\T^d}q\ud x=0$.
\end{Proof}

\begin{rmk}
  By decomposing the linear operator $\mathbb S$ into real and imaginary part and by using Fourier series, one can study its spectrum. In particular, $\mathbb S$ has exactly one negative eigenvalue $\lambda_-=-2\gamma^2\kappa\avex{\Sigma}$ with eigenspace $\mathrm{span}_\R\{\mathbf 1\}$. Moreover, $\mathrm{Ker}(\mathbb S)=\mathrm{span}_\R\{i
   \mathbf 1\}$. Finally, if \eqref{cond_small_Hartree} is satisifed, then $\inf (\sigma (\mathbb S)\cap (0,\infty))\ge\frac{1-\delta}{2}$. Then, by applying the same arguments as in \cite[Section 6]{SRN2}, we can recover the coercivity of $\partial_u^2\mathscr L_\omega(\mathbf 1)$ on $T_{\mathbf 1}\mathscr S_{\mathbf 1}\cap (T_{\mathbf 1}\mathscr O_{\mathbf 1})^{\perp}$.
\end{rmk}

Finally, Proposition \ref{prop:coer_Hartree_hess_small} together with Lemma \ref{lem:coer_Hartree_hess} and Lemma \ref{lem:coer_Hartree_k}, gives Theorem \ref{prop:orbital_Hartree_k} and the orbital stability of the plane wave $u_\omega$.

\section{Stability analysis of the Schr\"odinger-Wave system: the case $k=0$}
\label{S:SW0}

Like in the case of the Hartree system, to study the stability of the plane wave solutions of the Schrödinger-Wave system \eqref{Schro-s}-\eqref{Schro-p}, it is useful to write its solutions in the form
\begin{equation*}
  U(t,x)=e^{i k\cdot x}u(t,x)
\end{equation*}
with $(t,x,z)\mapsto (u(t,x), \Psi(t,x,z))$ solution to 
\begin{equation}
  \label{sw_k}
  \begin{aligned}
    &i\partial_t u+\ds
\frac12\Delta_x u -\ds\frac{k^2}{2} u + ik\cdot \nabla_x u  =\left(\gamma\sigma_1\star\ds\int _{\mathbb R^n}\sigma_2\Psi\ud z\right)
u,\\
&\frac{1}{c^2}\partial^2_{tt}\Psi-\Delta_z \Psi = -\gamma \sigma_2\sigma_1\star |u|^2.
  \end{aligned}
\end{equation}
If $k\in \Z^d$ and $\omega\in\mathbb R$ satisfy the dispersion relation \eqref{eq:dispersion}, 
\begin{equation*}
  u_\omega(t,x)=e^{i\omega t}\mathbf 1(x),\quad \Psi_*(t,x,z)=-\gamma \Gamma(z)\avex{\sigma_1}, \quad \Pi_* (t,x,z)={\color{black}-\frac{1}{2c^2}}\partial_t\Psi_*(t,x,z)=0
\end{equation*}
with $\Gamma$ the solution of $-\Delta_z \Gamma=\sigma_2$ (see \eqref{formGamma}),
is a solution to \eqref{sw_k} with initial condition 
\begin{equation*}
  u_\omega(0,{\color{black} x})=\mathbf 1(x), \quad \Psi_*(0,x,z)=-\gamma \Gamma(z)\avex{\sigma_1}, \quad \Pi_* (0,x,z)=0.
\end{equation*}
For the time being, we stick to the framework identified for the study of the asymptotic Hartree equation.
Problem \eqref{sw_k} has a natural Hamiltonian symplectic structure  when considered on the \emph{real} Banach space  $H^1(\T^d)\times H^1(\T^d)\times L^2(\T^d;\overbigdot {H}^1(\R^n))\times L^2(\T^d\times \R^n  )$.
Indeed, if we write $u=q+ip$, with $p,q$ real-valued, we obtain 
\[\partial_t 
\begin{pmatrix}
q\\p\\\Psi\\ \Pi
\end{pmatrix}
=\begin{pmatrix}\mathbb  J & 0\\ 0 &-\mathbb  J\end{pmatrix} \nabla_{(q,p,\Psi,\Pi)}\mathscr H_{SW}(q,p,\Psi,\Pi)
\]
with 
\[
\mathbb J=\begin{pmatrix}
0 & 1 
\\
-1 & 0
\end{pmatrix}\]
and
\begin{align*}
  \mathscr H_{SW}(q,p,\Psi, \Pi)=&\,
\frac12\left(\ds\frac12\ds\int_{\mathbb T^d}  
|\nabla q|^2+|\nabla p|^2 \ud x+ \frac{k^2}{2}\int_{\mathbb T^d}(p^2+q^2)
\ud x-\int_{\T^d}pk\cdot \nabla q\ud x+\int_{\T^d}qk\cdot \nabla p\ud x\right)\\
&+\int_{\T^d\times\R^n}\left({\color{black} c^2}\Pi^2+{\color{black}\frac{1}{4}} |\nabla_z\Psi|^2\right)\ud x \ud z\\
  &+\frac{\gamma}{2} \ds\int_{\mathbb T^d}
  \left(\int_{\T^d\times \R^n}(\sigma_1(x-y)\sigma_2(z)\Psi(t,y,z)\ud y\ud z\right)(p^2+q^2)(x)\ud x.
\end{align*}
Coming back to $u=q+ip$, we can write
\begin{align}\label{ham_SW_k}
  \mathscr H_{SW}(u,\Psi,\Pi)=&\,
  \frac12\left(\ds\frac12\ds\int_{\mathbb T^d}  
  |\nabla u|^2\ud x+ \frac{k^2}{2}\int_{\mathbb T^d}|u(x)|^2
  \ud x+\int_{\T^d}k\cdot(-i \nabla u)\overline{u}\ud x\right)\nonumber\\
  &+\int_{\T^d\times\R^n}\left({\color{black} c^2}\Pi^2+{\color{black}\frac{1}{4}} |\nabla_z\Psi|^2\right)\ud x \ud z\nonumber\\
  &+\frac{\gamma}{2} \ds\int_{\mathbb T^d}
  \left(\int_{\T^d\times \R^n}(\sigma_1(x-y)\sigma_2(z)\Psi(t,y,z)\ud y\ud z\right)|u(x)|^2\ud x.
\end{align}
As a consequence, $\mathscr H_{SW}$ is a constant of the motion.
Moreover, it is clear that \eqref{sw_k} is invariant under multiplications by a phase factor of $u$ so that $F(u)=\frac{1}{2}\|u\|^2_{L^2}$ is conserved by the dynamics.
However, now,  the quantities 
\begin{equation}\label{eq:momentumSchro}
  G_{j}(u)=\frac12 \int_{\T^d}\left(\frac{1}{i}\partial_{x_j}u\right)\overline{u}\ud x
\end{equation}
are not constants of the motion:
\[
\ds\frac{\ud}{\ud t} G_j(u)(t)=\frac\gamma{2} \int_{\T^d}\int_{\T^d}
\partial_{x_j}\sigma_1(x-y)\left(\int_{\mathbb R^n}\sigma_2(z)\Psi(t,y,z)\ud z\right) |u|^2(t,x)\ud y \ud x.
\]
As a consequence, they cannot be used in the construction of the Lyapounov functional as we did for the Hartree system (see \eqref{eq:lyapounov_Hartree}).
% The relation 
% \begin{equation*}
%   H_{SW}(u,\Psi, \Pi)=\mathscr H_{SW}(u,\Psi,\Pi)-\frac{k^2}{2}F(u)-\sum_{j=1}^dk_jG_j(u),
% \end{equation*}
% still holds, but $H_{SW}$ and $\mathscr H_{SW}$ are equivalently conserved only in the case $k=0$. 

Finally, we consider the Banach space
$H^1(\T^d)\times H^1(\T^d)\times L^2(\T^d;\overbigdot {H}^1(\R^n))\times L^2(\T^d\times \R^n)$ endowed  with the inner 
product
\begin{equation*}
\left\langle \begin{pmatrix}
q\\p\\\Psi\\ \Pi\end{pmatrix} \Big|  \begin{pmatrix}
q'
\\
p'\\\Psi'\\ \Pi'\end{pmatrix}\right\rangle=\ds\int_{\mathbb T^d} \big( pp'+qq')\ud x + \int_{\mathbb T^d\times \R^n}(\nabla_z\Psi\nabla_z\Psi'+ \Pi\Pi')\ud x\ud z
\end{equation*}
that can be also interpreted as an inner product for complex valued functions:
\begin{equation}
  \label{eq:innerprodC_SW}
  \langle (u,\Psi,\Pi) | (u',\Psi',\Pi')\rangle =\mathrm{Re} \ds\int_{\mathbb T^d} u\overline{u'}\ud x+ \int_{\mathbb T^d\times \R^n}(\nabla_z\Psi\cdot\nabla_z\Psi'+ \Pi\Pi')\ud x\ud z.
\end{equation}
We denote by $\|\cdot\|$ the norm on $H^1(\T^d) \times L^2(\T^d;\overbigdot {H}^1(\R^n))\times L^2(\T^d\times \R^n)$ induced by this inner product.

\subsection{Preliminary results for the linearized problem: spectral stability when $k=0$}
\label{Prelim}

As before, we linearize the system \eqref{Schro-s}-\eqref{Schro-p} around the plane wave solution obtained in Section~\ref{PWsol}. Namely, we 
expand $$U(t,x)=U_k(t,x)(\mathbf 1+u(t,x)) ,\qquad
\Psi(t,x,z)=-
\gamma\avex{\sigma_1}\Gamma(z)+\psi(t,x,z)$$
 and, assuming that $u,\psi$ and their derivatives are small,
we are led to  the following equations for the fluctuation $(t,x)\mapsto u(t,x)\in\mathbb C$, $(t,x,z)\mapsto \psi(t,x,z)\in\mathbb R$
\begin{equation}
\label{Schro-l}\begin{array}{l}
 i\partial_t u + \ds\frac{ 1}{2}\Delta_x u +   ik\cdot \nabla_x u= \gamma \Phi,
\\[.3cm]
\Big(\ds\frac1{c^2}\partial_{tt}^2\psi-\Delta _z \psi\Big)(t,x,z)=-\gamma 
 
 \sigma_2(z)\sigma_1\star  \rho (t,x),
\\[.3cm]
\rho(t,x)=2\mathrm{Re}\big(u(t,x)\big),
\\[.3cm]
\Phi(t,x)=\ds\iint_{\mathbb T^d\times\mathbb R^n}\sigma_1(x-y)\sigma_2(z)\psi(t,y,z)\ud z\ud y.
\end{array}
\end{equation}
We split the solution into real and imaginary parts
\[
u(t,x)=q(t,x)+ip(t,x),\qquad
q(t,x)=\mathrm{Re}(u(t,x)),\qquad
p(t,x)=\mathrm{Im}(u(t,x)).
\]
We obtain
\begin{equation}
\label{eig2}\begin{array}{l}
(\partial_t q + \ds\frac{ 1}{2}\Delta_x p +   k\cdot \nabla_x q)(t,x)= 0,
\\[.3cm]
(\partial_t p- \ds\frac{ 1}{2}\Delta_x q+   k\cdot \nabla_x p)(t,x)= -\gamma \left(\sigma_1\star   
 \ds\int_{\mathbb R^n}\sigma_2(z)\psi(t,\cdot,z)\ud z\right)(x) ,
\\[.3cm]
\Big(\ds\frac1{c^2}\partial_{tt}^2\psi-\Delta _z \psi\Big)(t,x,z)=-2\gamma    
\sigma_2(z)\sigma_1\star   q(t,x).
\end{array}
\end{equation}

 It is convenient to set 
\[\pi=-\frac{1}{2c^2}\partial_t \psi,\]
in order to rewrite the wave equation as a first order system.
We obtain
\begin{equation}
\label{syst_lin_sw}
\partial_t 
\begin{pmatrix}
q 
\\
 p \\
 \psi
 \\
 \pi
 \end{pmatrix}=
 \mathbb L_k\begin{pmatrix}
q 
\\
 p \\
 \psi
 \\
 \pi
 \end{pmatrix}
 \end{equation}
 where $\mathbb L_k$ is the operator defined by 
 \[
 \mathbb L_k : \begin{pmatrix}
q 
\\
 p \\
 \psi
 \\
 \pi
 \end{pmatrix}
 \longmapsto 
\begin{pmatrix}
-\ds\dfrac12\Delta_x p-k\cdot \nabla_x q
\\
 \ds\frac{ 1}{2}\Delta_x q-    k\cdot \nabla_x p -\gamma \sigma_1\star\left(
  \ds\int_{\mathbb R^n}\sigma_2\psi \ud z\right) 
 \\
 -2c^2\pi
 \\
-\ds\frac{1}{2}\Delta _z \psi+\gamma    
\sigma_2\sigma_1\star   q
\end{pmatrix}.
\]
For the next step, we proceed via Fourier analysis as before. We expand $q$, $p$, $\psi$, $\pi$ and $\sigma_1$ by means of their Fourier series:
\[\begin{array}{ll}
\psi(t,x,z)=\ds\sum_{m\in \mathbb Z^d} \psi_{m}(t,z)e^{im\cdot x},\qquad &
\psi_{m}(t,z)=\ds\frac{1}{(2\uppi)^d}\ds\int_{\mathbb T^d} \psi(t,x,z)e^{-im\cdot x} \ud x,
\\[.3cm]
\pi(t,x,z)=\ds\sum_{m\in \mathbb Z^d} \pi_{m}(t,z)e^{im\cdot x},\qquad &
\pi_{m}(t,z)=\ds\frac{1}{(2\uppi)^d}\ds\int_{\mathbb T^d} \pi(t,x,z)e^{-im\cdot x} \ud x.
\end{array}\]
Moreover, recall that $\sigma_1$ being real and radially symmetric, \eqref{condsig1} holds and, by definition, $\avex{\sigma_1}=(2\uppi)^d \sigma_{1,0}$.

As a consequence, since the Fourier modes are uncoupled, the Fourier coefficients $$(Q_m(t),P_m(t),\psi_m(t,z),\pi_m(t,z))$$ satisfy
\begin{equation}%\label{eq:1}
    \partial_t \begin{pmatrix}
Q_m
\\
P_m
\\
\psi_m
\\
\pi_m
\end{pmatrix}
=
\mathbb L_{k,m}\begin{pmatrix}
Q_m
\\
P_m
\\
\psi_m
\\
\pi_m
\end{pmatrix}\label{eigs3_3}\end{equation}
where $\mathbb L_{k,m}$ stands for the operator defined by 
\[
\mathbb L_{k,m}\begin{pmatrix}
Q_m
\\
P_m
\\
\psi_m
\\
\pi_m
\end{pmatrix}=
\begin{pmatrix}
-ik\cdot m Q_m+ \ds\frac{m^2}{2}  P_m 
\\
-\ds \frac{m^2}{2}Q_m  -ik\cdot m P_m
-
\gamma  (2\uppi)^{d}\sigma_{1,m} \ds\int_{\mathbb R^n} \sigma_2(z)\psi_m\ud z  
\\
-2c^2\pi_m
\\
\gamma  (2\uppi)^d 
 
%c^2 
\sigma_2(z)\sigma_{1,m} Q_m
- \ds\frac{1}{2}\Delta_z \psi_m 
\end{pmatrix}.
\]
Like for the Hartree equation, the behavior of the mode $m=0$ 
can be analysed explicitly.

\begin{lemma}[The mode $m=0$]\label{mod0SW}
For any $k\in\mathbb Z^d$, the kernel of $\mathbb L_{k,0}$ is spanned by $(0,1,0,0)$.
Moreover, equation \eqref{eigs3_3} for $m=0$ admits solutions which grow linearly with time.
\end{lemma}

\begin{Proof}
Let $(Q_0,P_0,\psi_0,\pi_0)\in \mathrm{Ker}(\mathbb L_{k,0})$. It means that 
\begin{equation*}
  \left\{
    \begin{aligned}
      &\gamma  (2\uppi)^{d}\sigma_{1,0} \ds\int_{\mathbb R^n} \sigma_2(z)\psi_0(z)\ud z=0,\\
      &\pi_0=0,\\
      &\Delta_z \psi_0=2\gamma  (2\uppi)^d 
      \sigma_2(z)\sigma_{1,0} Q_0,
    \end{aligned}
  \right.
\end{equation*}
which yields $\psi_0(z)=-2\gamma \avex{\sigma_{1}} Q_0 \Gamma(z)$ with $\Gamma(z)=(-\Delta)^{-1}\sigma_2(z)$ so that 
\begin{equation*}
  -2\gamma^2\avex{\sigma_1}^2\kappa Q_0=0.
\end{equation*}
It implies that $Q_0=0$, $\psi_0=0$ while $P_0$ is left undetermined.

For $m=0$, the first equation in \eqref{eigs3_3} tells us that $Q_0(t)=Q_0(0)\in \C$ is constant.
Next, we get $\partial_t \psi_0= -2c^2 \pi_0$ which leads to
\begin{equation}\label{eqlinmode0SW}
\partial^2_{tt}\psi_0-c^2\Delta_z \psi_0=-\sigma_2(z)
\underbrace{2\gamma 
c^2 \avex{\sigma_{1}}Q_0(0)}_{:=C_1}
\end{equation}
The solution of \eqref{eqlinmode0SW} with initial condition $(\psi_0(z),\pi_0(z)=-\frac{1}{2c^2}\partial_t\psi(0,z))\in  \overbigdot {H}^1( \mathbb R^n)\times L^2( \mathbb R^n)$ satisfies
\[
\widehat\psi_0(t,\xi)=
\widehat\psi_0(0,\xi)\cos(c|\xi|t)-2c^2 \widehat\pi_0(\xi)\ds\frac{\sin(c|\xi|t)}{c|\xi|}
- \ds\int_0^t \ds\frac{\sin(c|\xi|s)}{c|\xi|}\widehat \sigma_2(\xi)C_1\ud s
\]
where $\widehat\psi_0(t,\xi)$ and $\widehat\pi_0(t,\xi)$ are the Fourier transforms of $z\mapsto \psi(t,z)$ and $z\mapsto\pi(t,z)$ respectively. 
Finally, integrating 
\begin{equation*}
  \partial_t P_0=\underbrace{-\gamma\avex{\sigma_{1}}}_{:=C_2}\int_{\mathbb R^n} \sigma_2(z)\psi_0(z)\ud z
\end{equation*}
we obtain
\begin{align*}
  P_0(t)&=P_0(0)+C_2\int_{\mathbb R^n} \widehat \sigma_2(\xi)\widehat\psi_0(0,\xi) \ds\frac{\sin(c|\xi|t)}{c|\xi|}\ds\frac{\ud \xi}{(2\uppi)^n}
  -2c^2C_2\int_{\mathbb R^n} \widehat \sigma_2(\xi)\widehat\pi_0(0,\xi) \ds\frac{1-\cos(c|\xi|t)}{c^2|\xi|^2}\ds\frac{\ud \xi}{(2\uppi)^n}\\
  &\ -C_1C_2 \int_0^t \ds\int_0^s p_c(\tau) \ud \tau \ud s
\end{align*}
where
\[
p_c(\tau) = \ds\int_{\mathbb R^d} |\widehat \sigma_2(\xi)|^2 \ds\frac{\sin(c|\xi|\tau )}{c|\xi|} \ds\frac{\ud \xi}{(2\uppi)^n}.\]
This kernel already appears in the analysis performed in \cite{dBGV,Vi2}.
The contribution involving the initial data of the vibrational field can be uniformly bounded 
by
\[
\ds\frac{1}{(2\uppi)^n}\left( \ds\int_{\mathbb R^d} \ds\frac{|\widehat \sigma_2(\xi)|^2}{c^2|\xi|^2} \ud \xi \right)^{1/2}
\left\{
\left( \ds\int_{\mathbb R^d}|\widehat\psi_0(0,\xi)|^2\ud \xi\right)^{1/2}
+  4c^2\left( \ds\int_{\mathbb R^d}\ds\frac{| \widehat\pi_0(0,\xi)|^2}{c^2|\xi|^2}\ud \xi\right)^{1/2}
\right\}.
\]
Next, as a consequence of \ref{H2}, it turns out that $p_c$ is compactly supported, with $\int_0^\infty p_c(\tau)\ud \tau =\frac{\kappa}{c^2}$, see \cite[Lemma~14]{dBGV}
and \cite[Section~2.4]{Vi2}. 
It follows that 
\[
\ds\int_0^t \ds\int_0^s p_c(\tau) \ud \tau \ud s=
\ds\int_0^t  p_c(\tau)\left( \ds\int_\tau^t   \ud s\right) \ud \tau
=\ds\int_0^t (t-\tau)p_c(\tau) \ud \tau
\underset{t\to \infty}{\sim} t \ds\frac{\kappa}{c^2}- \ds\int_0^\infty \tau p_c(\tau)\ud \tau,\]
which concludes the proof.
\end{Proof}

When $k=0$, basic estimates based on the energy conservation
allow us to justify 
the stability of the solutions with zero mean. However, in contrast to what has been established for the Hartree system, this analysis does not extend 
to any mode $k\neq 0$, since the system is not Galilean invariant.

 \begin{theo}[Linearized stability for the Schr\"odinger-Wave system when $k=0$]
 \label{LinStabk0SW}
Let $k=0$. Suppose \eqref{small} and
let $(u,\psi,\pi)$ be the solution of  \eqref{Schro-l} associated to an initial data $u^{\mathrm{Init}}
\in H^1(\mathbb T^d), \psi^{\mathrm{Init}}\in L^2(\mathbb T^d; \overbigdot H^1( \mathbb R^n)), \pi^{\mathrm{Init}}\in L^2(\mathbb T^d\times \mathbb R^n)$ such that $\int_{\mathbb T^d}u^{\mathrm{Init}}\ud x=0$.
Then, there exists a constant $C>0$ such that $\sup_{t\geq 0}\|u(t,\cdot)\|_{H^1}\leq C$.
\end{theo}

\begin{Proof}
Again, we use the energetic properties of the linearized equation \eqref{Schro-l}. 
We have already remarked that $\int_{\mathbb T^d}u(t,x) \ud x=0$ for any $t\geq 0$ when $\int_{\mathbb T^d}u^{\mathrm{Init}}\ud x=0$.
We start by computing
\[
\begin{array}{l}
\ds\frac{\ud}{\ud t}\left\{\ds\frac12\ds\int_{\mathbb T^d}|\nabla_x u|^2\ud x
+\ds\frac12 
 \ds\int_{\mathbb T^d\times \mathbb R^n}
\Big(\ds\frac{|\partial_t \psi |^2}{c^2} + |\nabla_z\psi|^2\Big) \ud z\ud x\right\}
\\[.3cm]
\quad
= -\ds\frac{i\gamma}{2}\ds\int_{\mathbb T^d} \Phi\Delta_x(u-\overline u)\ud x
- \gamma  
 \ds\int_{\mathbb T^d\times \mathbb R^n} \partial_t \psi \sigma_2 \sigma_1\star  (u+\overline u)\ud z\ud x.
\end{array}\]
Next, we get
\[\begin{array}{lll}
\ds\frac{\ud}{\ud t}
\ds\int_{\mathbb T^d} \Phi (u+\overline u)\ud x
&=&
\ds\int_{\mathbb T^d\times \mathbb R^n} \partial_t \psi \sigma_2 \sigma_1\star  (u+\overline u)\ud z\ud x
\\[.3cm]
&&
+\ds\frac{i}{2}\ds\int_{\mathbb T^d} \Phi\Delta_x(u-\overline u)\ud x
-
\ds\int _{\mathbb T^d} \Phi k\cdot\nabla_x(u+\overline u)\ud x.
\end{array}\]
We get rid of the last term by assuming $k=0$ and we arrive in this case at 
\[
\ds\frac{\ud}{\ud t}\left\{\ds\frac12\ds\int_{\mathbb T^d}|\nabla_x u|^2\ud x
+\ds\frac12
 \ds\int_{\mathbb T^d\times \mathbb R^n}
\Big(\ds\frac{|\partial_t \psi |^2}{c^2} + |\nabla_z\psi|^2\Big) \ud z\ud x
+\gamma\ds\int_{\mathbb T^d} \Phi (u+\overline u)\ud x
\right\}=0.
\] 
We estimate the coupling term as follows
\begin{align*}
&\left|
\ds\int_{\mathbb T^d} \Phi (u+\overline u)\ud x
\right|
=\left|
\ds\int_{\mathbb T^d\times \mathbb R^n} \sigma_2(z) \psi(t,x,z) \sigma_1\star (u+\overline u)(t,x)\ud z\ud x
\right|\\
&\qquad \le \|\sigma_1\star (u+\overline u)\|_{L^2} \times \left(
  \ds\int_{\mathbb T^d} 
  \Big| \ds\int_{\mathbb R^n} \sigma_2(z) \psi(t,x,z)\ud z\Big|^2\ud x
  \right)^{1/2}\\
&\qquad \le \|\sigma_1\|_{L^1}\|u+\overline u\|_{L^2} \times \left(
  \ds\int_{\mathbb T^d} 
  \Big| \ds\int_{\mathbb R^n} \widehat \sigma_2(\xi) \overline{\widehat \psi(t,x,\xi)}\ds\frac{\ud \xi}{(2\uppi)^n}\Big|^2\ud x
  \right)^{1/2}\\
&\qquad \le 2\|\sigma_1\|_{L^1}\|u\|_{L^2} \times \left(
  \ds\int_{\mathbb T^d} 
  \Big| \ds\int_{\mathbb R^n} \ds\frac{\widehat \sigma_2(\xi)}{|\xi|} |\xi|\ds |\overline{\widehat \psi(t,x,\xi)}|\ds\frac{\ud \xi}{(2\uppi)^n}\Big|^2\ud x
  \right)^{1/2}
  \\
&\qquad \le 2\|\sigma_1\|_{L^1}\|u\|_{L^2} \times \left(\ds\int_{\mathbb R^n} \ds\frac{|\widehat \sigma_2(\xi)|^2}{|\xi|^2} \ud \xi\right)^{1/2} \times
\left(
\ds\int_{\mathbb T^d\times\mathbb R^n} 
|\xi|^2 |\widehat \psi(t,x,\xi)|^2 \ds\frac{\ud \xi}{(2\uppi)^n}\ud x
\right)^{1/2}\\
&\qquad \le 2\sqrt \kappa \|\sigma_1\|_{L^1}\|u\|_{L^2} \times\left(
  \ds\int_{\mathbb T^d\times \mathbb R^n} 
   |\nabla_z  \psi(t,x,\xi)|^2 \ud z \ud x\right)^{1/2}= 2\sqrt\kappa \|\sigma_1\|_{L^1}\|u\|_{L^2}\|\nabla_z \psi\|_{L^2}\\
& \qquad 
\leq \ds\frac{1}{2\gamma} \|\nabla_z \psi\|_{L^2}^2 + 2\kappa  \gamma \|\sigma_1\|^2_{L^1}\|u\|^2_{L^2}.
\end{align*}
By using the Poincar\'e-Wirtinger inequality $\|u\|_{L^2}\leq \|\nabla_x u\|_{L^2}$, we deduce that 
\[
\ds\frac12\ds\int_{\mathbb T^d}|\nabla_x u(t,x)|^2\ud x\leq \ds\frac {E_0}{1-4\gamma^2\kappa \|\sigma_1\|^2_{L^1}}
,
\]
where $E_0$ depends on the energy of the initial state.
\end{Proof}

While it is natural to start with the linearized operator $\mathbb L_k$ in \eqref{syst_lin_sw},
it turns out that this formulation 
is not well-adapted to study the spectral stability issue.
The difficulties relies on the fact that the wave part of the system
induces an essential spectrum, reminiscent  to the fact that $\sigma_{\mathrm{ess}}(-\Delta_z)
=[0,\infty)$.
For instance, this is even an obstacle to set up a perturbation argument from the Hartree equation, in the spirit of \cite{HarG}.
We shall introduce  later on a more adapted formulation of the linearized equation, which will allow us to overcome these difficulties  (and also to go beyond a mere perturbation analysis).

\subsection{Orbital stability for the Schr\"odinger-Wave system when $k=0$}

In this subsection, we wish to establish the \emph{orbital stability} of the plane wave solution to \eqref{sw_k} obtained in Section \ref{PWsol}, namely
\begin{equation*}
  u_\omega(t,x)=e^{i\omega t}\mathbf 1(x),\quad \Psi_*(t,x,z)= - \gamma \Gamma(z)\avex{\sigma_1},\quad \Pi_*(t,x,z)=0
\end{equation*}
with $k\in \Z^d$ and $\omega\in \mathbb R$ that satisfy the dispersion relation \eqref{eq:dispersion} and $\Gamma(z)=(-\Delta)^{-1}\sigma_2(z)$.
The system \eqref{sw_k} being invariant under multiplications of $u$ by a phase factor, we define the corresponding orbit through $(\mathbf 1(x),- \gamma \Gamma(z)\avex{\sigma_1}, 0)$ by
\[
\mathscr O_{\mathbf 1}=\{(e^{i\theta},- \gamma \Gamma(z)\avex{\sigma_1}, 0) ,\ \theta \in \mathbb R\}.\]
As before, orbital stability intuitively means that the solutions 
of \eqref{sw_k}
associated to initial data close enough to $(\mathbf 1(x),- \gamma \Gamma(z)\avex{\sigma_1}, 0)$ remain at a close distance to the set $\mathscr O_{\mathbf 1}$.

Let us introduce, for any $k\in \Z^d$ and $\omega\in \mathbb R$ satisfying the dispersion relation \eqref{eq:dispersion}, the set 
\begin{equation*}
  \label{eq:surface_SW_k}
  \mathscr S_\omega=\Big\{(u,\Psi,\Pi)\in H^1(\T^d;\C)\times L^2(\T^d;\overbigdot {H}^1(\R^n))\times L^2(\T^d,L^2(\R^n)),\ F(u)=F(\mathbf 1)=\ds\frac{(2\uppi)^{d}  
}{2}\Big\},
\end{equation*}
 and the functional
\begin{equation}
  \label{eq:lyapounov_SW}
  \mathscr L_{\omega,k}(u,\Psi,\Pi)=\mathscr H_{SW}(u,\Psi,\Pi)+ 
\omega F(u), 
\end{equation}
intended to serve as a Lyapounov functional, where  $\mathscr H_{SW}$ is the constant of motion defined in \eqref{ham_SW_k}. For further purposes, we simply denote $\mathscr L_\omega=\mathscr L_{\omega,0}$.
Note that
\begin{equation*}
  \mathscr L_{\omega,k}(u,\Psi,\Pi)=H_{SW}(u,\Psi, \Pi)+ \underbrace{\ds\frac{1}{2i}\ds\int_{\mathbb T^d} k\cdot\nabla u\ \bar u\ud x}_{=\ds\sum_{j=1}^d k_jG_j(u)}+ 
\Big(\omega + \ds\frac{k^2}{2}\Big)F(u)
\end{equation*}
with $H_{SW}$ defined in \eqref{eq:hamiltonian} and $G_j(u)$ defined in \eqref{eq:momentumSchro}. Thanks to the dispersion relation \eqref{eq:dispersion}, only the second term of this expression depends on $k$. Unfortunately, as pointed out before, the quantities $G_j(u)$ are not constants of the motion so that
the dependence on $k$ of 
the Lyapounov functional \eqref{eq:lyapounov_SW} cannot be disregarded, in contrast 
to what we did for the Hartree system in \eqref{eq:lyapounov_Hartree}.

Next, as in subsection \ref{sec:orbitalHartree}, we need to evaluate the first and second order variations of $\mathscr L_{\omega,k}$.
We compute
\begin{align*}
  &\partial_{(u,\Psi,\Pi)} H_{SW}(u,\Psi,\Pi)(v,\phi,\tau)\\
  &=\mathrm{Re}\left(\ds\frac12
 \ds\int_{\mathbb T^d} (-\Delta u)\overline v\ud x
+\gamma \ds\int_{\mathbb T^d}\left(
\ds\iint_{\mathbb T^d\times\mathbb R^n}
 \sigma_1(x-y)\sigma_2(z)\Psi(t,y,z)\ud z\ud y\right)
u(x) \overline{v(x)}\ud x\right)\\
&\quad +\frac{\gamma}{2} \ds\int_{\mathbb T^d}\left(
  \ds\iint_{\mathbb T^d\times\mathbb R^n}
   \sigma_1(x-y)\sigma_2(z)\phi(t,y,z)\ud z\ud y\right)
  |u(x)|^2\ud y\ud x\\
&\quad +
\ds\frac12
 \ds\iint_{\mathbb T^d\times\mathbb R^n} 
 \Big(\ds{\color{black}4 c^2}
 \Pi\ \tau +(- \Delta_z\Psi ) \ \phi\ud z\Big)\ud x
\end{align*}
and
\begin{align*}
  &\partial^2_{(u,\Psi,\Pi)} H_{SW}(u,\Psi,\Pi)\big((v,\phi,\tau),(v',\phi',\tau')\big)\\
  &= \mathrm{Re}\left\{\ds\frac12
  \ds\int_{\mathbb T^d} (-\Delta v)\overline {v'}\ud x\right.\\
  &\quad \left.
    +\gamma \ds\int_{\mathbb T^d}\left(
    \ds\iint_{\mathbb T^d\times\mathbb R^n}
     \sigma_1(x-y)\sigma_2(z)(\phi(t,y,z) \overline{v'(x)}+\phi'(t,y,z)\overline{v(x)})\ud z\ud y\right)
    u(x) \ud x\right)\\
  &\quad \left.\left.
    +\gamma \ds\int_{\mathbb T^d}\left(
    \ds\iint_{\mathbb T^d\times\mathbb R^n}
     \sigma_1(x-y)\sigma_2(z)\Psi(t,y,z)\ud z\ud y\right)
    v(x) \overline{v'(x)}\ud x\right)\right\}\\
  &\quad +
  \ds\frac12
   \ds\iint_{\mathbb T^d\times\mathbb R^n} 
   \Big(\ds{\color{black}4 c^2}
   \tau\ \tau' + (- \Delta_z\phi ) \ \phi'\ud z\Big)\ud x.
\end{align*}
Besides, we have
\[\begin{array}{ll}
\partial_uF(u)(v)=\ds\mathrm{Re} \left(\ds\int_{\mathbb T^d} u\overline v\ud x\right),
\qquad & 
\partial^2_uF(u)(v,v')=\ds\mathrm{Re}\left( \ds\int_{\mathbb T^d} v\overline{ v'}\ud x\right),
\\
\partial_uG_j(u)(v)=
\mathrm{Im}\left(\int_{\mathbb T^d} \partial_{x_j}u\overline v \ud x\right),
\qquad & 
\partial^2_uG(u)(v,v')=
\mathrm{Im}\left(\int_{\mathbb T^d} \partial_{x_j}v'\overline v \ud x\right).
\end{array}
\]
Accordingly, we are led to
\begin{align*}
  \partial_{(u,\Psi,\Pi)}&\mathscr L_{\omega,k}(\mathbf 1,-\gamma \avex{\sigma_1}\Gamma,0)
(v,\phi,\tau)
\\&=\mathrm{Re}\left(
-\gamma^2\avex{\sigma_1}^2\kappa
 \ds\int_{\mathbb T^d}\overline{v}
 \ud x+\Big(\omega+\ds\frac{k^2}{2} \Big)\ds\int_{\mathbb T^d} \overline{v}\ud x
 +\frac{\gamma}{2}\avex{\sigma_1}
 \ds\iint_{\mathbb T^d\times\mathbb R^n} 
 ( \sigma_2+ \Delta_z\Gamma)  \ \phi\ud z \ud x
 \right)
 \\&=0
\end{align*}
thanks to the dispersion relation \eqref{eq:dispersion} and the definition of $\Gamma$.
Similarly, the second order derivative casts as
\begin{align*}
  &\partial^2_{(u,\Psi,\Pi)}\mathscr L_{\omega,k}(\mathbf 1,-\gamma \avex{\sigma_1}\Gamma,0)
\big((v,\phi,\tau),(v,\phi,\tau)\big)\\
&=\mathrm{Re}\left(\ds\frac12
\ds\int_{\mathbb T^d} (-\Delta v)\overline {v}\ud x +
\ds\frac12 
\ds\iint_{\mathbb T^d\times\mathbb R^n} 
 \Big(\ds{\color{black}4 c^2}\tau^2
 +( - \Delta_z\phi)  \ \phi\ud z\Big)\ud x
 \right.\\
 &\quad +\left.
  2\gamma \ds\int_{\mathbb T^d}\left(
  \ds\iint_{\mathbb T^d\times\mathbb R^n}
   \sigma_1(x-y)\sigma_2(z)\phi(t,y,z)\ud z\ud y\right)
   \overline{v(x)}\ud x\right.\\
  &\quad- \left.
\gamma^2 \avex{\sigma_1}\ds\int_{\mathbb T^d}\left(
\ds\iint_{\mathbb T^d\times\mathbb R^n}
 \sigma_1(x-y)\sigma_2(z)\Gamma(z)\ud z\ud y\right)
v(x) \overline{v(x)}\ud x
+\Big(\omega+\ds\frac{k^2}{2}\Big)\ds\int_{\mathbb T^d} v(x)\overline{ v(x)}\ud x\right)
\\
&\quad+\mathrm{Im} \left(\ds\sum_{j=1}^d k_j 
\ds\int_{\mathbb T^d} \partial_{x_j}v\overline v \ud x\right).
\end{align*}
The forth and fifth integrals combine as
\[
\ds\int_{\mathbb T^d} 
\Big(\omega+\ds\frac{k^2}{2}- \gamma^2\kappa \avex{\sigma_1}^2\Big)
v(x)\overline {v(x)}\ud x=0
\]
which cancels out by virtue of the dispersion relation \eqref{eq:dispersion}.
Hence we get
\begin{align*}
  \partial^2_{(u,\Psi,\Pi)}&\mathscr L_{\omega,k}(\mathbf 1,-\gamma \avex{\sigma_1}\Gamma,0)
\big((v,\phi,\tau),(v,\phi,\tau)\big)\\
&=\mathrm{Re}\left(\ds\frac12
\ds\int_{\mathbb T^d} (-\Delta v)\overline {v}\ud x +
\ds\frac12 
\ds\iint_{\mathbb T^d\times\mathbb R^n} 
 \Big(\ds{\color{black}4 c^2}\tau^2
 +( - \Delta_z\phi)  \ \phi\ud z\Big)\ud x
 \right.\\
 &\quad +\left.
  2\gamma \ds\int_{\mathbb T^d}\left(
  \ds\iint_{\mathbb T^d\times\mathbb R^n}
   \sigma_1(x-y)\sigma_2(z)\phi(t,y,z)\ud z\ud y\right)
   \overline{v(x)}\ud x
   -i\ds\int_{\mathbb T^d}k\cdot\nabla v\ \overline v\ud x
 \right).
\end{align*}

\begin{rmk}\label{rem:ubound_hess_SW}
  Note that 
  the following continuity estimate holds:
  for any $(v,\phi,\tau)\in H^1(\T^d;\C)\times L^2(\T^d;\overbigdot {H}^1(\R^n))\times L^2(\T^d\times \R^n)$,  
\begin{align*}
    &\partial^2_{(u,\Psi,\Pi)}\mathscr L_{\omega,k}(\mathbf 1,-\gamma \avex{\sigma_1}\Gamma,0)
\big((v,\phi,\tau),(v,\phi,\tau)\big)\le \frac12 \|\nabla v\|_{L^2}^2+{\color{black} 2c^2}\|\tau\|^2_{L^2}+\frac12 \|\phi\|^2_{L^2_x\overbigdot {H}^1_z}\\&\quad+2\gamma \kappa^{1/2}\|\sigma_1\|_{L^1}\|v\|_{L^2}\|\phi\|_{L^2_x\overbigdot {H}^1_z}+|k|\|\nabla v\|_{L^2}\|v\|_{L^2}\le \frac12 \left((1+|k|)\|v\|_{H^1}^2+{\color{black} 4c^2}\|\tau\|^2_{L^2}+C \|\phi\|^2_{L^2_x\overbigdot {H}^1_z}\right)\\
&\quad\le \frac{\max({\color{black} 4c^2},1+|k|,C)}{2}\|(v,\phi,\tau)\|^2
\end{align*}
with $C=1+4\gamma^2\kappa \|\sigma_1\|_{L^1}^2$.
\end{rmk}

The functional $\mathscr L_{\omega,k}$ is conserved
by the solutions of \eqref{sw_k}; however
the difficulty relies on justifying its coercivity.
We are only able to answer positively 
 in the specific case $k=0$.
Hence, 
the main result of this subsection restricts to this situation. 

\begin{theo}[Orbital stability for the Schrödinger-Wave system]
  \label{prop:orbital_SW_k} Let $k=0$ 
   and $\omega\in\mathbb R$ such that the dispersion relation \eqref{eq:dispersion} is satisfied. Suppose \eqref{small} holds. Then the plane wave solution $(e^{i\omega t}\mathbf 1(x), -\gamma\Gamma(z)\avex{\sigma}, 0)$ is orbitally stable, \emph{i.e.}
  \begin{align}
    \label{eq:orbital_SW_k}
    &\forall \varepsilon>0,\ \exists \delta>0,\
    \forall (v^{\mathrm {Init}},\phi^{\mathrm {Init}},\tau^{\mathrm {Init}}) \in H^1(\T^d;\C)\times L^2(\T^d;\overbigdot {H}^1(\R^n))\times L^2(\T^d\times \R^n), \nonumber\\
    &\|v^{\mathrm {Init}}-\mathbf 1\|_{H^1}+\|\phi^{\mathrm {Init}}+\gamma\Gamma\avex{\sigma}\|_{L^2_x\overbigdot {H}^1_z}+\|\tau^{\mathrm {Init}}\|_{L^2}<\delta \Rightarrow \sup_{t\ge 0} {\mathrm{dist}}((v(t),\phi(t),\tau(t)),\mathscr O_{\mathbf 1})<\varepsilon
  \end{align}
 {where $ {\mathrm{dist}}((v,\phi,\tau),\mathscr O_{\mathbf 1})=\inf_{{\theta}\in [0,2\uppi[}\|v-e^{i{\theta}}\mathbf 1\|_{H^1}+\|\phi+\gamma\Gamma\avex{\sigma}\|_{L^2_x\overbigdot {H}^1_z}+\|\tau\|_{L^2}$
  and $
  (t,x,z)\mapsto (v(t,x),\phi(t,x,z),\tau(t,x,z))$ stands for the solution of \eqref{sw_k} with Cauchy data 
  $ (v^{\mathrm {Init}},\phi^{\mathrm {Init}},\tau^{\mathrm {Init}})$.}
\end{theo}

Using the same argument as in the case of Theorem~\ref{prop:orbital_Hartree_k}, we can reduce the proof of Theorem \ref{prop:orbital_SW_k} to the following coercivity estimate on the Lyapounov functional (and this is where we use that $\mathscr L_{\omega,k}$ is a conserved quantity). 

\begin{lemma}\label{lem:coer_SW_k} Let $k\in \Z^d$
 and $\omega\in \mathbb R$ such that the dispersion relation \eqref{eq:dispersion} is satisfied.
  Suppose that there exist $\eta>0$ and $c>0$ such that $\forall (w,\psi,\chi)\in \mathscr S_\omega$,
  \begin{align}
    \label{eq:coer_SW}
      \mathrm{dist}((w,\psi,\chi),\mathscr O_{\mathbf 1})<\eta\Rightarrow \mathscr L_{\omega,k}((w,\psi,\chi))-\mathscr L_{\omega,k}((\mathbf 1(x), -\gamma\Gamma(z)\avex{\sigma}, 0))\ge c\ds  \mathrm{dist}((w,\psi,\chi),\mathscr O_{\mathbf 1})^2.
  \end{align}
  Then the the plane wave solution $(e^{i\omega t}\mathbf 1(x), -\gamma\Gamma(z)\avex{\sigma}, 0)$ is orbitally stable.
\end{lemma}

As we have seen before, since $\partial_{(u,\psi,\Pi)} \mathscr L_{\omega,k}((\mathbf 1,-\gamma\Gamma(z)\avex{\sigma}, 0))=0$, the coercivity estimate \eqref{eq:coer_SW} can be obtained from an estimate on the bilinear form $$\partial_{(u,\psi,\Pi)}^2 \mathscr L_{\omega,k}((\mathbf 1,-\gamma \avex{\sigma_1}\Gamma,0))((u,\phi,\tau),(u,\phi,\tau))$$ for any $(u,\phi,\tau)\in T_{\mathbf 1}\mathscr S_{\omega}\cap (T_{\mathbf 1}\mathscr O_{\mathbf 1})^{\perp}$. 
Here the tangent set to $\mathscr S_\omega$ is given by 
\[T_{ \mathbf 1}\mathscr S_\omega=\left\{u\in H^1(\T^d;\C), \mathrm{Re}\left(\int_{\T^d}u(x)\mathbf 1(x)\ud x\right)=0\right\}\times L^2(\T^d;\overbigdot {H}^1(\R^n))\times L^2(\T^d\times \R^n)
.\]
This set is the orthogonal to $(\mathbf 1,0,0)$ with respect to the inner product defined in \eqref{eq:innerprodC_SW}. The tangent set to $\mathscr O_{\mathbf 1}$ (which is the orbit generated by the phase multiplications of $\mathbf 1$) is
\begin{equation*}
  T_{ \mathbf 1}\mathscr O_{\mathbf 1}=\mathrm{span}_{\R}\{(i\mathbf 1,0,0)\}
\end{equation*}
so that 
\[(T_{ \mathbf 1}\mathscr O_\mathbf 1)^\perp=
\left\{u\in H^1(\T^d;\C), \mathrm{Re}\left(i\int_{\T^d}u(x)\mathbf 1(x)\ud x\right)=0\right\}\times L^2(\T^d;\overbigdot {H}^1(\R^n))\times L^2(\T^d\times \R^n).\]

\begin{lemma}\label{lem:coer_SW_hess} Let $k\in \Z^d$
and $\omega\in \mathbb R$ such that the dispersion relation \eqref{eq:dispersion} is satisfied. Suppose that there exists $\tilde c>0$ 
  \begin{equation}
    \label{eq:coer_SW_hess}
    \partial_{(u,\psi,\Pi)}^2 \mathscr L_{\omega,k}((\mathbf 1, -\gamma\Gamma(z)\avex{\sigma}, 0))((u,\phi,\tau),(u,\phi,\tau))\ge \tilde c(\|u\|^2_{H^1}+\|\phi\|^2_{L^2_x\overbigdot {H}^1_z}+\|\tau\|^2_{L^2})=\tilde c\|(u,\phi,\tau)\|^2
  \end{equation}
  for any $(u,\phi,\tau)\in T_{\mathbf 1}\mathscr S_{\mathbf 1}\cap (T_{\mathbf 1}\mathscr O_{\mathbf 1})^{\perp}$.
Then there exist $\eta>0$ and $c>0$ such that \eqref{eq:coer_SW} is satisfied. \end{lemma}

\begin{Proof}
  Let $(w,\psi,\chi)\in \mathscr{S}_{\omega}$ such that $ {\mathrm{dist}}((w,\psi,\chi),\mathscr O_{\mathbf 1})<\eta$ with $\eta>0$ small enough. Hence, $\inf_{{\theta} \in [0,2\uppi)}\|w-e^{i{\theta}}\mathbf 1\|<\eta$ and,
  by means of an implicit function theorem argument (see \cite[Section 9, Lemma 8]{SRN1}), we obtain that there exists ${\theta}  \in [0,2\uppi)$ and $v\in \left\{u\in H^1(\T^d;\C), \mathrm{Re}\left(i\int_{\T^d}u(x)\ud x\right)=0\right\}$ such that 
  \begin{equation*}
    e^{i{\theta} }w=\mathbf 1+v,\qquad \inf_{{\theta} \in [0,2\uppi)}\|w-e^{i{\theta} }\mathbf 1\|\le \|v\|_{H^1}\le C \inf_{{\theta} \in [0,2\uppi)}\|w-e^{i{\theta} }\mathbf 1\|
  \end{equation*}
  for some positive constant $C$. Denote by $\phi(x,z)=\psi(x,z)+\gamma\Gamma(z)\avex{\sigma_1}$. Then $(v,\phi,\chi)\in (T_{\mathbf 1}\mathscr O_{\mathbf 1})^{\perp}$ and $\|(v,\phi,\chi)\|\le C\eta$.  

  Next, we use the fact that $H^1(\T^d)=\left\{u\in H^1(\T^d;\C), \mathrm{Re}\left(\int_{\T^d}u(x)\ud x\right)=0\right\}\oplus \mathrm{span}_{\R}\{\mathbf{1}\}$ to write $(v,\phi,\chi)=(v_1,\phi,\chi)+(v_2,0,0)$ with $(v_1,\phi,\chi)\in T_{   \mathbf 1}\mathscr S_\omega\cap (T_{\mathbf 1}\mathscr O_{\mathbf 1})^{\perp}$ and $v_2\in \mathrm{span}_{\R}\{\mathbf 1\}
  $. Moreover, 
  \begin{equation*}
    \|v_2\|_{H^1}\le {\color{black} \frac{\|v\|_{H^1}^2}{2\|\mathbf 1\|_{L^2}}}
  \end{equation*}
  and $$\|v_1\|_{H^1}\ge \frac12 \|v\|_{H^1}$$
  provided $\|v\|_{H^1}\le \|\mathbf 1\|_{L^2}$. 
As a consequence, if $\|v\|_{H^1}$ is small enough, using that 
\begin{align*}
  \partial^2_{(u,\Psi,\Pi)}&\mathscr L_{\omega,k}(\mathbf 1,-\gamma \avex{\sigma_1}\Gamma,0)
\big((v,\phi,\tau),(v',\phi',\tau')\big)
\le C \|(v,\phi,\tau)\|\|(v',\phi',\tau')\|,
\end{align*}
we obtain
  \begin{align*}
    &\partial^2_{(u,\Psi,\Pi)}\mathscr L_{\omega,k}(\mathbf 1,-\gamma \avex{\sigma_1}\Gamma,0)
    \big((v_1,\phi,\chi),(v_2,0,0)\big)\le C \|(v,\phi,\chi)\|\,\|v\|^2_{H^1}\le C\|(v,\phi,\chi)\|^3,\\
    &\partial^2_{(u,\Psi,\Pi)}\mathscr L_{\omega,k}(\mathbf 1,-\gamma \avex{\sigma_1}\Gamma,0)
    \big((v_2,0,0),(v_2,0,0)\big)\le C\|v\|^4_{H^1}\le C\|(v,\phi,\chi)\|^4.
  \end{align*}
  This leads to 
  \begin{align*}
    \partial^2_{(u,\Psi,\Pi)}&\mathscr L_{\omega,k}(\mathbf 1,-\gamma \avex{\sigma_1}\Gamma,0)
\big((v,\phi,\chi),(v,\phi,\chi)\big)\\
&=\partial^2_{(u,\Psi,\Pi)}\mathscr L_{\omega,k}(\mathbf 1,-\gamma \avex{\sigma_1}\Gamma,0)
\big((v_1,\phi,\chi),(v_1,\phi,\chi)\big)+o(\|(v,\phi,\chi)\|^2).
  \end{align*}
  Finally, let $(w,\psi,\chi)\in \mathscr{S}_{\omega}$ such that $d((w,\psi,\chi),\mathscr O_{\mathbf 1})<\eta$, we have
  \begin{align*}
    \mathscr L_{\omega,k}((w,\psi,\chi))&-\mathscr L_{\omega,k}((\mathbf 1(x), -\gamma\Gamma(z)\avex{\sigma}, 0))=\mathscr L_{\omega,k}((e^{i \theta }w,\psi,\chi))-\mathscr L_{\omega,k}((\mathbf 1(x), -\gamma\Gamma(z)\avex{\sigma}, 0))\\
    &={\color{black} \dfrac12} \partial^2_{(u,\Psi,\Pi)}\mathscr L_{\omega,k}(\mathbf 1,-\gamma \avex{\sigma_1}\Gamma,0)
    \big((v,\phi,\chi),(v,\phi,\chi)\big)+o(\|(v,\phi,\chi)\|^2)\\
    &={\color{black} \dfrac12}\partial^2_{(u,\Psi,\Pi)}\mathscr L_{\omega,k}(\mathbf 1,-\gamma \avex{\sigma_1}\Gamma,0)
    \big((v_1,\phi,\chi),(v_1,\phi,\chi)\big)+o(\|(v,\phi,\chi)\|^2)\\
    &\ge {\color{black} \dfrac{\tilde c}{2}}\|(v_1,\phi,\tau)\|^2+o(\|(v,\phi,\chi)\|^2)\ge {\color{black} \frac{\tilde c}{4}}\|(v,\phi,\tau)\|^2+o(\|(v,\phi,\chi)\|^2)\\
    &\ge {\color{black} \frac{\tilde c}{8}} d((w,\psi,\chi),\mathscr O_{\mathbf 1})^2
  \end{align*}
  where we use $\partial_{(u,\Psi,\Pi)}\mathscr L_{\omega,k}(\mathbf 1,-\gamma \avex{\sigma_1}\Gamma,0)=0$ and $(v_1,\phi,\chi)\in T_{   \mathbf 1}\mathscr S_\omega\cap (T_{\mathbf 1}\mathscr O_{\mathbf 1})^{\perp}$.
\end{Proof}

As before, to prove the orbital stability of the plane solution $(e^{i\omega t}\mathbf 1(x), -\gamma\Gamma(z)\avex{\sigma}, 0)$ it is enough to prove \eqref{eq:coer_SW_hess} for any $(u,\phi,\tau)\in T_{\mathbf 1}\mathscr S_{\mathbf 1}\cap (T_{\mathbf 1}\mathscr O_{\mathbf 1})^{\perp}$. 
Let $(u,\phi,\tau)\in T_{\mathbf 1}\mathscr S_{\mathbf 1}\cap (T_{\mathbf 1}\mathscr O_{\mathbf 1})^{\perp}$ and write $u=q+ip$ with $q,p\in H^{1}(\T^d;\R)$. Then
\begin{align}\label{eq:SW_hess}
  \partial^2_{(u,\Psi,\Pi)}&\mathscr L_{\omega,k}(\mathbf 1,-\gamma \avex{\sigma_1}\Gamma,0)
\big((u,\phi,\tau),(u,\phi,\tau)\big)\nonumber\\
&=\mathrm{Re}\left(\ds\frac12
\ds\int_{\mathbb T^d} (-\Delta u)\overline {u}\ud x +
\ds\frac12 
\ds\iint_{\mathbb T^d\times\mathbb R^n} 
 \Big(\ds{\color{black} 4c^2}\tau^2
 +( - \Delta_z\phi)  \ \phi\ud z\Big)\ud x
\right.\nonumber\\
 &\quad +\left.
  2\gamma \ds\int_{\mathbb T^d}\left(
  \ds\iint_{\mathbb T^d\times\mathbb R^n}
   \sigma_1(x-y)\sigma_2(z)\phi(t,y,z)\ud z\ud y\right)
   \overline{u(x)}\ud x- i\ds\int_{\mathbb T^d}k\cdot \nabla u\ \overline u \ud x\right)
\end{align}
can be reinterpreted as a quadratic form acting on the 4-uplet  $W=(q,p,\phi,\tau)$. To be specific, it expresses as the following quadratic form on $W$,
\begin{align*}
  \mathscr Q(W,W)=&
\ds\frac12\ds\int_{\mathbb T^d} |\nabla p|^2\ud x
+\ds{\color{black}2c^2}\ds\iint_{\mathbb T^d\times\mathbb R^n} |\tau|^2\ud z\ud x+\frac12\ds\int_{\mathbb T^d} |\nabla q|^2\ud x+\frac12\ds \iint_{\mathbb T^d\times\mathbb R^n}( - \Delta_z\phi)  \ \phi\ud x\ud z\\
&+
  2\gamma \ds\int_{\mathbb T^d}\left(
  \ds\iint_{\mathbb T^d\times\mathbb R^n}
   \sigma_1(x-y)\sigma_2(z)\phi(t,y,z)\ud z\ud y
   q(x)\ud x\right)
 +2\ds\int_{\mathbb T^d} qk\cdot\nabla p\ud x.
\end{align*}
The crossed 
term $\int_{\mathbb T^d} qk\cdot\nabla p\ud x$ is an obstacle for proving a coercivity on $\mathscr Q$.

For this reason, let us focus on the case $k=0$.
Since  $(u,\phi,\tau)\in T_{\mathbf 1}\mathscr S_{\mathbf 1}\cap (T_{\mathbf 1}\mathscr O_{\mathbf 1})^{\perp}$, we have
\begin{equation*}
  \int_{\T^d}q\ud x=0\text{ and } \int_{\T^d}p\ud x=0.
\end{equation*}
As a consequence, thanks to the Poincaré-Wirtinger inequality, we deduce, when $k=0$
\begin{align}\label{proof:coerc1_SW}
  \mathscr Q(W,W)\ge& \frac14\|p\|^2_{H^1}+{\color{black}2c^2}\|\tau\|_{L^2}^2 +\frac12\int_{\T^d}|\nabla q|^2\ud x+\frac12\ds \iint_{\mathbb T^d\times\mathbb R^n}( - \Delta_z\phi)  \ \phi\ud x\ud z\nonumber\\
  &+2\gamma \ds\int_{\mathbb T^d}\left(
    \ds\iint_{\mathbb T^d\times\mathbb R^n}
     \sigma_1(x-y)\sigma_2(z)\phi(t,y,z)\ud z\ud y  \right)
     q(x)\ud x
\end{align}
Next, we expand $q$, $\sigma_1$ and $\phi(\cdot,z)$ in Fourier series, \emph{i.e.} $$q(x)=\ds\sum_{m\in \mathbb Z^d} q_{m}e^{im\cdot x},\ \phi(x,z)=\ds\sum_{m\in \mathbb Z^d} \phi_{m}(z)e^{im\cdot x}  \text{ and } \sigma_1(x)=\ds\sum_{m\in \mathbb Z^d} \sigma_{1,m}e^{im\cdot x}.$$ Note that $\overline{\sigma_{1,m}}=\sigma_{1,m}=\sigma_{1,-m}$ since $\sigma_1$ is real and radially symmetric. Moreover, $\int_{\T^d}q\ud x=0$ implies $q_0=0$.
Hence, 
\begin{align*}
  \int_{\mathbb T^d}\left(
    \int_{\mathbb T^d\times\mathbb R^n}
     \sigma_1(x-y)\sigma_2(z)\phi(t,y,z)\ud z\ud y\right)&
     q(x)\ud x\\
     &=(2\uppi)^{2d}\mathrm{Re}\left(\sum_{m\in \Z^d\smallsetminus\{0\}}\sigma_{1,m}q_m\int_{\R^n}\sigma_2(z)\overline{\phi_m(z)}\ud z\right)
\end{align*}
which implies
\begin{align*}
  \frac12\int_{\T^d}|\nabla q|^2&\ud x+\frac12\ds \iint_{\mathbb T^d\times\mathbb R^n}( - \Delta_z\phi)  \ \phi\ud x\ud z\nonumber\\
  &+2\gamma \ds\int_{\mathbb T^d}\left(
    \ds\iint_{\mathbb T^d\times\mathbb R^n}
     \sigma_1(x-y)\sigma_2(z)\phi(t,y,z)\ud z\ud y\right)
     q(x)\ud x  \\
  &=(2\uppi)^{d}\sum_{m\in \Z^d\smallsetminus\{0\}}\mathrm{Re}\left(\frac {m^2}{2}q_m^2+\frac12\int_{\R^n}|\nabla_z \phi_m|^2\ud z+2(2\uppi)^d\gamma\sigma_{1,m}q_m\int_{\R^n}\sigma_2(z)\overline{\phi_m(z)}\ud z\right).
\end{align*}

Next, we remark that for any $m\in \Z^d$, 
\begin{align*}
  \left|\mathrm{Re}\left(2(2\uppi)^d\gamma\sigma_{1,m}q_m\int_{\R^n}\sigma_2(z)\overline{\phi_m(z)}\ud z\right)\right|&\le 2(2\uppi)^d\gamma\sigma_{1,m} |q_m| \sqrt{\kappa} \|\nabla \phi_m\|_{L^2}\\
  &\le \frac{1}{2\tilde\delta}(4\gamma^2\kappa(2\uppi)^{2d}\sigma_{1,m}^2)q_m^2+\frac{\tilde\delta}{2}\|\nabla \phi_m\|_{L^2}^2
\end{align*}
for any $\tilde \delta>0$.
Finally, for any $\tilde \delta\in (0,1)$, we get
\begin{align}\label{proof:coerc2_SW}
  \frac12\int_{\T^d}|\nabla q|^2&\ud x+\frac12\ds \iint_{\mathbb T^d\times\mathbb R^n}( - \Delta_z\phi)  \ \phi\ud x\ud z\nonumber\\
  &+2\gamma \ds\int_{\mathbb T^d}\left(
    \ds\iint_{\mathbb T^d\times\mathbb R^n}
     \sigma_1(x-y)\sigma_2(z)\phi(t,y,z)\ud z\ud y\right)
     q(x)\ud x \nonumber\\
  &\ge (2\uppi)^d \sum_{m\in \Z^d}\left(\left(\frac{m^2}{2}-\frac{1}{2\tilde\delta}(4\gamma^2\kappa(2\uppi)^{2d}\sigma_{1,m}^2)\right)q^2_m+\frac{1-\tilde\delta}{2}\|\nabla \phi_m\|_{L^2}^2\right)
\end{align}
As a consequence, we obtain the following statement.

\begin{proposition}
  \label{prop:coer_SW_hess_small} Let $k=0$ 
   and $\omega\in\mathbb R$ such that the dispersion relation \eqref{eq:dispersion} is satisfied. Suppose that there exists $\delta\in (0,1)$ such that
  \begin{equation}
    \label{cond_small_SW}
    4\gamma^2\kappa(2\uppi)^{2d}\frac{\sigma^2_{1,m}}{m^2}\le \delta
  \end{equation}
  for all $m\in \Z^d\smallsetminus\{0\}$.
  Then, there exists $\tilde c>0$ such that  
  \begin{equation}
    \label{eq:coer_SW2}
    \partial^2_{(u,\Psi,\Pi)}\mathscr L_\omega(\mathbf 1,-\gamma \avex{\sigma_1}\Gamma,0)
\big((u,\phi,\tau),(u,\phi,\tau)\big)\ge \tilde c\|(u,\phi,\tau)\|^2
  \end{equation}
  for any $(u,\phi,\tau)\in T_{\mathbf 1}\mathscr S_{\mathbf 1}\cap (T_{\mathbf 1}\mathscr O_{\mathbf 1})^{\perp}$.
\end{proposition}

\begin{Proof}
  If \eqref{cond_small_SW} holds, then, for any $\tilde \delta \in (\delta,1)$, \eqref{eq:SW_hess}-\eqref{proof:coerc1_SW}-\eqref{proof:coerc2_SW} lead to 
  \begin{align*}
    \partial^2_{(u,\Psi,\Pi)}&\mathscr L_\omega(\mathbf 1,-\gamma \avex{\sigma_1}\Gamma,0)
\big((u,\phi,\tau),(u,\phi,\tau)\big)\ge \frac14\|p\|^2_{H^1}+{\color{black}2c^2}\|\tau\|_{L^2}^2 \nonumber\\
    &+ \frac{\tilde \delta-{\delta}}{2{\tilde \delta}} (2\uppi)^d\sum_{m\in \Z^d\smallsetminus \{0\}}{m^2}q_m^2+ \frac{1-\tilde \delta}{2}(2\uppi)^d\sum_{m\in \Z^d}\|\nabla \phi_m\|_{L^2}^2 \\
    &=\frac14\|p\|^2_{H^1}+\frac{1}{2c^2}\|\tau\|_{L^2}^2+\frac{\tilde \delta-{\delta}}{2{\tilde \delta}}\|\nabla q\|^{{\color{black} 2}}_{L^2}+\frac{1-\tilde \delta}{2}\|\phi\|^{{\color{black} 2}}_{L^2_x{H}^1_z}\ge \tilde c \|(u,\phi,\tau)\|^2
  \end{align*}
  where in  the last inequality we used the Poincaré-Wirtinger inequality together with the fact that $\int_{\T^d}q\ud x=0$.
\end{Proof}

Finally, Proposition \ref{prop:coer_SW_hess_small} together with Lemma \ref{lem:coer_SW_hess} and Lemma \ref{lem:coer_SW_k} gives Theorem \ref{prop:orbital_SW_k} and the orbital stability of the plane wave solution $(e^{i\omega t}\mathbf 1(x), -\gamma\Gamma(z)\avex{\sigma}, 0)$ in the case $k=0$.

\begin{rmk} 
  The coercivity of $\partial^2_{(u,\Psi,\Pi)}\mathscr L_\omega(\mathbf 1,-\gamma \avex{\sigma_1}\Gamma,0)
  \big((u,\phi,\tau),(u,\phi,\tau)\big)$ on $T_{\mathbf 1}\mathscr S_{\mathbf 1}\cap (T_{\mathbf 1}\mathscr O_{\mathbf 1})^{\perp}$ can be recovered from the spectral properties of a convenient unbounded linear operator $\mathbb S$. Indeed, as we have seen before, by decomposing $u$ into real and imaginary part, the quadratic form defined by \eqref{eq:SW_hess} (with $k=0$) can be written as
  \begin{equation*}
    \mathscr Q(W,W)=
\ds\frac12\ds\int_{\mathbb T^d} |\nabla p|^2\ud x
+\ds{\color{black}2c^2}\ds\iint_{\mathbb T^d\times\mathbb R^n} |\tau|^2\ud z\ud x
+\left \langle \mathbb S\begin{pmatrix}q\\\phi\end{pmatrix}\Big|
\begin{pmatrix}q\\\phi\end{pmatrix}\right\rangle
  \end{equation*}
  with $\mathbb S: H^2(\T^d)\times L^2(\T^d;\overbigdot {H}^1(\R^n))\subset L^2(\T^d)\times L^2(\T^d;\overbigdot {H}^1(\R^n))\to L^2(\T^d)\times L^2(\T^d;\overbigdot {H}^1(\R^n))$ the unbounded linear operator given by
  \[
\mathbb S\begin{pmatrix}q\\\phi\end{pmatrix}=
\begin{pmatrix}-\ds\frac12\Delta_x q +\gamma\sigma_1\star \ds\int_{\mathbb R^n}\sigma_2 \phi\ud z
\\
\ds\frac12
\phi+\gamma \Gamma \sigma_1\star q
 \end{pmatrix}
\]
(where we remind the reader that $\Gamma=(-\Delta)^{-1}\sigma_2)$) and the inner product
\[
\left \langle \begin{pmatrix}q\\\phi\end{pmatrix}\Big|
\begin{pmatrix}q'\\\phi'\end{pmatrix}\right\rangle
=\ds\int_{\mathbb T^d} qq'\ud x+\ds\int_{\mathbb T^d\times\mathbb R^n} \nabla_z\phi \cdot \nabla_z\phi' \ud z\ud x=\int_{\mathbb T^d} qq'\ud x+\ds\int_{\mathbb T^d\times\mathbb R^n} \hat\phi(x,\xi)
 \overline{\hat\phi'(x,\xi)} \frac{|\xi|^2\ud \xi}{(2\uppi)^n}\ud x.
\]
Note that $L^2(\T^d)\times L^2(\T^d;\overbigdot {H}^1(\R^n))$ is an Hilbert space with this inner product since $n\ge 3$.

Since
\begin{align*}
  \int_{T^d}&\left(\sigma_1\star\int_{\mathbb R^n}\sigma_2 \phi\ud z\right)(x)q'(x)\ud x=\int_{\T^d}\left(\int_{\T^d\times\R^n}\sigma_1(x-y)\sigma_2(z)\phi(y,z)\ud z\ud y\right)q'(x)\ud x\\
  &=\int_{\T^d\times \R^n}\phi(x,z)\sigma_2(z)(\sigma_1\star q')(x)\ud x\ud z= \int_{\T^d\times \R^n}\hat\phi(x,\xi)\frac{ \overline{\hat\sigma_2(\xi)}}{|\xi|^2}(\sigma_1\star q')(x)\ud x\frac{|\xi|^2\ud \xi}{(2\uppi)^n}
\end{align*}
we can  check that $\mathbb S$ is a self-adjoint operator on $L^2(\T^d)\times L^2(\T^d;\overbigdot {H}^1(\R^n))$. In particular, $\sigma(\mathbb S)\subset\R$ and one can easily study the spectrum of $\mathbb S$.

More precisely, using Fourier series, we find that if $\lambda$ is an eigenvalue of $\mathbb S$ then there exists at least one $m\in \Z^d$ such that for some $(q_m,\phi_m)\neq(0,0)$ there holds
\begin{equation*}
  \left\{\begin{aligned}
    &\left(\frac{m^2}{2}-\lambda\right)q_m+\gamma (2\uppi)^d\sigma_{1,m}\int_{R^n}\sigma_2(z)\phi_m(z)\ud z=0,\\
    &\left(\frac{1}{2}-\lambda\right)\phi_m(z)+\gamma (2\uppi)^d \Gamma(z)\sigma_{1,m}q_m=0.
  \end{aligned}\right.
\end{equation*}

Let $\lambda\neq \frac12$. Hence, for any $m\in \Z^d$, $q_m=0$ implies $\phi_m(z)=0$ for any $z\in \R^n$. As a consequence, we may assume $q_m\neq 0$. This leads to $\phi_m(z)=-\frac{\gamma (2\uppi)^d\sigma_{1,m}q_m}{{1/2-\lambda}}\Gamma(z)$ and 
\begin{equation*}
  \left(\frac{m^2}{2}-\lambda\right)\left(\frac{1}{2}-\lambda\right)-\gamma^2 (2\uppi)^{2d}\sigma_{1,m}^2\kappa=0.
\end{equation*}
By solving this equation, we obtain
\begin{equation*}
  \lambda_{\pm,m}=\frac{\left(\frac{m^2+1}{2}\right)\pm\sqrt{\left(\frac{m^2-1}{2}\right)^2+4\gamma^2 (2\uppi)^{2d}\sigma_{1,m}^2\kappa}}{2}
\end{equation*}
so that $\lambda_{+,m}\ge \frac{1}{4}$ for any $m\in \Z^d$. Next, we remark that 
\begin{equation*}
  \lambda_{-,0}=\frac{\frac{1}{2}-\sqrt{\frac{1}{4}+4\gamma^2 (2\uppi)^{2d}\sigma_{1,0}^2\kappa}}{2}<0
\end{equation*}
since $4\gamma^2\kappa (2\uppi)^{2d}\sigma_{1,0}^2
>0$. This eigenvalue corresponds to an eigenfunction $(\tilde q, \tilde \phi)$ with $\tilde q\in \mathrm{span}_\R\{\mathbf 1\}$. In particular, $\int_{\T^d} \tilde q(x)\ud x\neq 0$.
Finally, if \eqref{cond_small_Hartree} holds, 
\begin{equation*}
  \lambda_{-,m}\ge \frac{\left(\frac{m^2+1}{2}\right)-\sqrt{\left(\frac{m^2-1}{2}\right)^2+\delta m^2}}{2}\ge \frac{1-\delta}{5}
\end{equation*}
for any $m\in \Z^d\smallsetminus \{0\}$.

We conclude  that 
\begin{equation*}
 \left\langle \mathbb S\begin{pmatrix}q\\\phi\end{pmatrix}\Big| \begin{pmatrix}q\\\phi\end{pmatrix}\right\rangle =
\left\langle\begin{pmatrix}-\ds\frac12\Delta_x q +\gamma\sigma_1\star \ds\int_{\mathbb R^n}\sigma_2 \phi\ud z
\\
\ds\frac12
\phi+\gamma\Gamma \sigma_1\star q
 \end{pmatrix}
 \Big| \begin{pmatrix}q\\\phi\end{pmatrix}\right\rangle 
 \ge \min\left(\frac12,\frac{1-\delta}{5}\right)(\|q\|^2_{L^2}+\|\phi\|_{L^2_x\overbigdot {H}^1_z})
\end{equation*}
for all $(q,\phi)\in \{q\in L^2(\T^d),\int_{T^d}q\ud x=0\}\times L^2(\T^d;\overbigdot H^1(\R^n))$. This, together with the Poincaré-Wirtinger inequality, proves the coercivity of $\partial^2_{(u,\Psi,\Pi)}\mathscr L_\omega(\mathbf 1,-\gamma \avex{\sigma_1}\Gamma,0)
\big((u,\phi,\tau),(u,\phi,\tau)\big)$ on $T_{\mathbf 1}\mathscr S_{\mathbf 1}\cap (T_{\mathbf 1}\mathscr O_{\mathbf 1})^{\perp}$.
\end{rmk}

\section{Discussion about the case $k\neq 0$}
\label{S:kn0}

\subsection{A new symplectic form of the linearized Schr\"odinger-Wave system}

We go back to the linearized 
problem. The viewpoint presented in Section~\ref{Prelim} looks quite natural; however, it misses some structural properties of the problem. In order to work in a unified functional framework, we find convenient to change
the wave unknown $\psi$, which is naturally valued in $\overbigdot H^1(\mathbb R^n)$, 
  into  $ (-\Delta)^{-1/2}\phi$, where the new unknown $\phi$ now lies in $L^2(\mathbb R^n)$. 
  The last component of the unknown vector $X$ becomes 
  $\pi =- \frac{ (-\Delta)^{-1/2}\partial_t \phi}{c}$.
  (The change of unknowns allows us to work in a convenient unified functional framework, based on $L^2$ spaces; the constants are chosen in order to make symmetry properties appear, see Lemma~\ref{basic} and the continuity estimate after \eqref{Duhamel} below.)
  Hence, the linearized problem is rephrased as
  \[\partial _t X=\mathbb L X,\]
  where $X$ stands for the $4$-uplet $(q,p,\phi,\pi)$ and 
  \begin{equation}\label{defopLd}
  \mathbb L X=
\begin{pmatrix}
-\ds\dfrac12\Delta_x p-k\cdot \nabla_x q
\\
 \ds\frac{ 1}{2}\Delta_x q-    k\cdot \nabla_x p -\gamma \sigma_1\star\left(
  \ds\int_{\mathbb R^n}(-\Delta)^{-1/2}\sigma_2\phi \ud z\right) 
 \\
 - c (-\Delta)^{1/2} \pi
 \\
 c(-\Delta)^{1/2} \phi
+ 2c\gamma    
\sigma_2\sigma_1\star   q
\end{pmatrix}.\end{equation}
  The operator $\mathbb L$ is seen as an operator on the  Hilbert space
  $$\mathscr V=L^2(\mathbb T^d)\times L^2(\mathbb T^d)\times L^2(\mathbb T^d;L^2(\mathbb R^n))
  \times  L^2(\mathbb T^d;L^2(\mathbb R^n)),$$ with domain 
  $D(\mathbb L)
  =H^2(\mathbb T^d)\times H^2(\mathbb T^d)\times L^2(\mathbb T^d;H^1(\mathbb R^n))
  \times  L^2(\mathbb T^d;H^1(\mathbb R^n)).
  $
  {\color{black} 
  The considered functional framework is  now
made of complex valued functions, which makes 
   the space $\mathscr V$ a complex Hilbert space 
   when
   endowed with the norm $\|\cdot \|_{\mathscr V}$ based on the $L^2$ inner product on each component.
   We are thus going to study the spectral properties of $\mathbb L$ on the space $\mathscr V$. }
  We can start with the following basic information, which has the consequence 
  that the spectral stability amounts to justify that $\sigma(\mathbb L)\subset i\mathbb R$.

  \begin{lemma}\label{basic}
  Let $(\lambda,X)$ be an eigenpair of $\mathbb L$. 
  Let $Y: (x,z)\mapsto (q(-x), -p(-x), \phi(-x,z), -\pi(-x,z))$.
  Then, $(\overline \lambda, \overline X)$, $(-\lambda, Y)$ and $(-\overline \lambda, \overline Y)$ are equally  eigenpairs of $\mathbb L$.
  \end{lemma}
  
  \begin{Proof}
  Since $\mathbb L$ has real coefficients, $\mathbb LX=\lambda X$ implies $\mathbb L\overline X=\overline\lambda \overline X$. Next, we check that
  \[\begin{array}{lll}
  \mathbb LY(x,z)&=&\begin{pmatrix}
\ds\dfrac12\Delta_x p+k\cdot \nabla_x q
\\
 \ds\frac{ 1}{2}\Delta_x q-    k\cdot \nabla_x p -\gamma \sigma_1\star\left(
  \ds\int_{\mathbb R^n}(-\Delta)^{-1/2}\sigma_2\phi \ud z'\right) 
 \\
 c (-\Delta)^{1/2} \pi
 \\
 c (-\Delta)^{1/2} \phi
+ 2c \gamma    
\sigma_2\sigma_1\star   q
\end{pmatrix}(-x,z) \\
&=&\lambda \begin{pmatrix}
-q(-x,z)\\
p(-x,z)
\\
-\phi(-x,z)
\\
\pi(-x,z)\end{pmatrix}
=-\lambda Y(x,z)
.\end{array}\]
  \end{Proof}

  Next, we make a new symplectic structure appear.
  To this end, let us introduce the blockwise operator
  \[
  \mathscr J=\begin{pmatrix}
  \mathscr J_1 & 0 \\
  0 & \mathscr J_2
  \end{pmatrix},
  \qquad
   \mathscr J_1=\begin{pmatrix}
  0 & 1 \\
  -1 &0
  \end{pmatrix},
  \qquad
   \mathscr J_2= 2c \begin{pmatrix}
0 & -(-\Delta)^{1/2}\\
  (-\Delta)^{1/2} & 0
  \end{pmatrix}.
  \]
  We are thus led to 
  \[
  \mathbb L=\mathscr J\mathscr L\]
  with 
 \begin{equation}\label{lin_self_adj_opSW}
  \mathscr LX=
  \begin{pmatrix}
 -  \ds\frac{ 1}{2}\Delta_x q+    k\cdot \nabla_x p +
 \gamma \sigma_1\star\left(
  \ds\int_{\mathbb R^n}(-\Delta)^{-1/2}\sigma_2\phi \ud z\right) 
\\
-\ds\dfrac12\Delta_x p-k\cdot \nabla_x q
 \\
 \ds\frac  \phi 2
+\gamma    
(-\Delta)^{-1/2}\sigma_2\sigma_1\star   q\\
 \ds\frac \pi2
\end{pmatrix}.
\end{equation}
For further purposes, we also set 
\begin{equation}
  \label{sympl_op_SW}
  \widetilde {\mathscr J}=
\begin{pmatrix}
\tilde {\mathscr J}_1&0  \\
  0 &\tilde {\mathscr J}_2
  \end{pmatrix},
  \qquad
  \tilde {\mathscr J}_1=
 \begin{pmatrix}
 0 & -1
 \\
 1 & 0
 \end{pmatrix},\qquad
\tilde {\mathscr J}_2=  \ds\frac 1{2c} \begin{pmatrix}
 0 & (-\Delta )^{-1/2}
 \\
 - (-\Delta )^{-1/2} & 0 \end{pmatrix}.
\end{equation}
The operator $\mathscr J$ has 0 in its essential spectrum; nevertheless 
$\widetilde{\mathscr J}$ plays the role of its 
inverse since $\mathscr J\widetilde{\mathscr J}=\mathbb I=\widetilde{\mathscr J}\mathscr J$.

\begin{lemma}
The operator $\mathscr L$ is an unbounded self adjoint operator on $\mathscr V$
with domain $D(\mathscr L)=
H^2(\mathbb T^d)\times H^2(\mathbb T^d)\times L^2(\mathbb T^d;L^2(\mathbb R^n))\times L^2(\mathbb T^d;L^2(\mathbb R^n))$, and the operator $\mathscr J$ is skew-symmetric.
\end{lemma}

\begin{Proof}
The space $\mathscr V$ is endowed with the standard $L^2$ inner product
\[
\big(X|X')=\ds\int_{\mathbb T^d} (q\overline{q'} +p\overline{p'})\ud x
+\ds\iint_{\mathbb T^d\times\mathbb R^n} (\phi\overline{\phi'} +\pi\overline{\pi'})\ud x\ud z
.\]
 We get
 \[\begin{array}{lll}
 \big(\mathscr L X| X'\big)&=&
\ds\int_{\mathbb T^d}
\Big\{\Big( -  \ds\frac{ 1}{2}\Delta_x q+    k\cdot \nabla_x p \Big)\overline{q'}
+ \Big(-\ds\dfrac12\Delta_x p-k\cdot \nabla_x q\Big)\overline{p'}\Big\}\ud x
\\
[.3cm]
&&+
 \gamma 
 \ds\int_{\mathbb T^d} \sigma_1\star\left(
  \ds\int_{\mathbb R^n}(-\Delta)^{-1/2}\sigma_2\phi \ud z\right) \overline{q'}\ud x
\\[.3cm]
&&+  \ds\frac12
\ds\iint_{\mathbb T^d\times \mathbb R^n}\Big(
  \phi \overline{\phi'}+ %2c^2 
  \pi\overline{\pi'}\Big)\ud x\ud z
\\[.3cm]
&&+\gamma
\ds\iint_{\mathbb T^d\times \mathbb R^n}
    
\Big((-\Delta)^{-1/2}\sigma_2\sigma_1\star   q \Big) \overline{\phi'}\ud x\ud z
\\
[.3cm]
&=&
\ds\int_{\mathbb T^d}\Big\{
q
\Big( -  \ds\frac{ 1}{2}\Delta_x \overline{q'}+    k\cdot \nabla_x \overline{p'} \Big)
+p  \Big(-\ds\dfrac12\Delta_x \overline{p'}-k\cdot \nabla_x \overline{q'}\Big)\Big\}\ud x
\\
[.3cm]
&&+
 \gamma 
 \ds\iint_{\mathbb T^d\times\mathbb R^n} 
 \phi (-\Delta)^{-1/2}\sigma_2 \sigma_1\star \overline{q'}
   \ud z \ud x
\\[.3cm]
&&+  \ds\frac12
\ds\iint_{\mathbb T^d\times \mathbb R^n}\Big(
  \phi \overline{\phi'}+ %2c^2
   \pi\overline{\pi'} \Big)\ud x\ud z
\\[.3cm]
&&+\gamma
\ds\int_{\mathbb T^d} q \sigma_1\star\left(\ds\int_{\mathbb R^n}
(-\Delta)^{-1/2}\sigma_2 \overline{\phi'}\ud z\right)\ud x
\\[.3cm]
&=&
  \big( X| \mathscr L X'\big),
 \end{array}\]
 and
 \[\begin{array}{lll}
  \big( \mathscr J X|  X'\big)
  &=&
   \ds\iint_{\mathbb T^d}
  \Big( p\overline{q'}-q
    \overline{p'}\Big)
   \ud x
  +
 2c  \ds\iint_{\mathbb T^d\times\mathbb R^n} 
    \Big( -(-\Delta)^{1/2}\pi\overline{\phi'}+
   (-\Delta)^{1/2}\phi \overline{\pi'}\Big)
\ud x\ud z
\\[.3cm]
 &=&
 -  \ds\iint_{\mathbb T^d}
  \Big( q
    \overline{p'} - p\overline{q'}\Big)
   \ud x
  -
 2c \ds\iint_{\mathbb T^d\times\mathbb R^n} 
    \Big( -
   \phi \overline{(-\Delta)^{1/2}\pi'}+
   \pi\overline{(-\Delta)^{1/2} \phi'}\Big)
\ud x\ud z
  \\[.3cm]
  &=&
 - \big(  X|  \mathscr J X'\big)
  \end{array}\]
\end{Proof}

As said above, justifying the spectral stability for the Schr\"odinger-Wave equation
reduces to verify that the spectrum $\sigma(\mathbb L)$ 
is purely  imaginary. However, the coupling 
with the wave equation induces delicate subtleties 
and a direct approach is not obvious. Instead, 
based on the expression $\mathbb L=\mathscr J\mathscr L$,
we can take advantage of stronger structural properties.
In particular, the functional framework adopted here allows us to overcome the difficulties
related to the essential spectrum induced by the wave equation, 
which ranges over all the imaginary axis. 
This approach is strongly inspired 
by the methods introduced by D. Pelinovsky and M. Chugunova
\cite{ChouPel,PelBk,PelProc}.
The workplan can be summarized as follows.
It can be shown that the eigenproblem $\mathbb L X=\lambda X$ for $\mathbb L$ 
is equivalent to a generalized eigenvalue problem $\mathbb AW=\alpha \mathbb KW$, with $\alpha=-\lambda^2$, see Proposition~\ref{equiv1} and~\ref{equiv2} below, where the 
auxiliary operators $\mathbb A$ and $\mathbb K$ 
depend on $\mathscr J, \mathscr L$.
Then, we need to identify
negative  eigenvalues and complex but non real eigenvalues of the generalized eigenproblem.
To this end, we appeal to a counting statement due to \cite{ChouPel}.

\subsection{Spectral properties of the operator $\mathscr L$}

The stability analysis relies on the spectral properties of $\mathscr L$, collected in the following claim.

\begin{proposition}
\label{SpecL} Let $\mathscr{L}$ the linear operator defined by \eqref{lin_self_adj_opSW} on $D(\mathscr L)\subset \mathscr V$. Suppose \eqref{small}. Then, the following assertions hold:
\begin{enumerate}
\item  $\sigma_{\mathrm{ess}}(\mathscr L)=\{1/2\}$, %\left\{\frac{1}{2},2c^2\right\}$,
\item $\mathscr L$ has a finite number of \emph{negative} eigenvalues, with eigenfunctions in $D(\mathscr L)$, given by
\[\begin{array}{lll}
n(\mathscr L)&=&1+\#\{m\in \mathbb Z^d\smallsetminus\{0\}, m^4-4(k\cdot m)^2<0 \text{ and } 
\sigma_{1,m}=0\}
\\
&&+
\#\{m\in \mathbb Z^d\smallsetminus\{0\}, m^4-4(k\cdot m)^2\leq 0 \text{ and } 
\sigma_{1,m}\neq 0\}.
\end{array}\]
In particular, $n(\mathscr L)=1$ when $k=0$.
The eigenspaces associated to the negative eigenvalues are all finite-dimensional.
\item With $X_0=(0,\mathbf 1,0,0)$, we have 
$\mathrm{span}_{\mathbb R}\{X_0\}\subset \mathrm{Ker}(\mathscr L)$.
Moreover, given $k\in \mathbb Z^d\smallsetminus\{0\}$, 
let $\mathscr K_*=\{m\in \mathbb Z^d\smallsetminus\{0\},\ 
m^4-4(k\cdot m)^2= 0 \text{ and } 
\sigma_{1,m}= 0\}$. Then, we get
$\mathrm{dim}(\mathrm{Ker}(\mathscr L))=1+\#\mathscr K_*$.
\end{enumerate}
\end{proposition}

We remind the reader that $\sigma_1$ is assumed radially symmetric, see \ref{H1}.
Consequently $\sigma_{1,m}=\sigma_{1,-m}=\overline{\sigma_{1,\pm m}}$
and both $\#\mathscr K_*$ and 
$\#\{m\in \mathbb Z^d\smallsetminus\{0\}, m^4-4(k\cdot m)^2\leq 0 \text{ and } 
\sigma_{1,m}\neq 0\}$ are necessarily even.
\\

\begin{Proof}
Since $\mathscr L$ is self-adjoint, $\sigma(\mathscr L)\subset \mathbb R$.
Let us study the eigenproblem for $\mathscr L$: $\lambda X=\mathscr LX$ means
\begin{equation}
  \label{eigenpb_L_SW}
  \left\{
  \begin{aligned}
    &\lambda q=
 -  \ds\frac{ 1}{2}\Delta_x q+    k\cdot \nabla_x p +
 \gamma \sigma_1\star\left(
  \ds\int_{\mathbb R^n}(-\Delta)^{-1/2}\sigma_2\phi \ud z\right) ,
\\
&\lambda p=
-\ds\dfrac12\Delta_x p-k\cdot \nabla_x q,
 \\
 &\lambda \phi=
 \ds\frac12 \phi
+\gamma    
(-\Delta)^{-1/2}\sigma_2\sigma_1\star   q,
\\
&\lambda \pi=
%2c^2 
\ds\frac12 \pi.
  \end{aligned}
  \right.
\end{equation}

Clearly   $\lambda=\frac12%2c^2
$ is an eigenvalue with  eigenfunctions of the form $(0,0,0,\pi)$, $\pi\in L^2(\mathbb T^d\times\mathbb R^n)$. As a consequence,  $\dim(\Ker(\mathscr L-\frac12%2c^2 I
))$ is not finite and %$2c^2
$\frac12\in \sigma_{\mathrm{ess}}(\mathscr L)$.

We turn to the case   $ \lambda\neq \frac12%2c^2
$, where the last equation imposes $\pi=0$.
Using Fourier series, we obtain 
\begin{equation}\label{SysFou}
\begin{array}{l}
\lambda q_m=
   \ds\frac{ m^2}{2} q_m+    ik\cdot m p_m +
 \gamma (2\uppi)^d \sigma_{1,m}\left(
  \ds\int_{\mathbb R^n}(-\Delta)^{-1/2}\sigma_2\phi_m \ud z\right) ,
\\[.3cm]
\lambda p_m=
\ds\dfrac{m^2}2 p_m-ik\cdot m q_m,
 \\[.3cm]
 \lambda \phi_m=
 \ds\frac12  \phi_m
+\gamma  (2\uppi)^d  
(-\Delta)^{-1/2}\sigma_2\sigma_{1,m}   q_m.
\end{array}\end{equation}
where $q_m, p_m\in \C$ are the Fourier coefficients of $q,p \in L^2(\T^d)$ while $\phi_m(z)=\frac{1}{(2\uppi)^d}\int_{\T^d}\phi(x,z)e^{-im\cdot x}\ud x$ for all $z\in \R^n$ and $\phi\in L^2(\T^d; L^2(\R^n))$.

We split the discussion into several cases. 
\\

{\bf Case $m=0$.}
For $m=0$, the equations \eqref{SysFou} degenerate to
\[\begin{array}{l}
\lambda q_0=
 \gamma (2\uppi)^d \sigma_{1,0}\left(
  \ds\int_{\mathbb R^n}(-\Delta)^{-1/2}\sigma_2\phi_0 \ud z\right) ,
\\
\lambda p_0=0,
 \\
  \Big(\lambda-\ds\frac12
\Big)
 \phi_0=
\gamma   (2\uppi)^d 
(-\Delta)^{-1/2}\sigma_2\sigma_{1,0}   q_0.
\end{array}\]
Combining the first and the third equation yields
\[  \lambda 
\Big(\lambda-
\ds\frac12\Big) 
q_0=\gamma^2(2\uppi)^{2d}\sigma_{1,0}^2 \kappa q_0,\]
still with $\kappa =\int (-\Delta)^{-1}\sigma_2\sigma_2\ud z$.
It permits us to identify the following eigenvalues:
\begin{itemize}
\item 
 $\lambda=0$ is an eigenvalue associated to the eigenfunction
$(0,\mathbf 1,0,0)$,
\item
since $\sigma_{1,0}=\frac{1}{(2\uppi)^d}\int_{\mathbb T^d}\sigma_1\ud x\neq 0$, and 
$(-\Delta)^{-1/2}\sigma_2\neq 0$, $\lambda=1/2$ is an eigenvalue associated to eigenfunctions 
$(0,0,\phi,0)$, for any function $z\mapsto \phi(z)$ orthogonal to $(-\Delta)^{-1/2}\sigma_2$.  
 We find another infinite dimensional eigenspace  associated to the eigenvalue $\lambda=\frac12$.
%As before, since $\dim(\Ker(\mathscr L-1%\tfrac12 I
%))$ is not finite, $\tfrac12 \in \sigma_{\mathrm{ess}}(\mathscr L)$.
\item the roots of 
\[\lambda   \Big(\lambda-
\ds\frac12\Big)
-\gamma^2(2\uppi)^{2d}\sigma_{1,0}^2 \kappa 
=
\lambda^2 -\ds\frac\lambda 2
 -\gamma^2(2\uppi)^{2d}\sigma_{1,0}^2 \kappa=0,\]
provide two additional eigenvalues
\[
\lambda_\pm=  \ds\frac{1/2
\pm \sqrt{1/4
 +4\gamma^2(2\uppi)^{2d}\sigma_{1,0}^2 \kappa}}{2},
\]
associated to the eigenfunctions 
$(\mathbf 1,0,\frac{\gamma (2\uppi)^{d}\sigma_{1,0}(-\Delta)^{-1/2}\sigma_2}{\lambda_\pm-1/2
},0)$, respectively.
\end{itemize}
To sum up, the Fourier mode $m=0$ gives rise to two positive eigenvalues (1/2 and $\lambda_+$), one negative eigenvalue ($\lambda_-$) and the eigenvalue 0, the last two being both one-dimensional.
It tells us that 
$$
\mathrm{dim}(\mathrm {Ker}(\mathscr L))\geq 1 \text{ and } n(\mathscr L)\geq 1.$$
\\

{\bf Case $m\neq 0$ with $\sigma_{1,m}=0$.}
In this case, the $m$-mode equations \eqref{SysFou} for the particle and the wave are uncoupled 
\[
 (\lambda-1/2)
 \phi_m=0,\qquad (M_m-\lambda)\begin{pmatrix}
q_m
\\
p_m
\end{pmatrix}=0\]
where we have introduced the $2\times 2$ matrix
\begin{equation}\label{defmatM}
M_m=\begin{pmatrix}
m^2/2 & ik\cdot m
\\
-ik\cdot m & m^2/2
\end{pmatrix}.\end{equation}
We 
identify the following eigenvalues:
\begin{itemize}
\item 
 $\lambda=1/2$
   is an eigenvalue associated to the eigenfunction
$(0,0,e^{im\cdot x}\phi(z), 0)$,
for any $\phi\in L^2(\mathbb R^n)$. Once again, this tells us that 
 $\tfrac12 \in \sigma_{\mathrm{ess}}(\mathscr L)$.
\item the eigenvalues
$\lambda_\pm=\frac{m^2\pm 2k\cdot m}{2}\in\mathbb R$ of the $2\times 2$ matrix $M_m$, associated to the eigenfunctions $(e^{im\cdot x}, \mp i e^{im\cdot x}, 0, 0)$,
respectively.
Since $\mathrm{tr}(M_m)>0$, at most only one of these eigenvalues can be negative, which occurs when $\mathrm{det}(M_m)=\frac{m^4}4-(k\cdot m)^2<0$.
\end{itemize}
Given $k\in\mathbb Z^d$, we conclude this case by asserting
\[
n(\mathscr L)\geq 1+\#\{m\in \mathbb Z^d\smallsetminus\{0\},\ m^4-4(k\cdot m)^2<0,\ \sigma_{1,m}=0\},
\]
and
\[
\mathrm{dim}(\mathrm{Ker}(\mathscr L))\geq 1+ 
\#\{m\in \mathbb Z^d\smallsetminus\{0\},\ m^2= \pm 2k\cdot m,\ \sigma_{1,m}=0\}.\]
\\

{\bf Case $m\neq 0$ with $\sigma_{1,m}\neq 0$.}
Again, we distinguish several subcases.
\begin{itemize}
\item if  $\lambda=\tfrac12$,
 the third equation on \eqref{SysFou} imposes $q_m=0$, 
and we are led to 
\[
\ds\frac{1-m^2}{2}p_m=0,\qquad
ik\cdot m p_m+\gamma (2\uppi)^d \sigma_{1,m}\left(
  \ds\int_{\mathbb R^n}(-\Delta)^{-1/2}\sigma_2\phi_m \ud z\right)=0.\]
Thus,  $\lambda=\tfrac12$
 is an eigenvalue associated to the eigenfunctions:  
\[ \text{$
(0,0,e^{im\cdot x}\phi(z),0)$, for any function $z\mapsto \phi(z)$ orthogonal to $(-\Delta)^{-1/2}\sigma_2$,}\]
(we recover the same eigenfunctions as for the case $m=0$),
\[(0,e^{im\cdot x},0,0)\textrm{ if $k\cdot m=0$, $m^2=1$}, 
\]
  and
  \[
  \Big(0,-\ds\frac{\gamma(2\uppi)^d \kappa \sigma_{1,m}}{i k\cdot m}e^{im\cdot x},(-\Delta)^{-1/2}\sigma_2(z)e^{im\cdot x},0\Big)\textrm{ if $k\cdot m\neq 0$, $m^2=1$}.
  \]
  \item if $\lambda =\frac{m^2}{2}\neq \frac12$, \eqref{SysFou} becomes
  \[\begin{array}{l}
 0=
      ik\cdot m p_m +
 \gamma (2\uppi)^d \sigma_{1,m}\left(
  \ds\int_{\mathbb R^n}(-\Delta)^{-1/2}\sigma_2\phi_m \ud z\right) ,
\\
0=
-ik\cdot m q_m,
 \\
 \ds\frac{m^2-1}{2} \phi_m=
\gamma  (2\uppi)^d  
(-\Delta)^{-1/2}\sigma_2\sigma_{1,m}   q_m.
  \end{array}\]
There is no non-trivial solution when $k\cdot m\neq 0$.
Otherwise, we see that  $\lambda=m^2/2$ is an eigenvalue associated to the eigenfunctions:  
$(0,e^{im\cdot x},0,0)$ 
\item if   $\lambda\notin\{ %\frac12
\tfrac12,\frac{m^2}2\}$, we set $\mu=\lambda- \frac{m^2}2$.
We see that a non trivial solution of \eqref{SysFou} exists if its component $q_m$ does not vanish. We combine the equations in \eqref{SysFou} to obtain
\[
P(\mu)q_m=0\]
where $P$ is the third order polynomial
\[\begin{array}{l}
P(\mu)=
\mu^3 +b\mu^2+c\mu +d,
\\
b=\ds\frac{m^2-1}{2}\geq 0,\qquad
c=-((k\cdot m)^2 +\gamma^2\kappa (2\pi)^{2d}\sigma_{1,m}^2)< 0,\qquad
d=- (k\cdot m)^2
\ds\frac{m^2-1}{2}\leq 0.
\end{array}.
\]
Observe that $d=- (k\cdot m)^2b$ and $(k\cdot m)^2<|c|<(k\cdot m)^2+\frac14$.
We thus need to examine the roots of this polynomial.
To this end, we compute the discriminant
\[
\mathcal D=18 bcd-4b^3 d+b^2 c^2 -4c^3-27  d^2.\]
A tedious, but elementary, computation allows us to reorganize terms as follows
\[\begin{array}{lll}
\mathcal D&=&
4(k\cdot m)^2 \big((k\cdot m)^2-b^2\big)^2 
+ b^2\sigma_{1,m}^2 \gammaup(20 (k\cdot m)^2   +\gammaup\sigma_{1,m}^2)
\\[.3cm]&&
+ 4(k\cdot m)^2 \sigma_{1,m}^2 \gammaup (2(k\cdot m)^2   +\gammaup\sigma_{1,m}^2)
+ 4\sigma_{1,m}^2 \gammaup\big(
(k\cdot m)^4+ 2(k\cdot m)^2  \sigma_{1,m}^2 \gammaup+\sigma_{1,m}^4 \gammaup^2
\big)
,\end{array}\]
where we have set $ \gammaup=\gamma^2\kappa(2\uppi)^{2d}$.
Since $\sigma_{1,m}\neq 0$, we thus have $\mathcal D>0$ and $P$ has 3 distinct real roots, 
$\mu_1<\mu_2<\mu_3$. In order to bring further information about the location of the roots, 
we observe that $\lim_{\mu\to \pm\infty}P(\mu)=\pm \infty$ while
$P(0)=d\leq 0$ and $P'(0)=c<0$. 
Moreover, studying the zeroes of $P'(\mu)=3\mu^2 +2b\mu +c$, we see that 
$\mu_{\mathrm{max}}=\frac{-b-\sqrt{b^2-3c}}{3} <0$ is a local maximum  and 
$\mu_{\mathrm{min}}=\frac{-b+\sqrt{b^2-3c}}{3} >0$ is a local minimum.
Moreover, $P''(\mu)=6\mu+2b$, showing that $P$ is convex on the domain 
$(-(m^2-1)/6,+\infty)$, concave on $(-\infty,-(m^2-1)/6)$.
A typical shape of the polynomial $P$ is depicted in Figure~\ref{polP}.
From this discussion, we infer 
$$\mu_1<\mu_{\mathrm{max}}<\mu_2\leq 0<\mu_{\mathrm{min}}<\mu_3.$$ 

\begin{figure}[!h]
\begin{center}
\includegraphics[height=6cm]{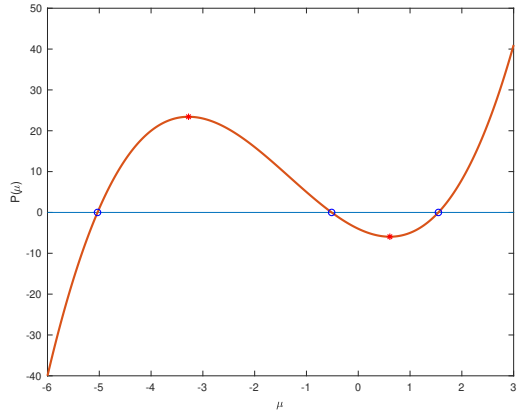}
\end{center}
\caption{Typical graph for $\mu\mapsto P(\mu)$, with its roots $\mu_1<\mu_2<\mu_3$ and local extrema $\mu_{\mathrm{max}}$, $\mu_{\mathrm{min}}$}
\label{polP}
\end{figure}

Coming back to the issue of counting the negative eigenvalues of $\mathscr L$, we are thus wondering whether or not $\lambda_j=\mu_j+m^2/2$ is negative.
We already know that $\mu_3>\mu_{\mathrm{min}}>0$, hence $\mu_3>-m^2/2$ and  we have at most 2 negative eigenvalues for each Fourier mode $m\neq 0$ such that $\sigma_{1,m}\neq 0$.
To decide how many negative eigenvalues should be counted, 
we look at the sign of $P(-m^2/2)$ (see Fig.~\ref{polP}):
\begin{itemize}
\item [i)]
 if $P(-m^2/2)>0$ then $\mu_1<-m^2/2<\mu_2$,

\item [ii)]  if $P(-m^2/2)=0$ then 
either $-m^2/2<\mu_{\mathrm{max}}$, in which case $\mu_1=-m^2/2<\mu_2$,
or $-m^2/2>\mu_{\mathrm{max}}$, in which case $\mu_2=-m^2/2>\mu_1$,

\item [iii)]   if $P(-m^2/2)<0$ then either $-m^2/2<\mu_{\mathrm{max}}$, in which case $-m^2/2<\mu_1<\mu_2$, or
$-m^2/2>\mu_{\mathrm{max}}$, in which case $\mu_1<\mu_2<-m^2/2$.
\end{itemize}
However, we remark that 
\begin{equation}\label{sgpm2}\begin{array}{lll}
P(-m^2/2)&=&
-\ds\frac{m^6}{8}+\ds\frac{m^4(m^2-1)}{8}+\ds\frac{m^2}{ 2}((k\cdot m)^2+\gammaup\sigma_{1,m}^2)
-\ds\frac{m^2-1}{2}(k\cdot m)^2
\\
[.3cm]&=&
-\ds\frac{m^4}{8}\Big(1-4\ds\frac{\gammaup \sigma_{1,m}^2}{m^2}
\Big)+
\ds\frac{(k\cdot m)^2}{2}
=-\ds\frac18
(m^4-4(k\cdot m)^2 -4m^2 \gammaup \sigma^2_{1,m}),
\end{array}\end{equation}
where, by virtue of \eqref{small}, $m\neq 0$ and $\sigma_{1_m}\neq 0$, $1>4\frac{\gammaup \sigma_{1,m}^2}{m^2}> 0$.

This can be combined together with
\begin{align*}
  P'(-m^2/2)=&
3\ds\frac{m^4}{4}-\ds\frac{m^2(m^2-1)}{2}-(k\cdot m)^2- \gammaup\sigma_{1,m}^2
=\ds\frac{m^4}{4}+\ds\frac{m^2}{2}-(k\cdot m)^2- \gammaup\sigma_{1,m}^2
\\
=& \frac{1}{4}\left(m^4 -4(k\cdot m)^2-4m^2\gammaup\sigma_{1,m}^2\right)+m^2\gammaup\sigma_{1,m}^2+\frac{m^2}{2}-\gammaup\sigma_{1,m}^2\\
=&-2P(-m^2/2)+\frac{m^2}{2}+(m^2-1)\gammaup\sigma_{1,m}^2>-2P(-m^2/2).
\end{align*}
Finally,
\begin{equation*}
  P''(-m^2/2)=-2m^2-1<0.
\end{equation*}

As a consequence,  $P(-m^2/2)<0$ implies $P'(-m^2/2)>0$, while 
$P''(-m^2/2)<0$. This shows that $-m^2/2<\mu_1$.
Therefore, in case iii), the only remaining possibility is the situation where $P(-m^2/2)<0$ with $-m^2/2<\mu_1<\mu_2$. As a conclusion, if  $P(-m^2/2)<0$, all eigenvalues $\lambda_j$ are positive.

Next, we claim that case ii) cannot occur. Indeed, $P(-m^2/2)=0$ if and only if $$(m^2-2k\cdot m )(m^2+2k\cdot m) =4m^2 \gammaup \sigma^2_{1,m}.$$
In particular, the term on the left hand side of this equality has to be positive. This is possible if and only if both factors, which belong to $\mathbb Z$, are positive. In this case, according to the sign of $k\cdot m$, one of them is $\geq m^2$ so that
$$
m^2 \le 4m^2 \gammaup \sigma^2_{1,m}.
$$ 
This contradicts the smallness condition \eqref{small}. 
Note that $P(-m^2/2)\neq 0$ implies $\lambda_j\neq 0$, \emph{i.e.} $m$-modes with $m\neq 0$ and $\sigma_{1,m}\neq 0$ cannot generate elements of $\Ker(\mathscr L)$.

As a conclusion, negative eigenvalues only come from  case i) and for each $m$-mode such that $P(-m^2/2)>0$  we have exactly one negative eigenvalue. Going back to \eqref{sgpm2},
in this case, we have
\[(m^4-4(k\cdot m)^2)=(m^2-2k\cdot m )(m^2+2k\cdot m)<m^2 4\gammaup \sigma_{1,m}^2
<m^2\]
owing to \eqref{small}. This excludes the possibility that $m^4-4(k\cdot m)^2>0$, since we noticed above that whenever this term is positive, it is $\geq m^2$.
Hence, case i) holds if and only if $m^4-4(k\cdot m)^2\leq 0$. 
\end{itemize}

This ends the counting of the negative eigenvalues of $\mathscr L$ in Proposition~\ref{SpecL}.
Note that the associated eigenspaces are spanned by
$$\Big(e^{im\cdot x}, - \frac{ik\cdot m}{\lambda-m^2/2}e^{im\cdot x},
e^{im\cdot x}\frac{\sigma_{1,m}\gamma(2\uppi)^d(-\Delta_z)^{-1/2}\sigma_2}{\lambda-1/2},0\Big).$$
\\

The discussion has permitted us to find the elements of $\mathrm{Ker}(\mathscr L)$.
To be specific, the equation $\mathscr LX=0$ yields  $\pi=0$ and the following relations for the Fourier coefficients 
\[\begin{array}{l}
\ds\frac{m^2}{2} p_m-ik\cdot mq_m=0,
\\ [.3cm]
 \ds\frac {\phi_m}{2}+(2\uppi)^{d}\gamma(-\Delta)^{-1/2}\sigma_2\sigma_{1,m} q_m=0,
\\[.3cm]
\ds\frac {m^2}{2}q_m +ik\cdot m p_m+(2\uppi)^{d}\gamma \sigma_{1,m}\ds\int 
(-\Delta)^{-1/2}\sigma_2\phi_m \ud z=0.\end{array}\]
We have seen that the mode $m=0$ gives rise the eigenspace spanned by $(0,\mathbf 1, 0, 0)$. 
For $m\neq 0$, elements of $\Ker(\mathscr L)$ can be obtained only in the case $\sigma_{1,m}=0$. Moreover, the condition $m^2=\pm 2 k\cdot m$ has to be fulfilled. In such a case,  $(e^{im\cdot x}, \mp i e^{im\cdot x}, 0, 0)\in \Ker(\mathscr L)$.\\

Finally, it remains to prove that 
 $\sigma_{\mathrm{ess}}(\mathscr L)=\left\{\frac{1}{2} %,2c^2
\right\}$.  We have already noticed that %$\left\{\frac{1}{2},2c^2\right\}\subset
 $\frac12$ lies in $\sigma_{\mathrm{ess}}(\mathscr L)$. Suppose, by contradiction, that there exists $\lambda \in \sigma_{\mathrm{ess}}(\mathscr L)$ with $\lambda \neq \frac12$.  %\smallsetminus \left\{\frac{1}{2},2c^2\right\}$. 
 Hence, by Weyl's criterion \cite[Theorem~B.14]{PelBk}, there exists a sequence  $(X_\nu)_{\nu \in\mathbb N}$ with $X_\nu=(q_\nu,p_\nu,\phi_\nu,\pi_\nu)\in D(\mathscr L)$ such that,  as $\nu$ goes to $\infty$,
\begin{equation}
  \label{weyl_ess_spectrum}
  \|(\mathscr L-\lambda I)X_\nu\|\to 0,\quad \|X_\nu\|=1 \text{ and } X_\nu \rightharpoonup 0
\textrm{ weakly in $\mathscr V$}.
\end{equation}
Since $\lambda\neq\tfrac12$ and $\lambda\neq 2c^2$, from  \eqref{eigenpb_L_SW} and \eqref{weyl_ess_spectrum} we have 
\begin{equation*}
  \|\pi_\nu\|_{L^2(\T^d;L^2(\R^n))}\to 0 \text{ and } 
 \phi_\nu=- \left(\frac12-\lambda\right)^{-1} \gamma    
 (-\Delta)^{-1/2}\sigma_2\sigma_1\star   q_\nu +\varepsilon_\nu
\end{equation*}
with $\varepsilon_\nu\in L^2(\T^d;L^2(\R^n))$ such that $\lim_{\nu\to \infty}\|\varepsilon_\nu\|_{L^2(\T^d;L^2(\R^n))}= 0$. This leads to 
\begin{equation*}
  \begin{aligned}
  &\left\|-  \ds\frac{ 1}{2}\Delta_x q_\nu -\lambda q_\nu +    k\cdot \nabla_x p_\nu
    -\frac{\gamma^2\kappa}{1/2-\lambda}\Sigma\star q_\nu
  +
  \gamma \sigma_1\star\left(
   \ds\int_{\mathbb R^n}(-\Delta)^{-1/2}\sigma_2\varepsilon_\nu \ud z\right)\right\|_{L^2(\T^d)}\xrightarrow[\nu\to\infty]{} 0,\\
  &\left\|-\ds\dfrac12\Delta_x p_\nu-\lambda p_\nu-k\cdot \nabla_x q_\nu\right\|_{L^2(\T^d)}\xrightarrow[\nu\to\infty]{}  0.
  \end{aligned}
\end{equation*}
Using the fact that the sequence $((q_\nu,p_\nu,\varepsilon_\nu))_{\nu\in \mathbb N}$ is bounded in $L^2(\T^d)\times L^2(\T^d)\times L^2(\T^d;L^2(\R^n))$, we deduce that $(q_\nu,p_\nu)_{\nu\in \mathbb N}$ is bounded in $H^2(\T^d)\times H^2(\T^d)$.
 Indeed, reasoning on Fourier series, this amounts to estimate
 $$\begin{array}{l}
 \ds\sum_{m\in \mathbb Z^d} |m|^4(|q_{\nu,m}|^2 + |p_{\nu,m}|^2) 
 \\
 \qquad\leq
  2 \ds\sum_{m\in \mathbb Z^d}  \big(| m^2 q_{\nu,m} + 2ik\cdot m p_{\nu,m} |^2 + | m^2 p_{\nu,m}-2 ik\cdot m q_{\nu,m}|^2) 
   \\
    \qquad \qquad \qquad+ 8\ds\sum_{m\in \mathbb Z^d} (|k\cdot m p_{\nu,m}|^2 + |k\cdot m q_{\nu,m}|^2)
   \\
   \qquad \leq 
   2\big\| -\Delta_xq_\nu + 2k\cdot\nabla_xp_\nu\big\|_{L^2(\mathbb T^d)}
   + 2\big\| -\Delta_xp_\nu - 2k\cdot\nabla_xq_\nu\big\|_{L^2(\mathbb T^d)}
   \\
    \qquad \qquad \qquad+ \ds\frac4\delta
    |k|^4 \ds\sum_{m\in \mathbb Z^d}\big(|q_{\nu,m}|^2+|p_{\nu,m}|^2\big)
    +4\delta
    \ds\sum_{m\in \mathbb Z^d} |m|^4 (|q_{\nu,m}|^2 + | p_{\nu,m}|^2).
 \end{array}.
 $$
 Choosing $0<\delta<1/4$ and using the already known estimates, we conclude that 
 $\|\Delta_x q_\nu\|^2_{L^2}+\|\Delta_x p_\nu\|^2_{L^2}=
 \sum_{m\in\mathbb Z^d} |m|^4\big(|q_{\nu,m}|^2+|p_{\nu,m}|^2\big)$ is bounded, uniformly with respect to $\nu$.
 Hence, because of the compact Sobolev embedding of $H^2(\T^d)$ into $L^2(\T^d)$, we have that $(q_\nu,p_\nu)_{\nu\in \mathbb N}$ has a (strongly) convergent subsequence in $L^2(\T^d)\times L^2(\T^d)$. As a consequence, the sequence $(X_\nu)_{\nu\in \mathbb N}$ has a convergent subsequence in $\mathscr V$, which contradicts \eqref{weyl_ess_spectrum}.
\end{Proof}

A consequence of Proposition~\ref{SpecL} is that $0$ is an isolated eigenvalue of $\mathscr L$.
Since the restriction of $\mathscr L$ to the subspace  $(\mathrm{Ker}(\mathscr L))^\perp $ is, by definition,  injective, it makes sense to define on it
its inverse
$\mathscr L^{-1}$, with domain $\mathrm{Ran}(\mathscr L)\subset (\mathrm{Ker}(\mathscr L))^\perp\subset \mathscr V$. 
In fact, 0 being an isolated eigenvalue, $\mathrm{Ran}(\mathscr L)$ is closed and coincides with 
$(\mathrm{Ker}(\mathscr L))^\perp$, \cite[Section~B.4]{PelBk}.
This can be shown by means of spectral measures.
Given $X\in (\mathrm{Ker}(\mathscr L))^\perp$, the support of the associated spectral measure 
$\ud\mu_X$ does not meet the interval $(-\epsilon,+\epsilon)$ for $\epsilon>0$ small enough, independent of $X$.
Accordingly, we get
\[\|\mathscr L X\|^2 = \ds\int_{-\infty}^{+\infty} \lambda^2 \ud \mu_X(\lambda)
=\ds\int_{|\lambda|\geq \epsilon} \lambda^2 \ud \mu_X(\lambda)\geq \epsilon^2\|X\|^2.\]
In particular, the Fredholm alternative applies: for any $Y\in (\mathrm{Ker}(\mathscr L))^\perp$, there exists a unique $X\in (\mathrm{Ker}(\mathscr L))^\perp$ such that
$\mathscr LX=Y$. We will denote $X=\mathscr L^{-1}Y$.

For further purposes, let us set 
\[X_0= (0,\mathbf 1, 0, 0)\in \mathrm{Ker}(\mathscr L)
\textrm { and } 
Y_0={\color{black} \mathscr JX_0}=(\mathbf 1,0,0,0).
\]
Note that $Y_0\in (\mathrm{Ker}(\mathscr L))^\perp$, so that it makes sense to consider the equation \[\mathscr LU_0=Y_0.\] We find
\[\pi_m=0,\quad \phi_m=-2\gamma(2\uppi)^d (-\Delta)^{-1/2}\sigma_2 \sigma_{1,m}q_m,\quad
m^2p_m=2ik\cdot m q_m,
\]
and
\[
m^2q_m +2ik\cdot m p_m + 2 \gamma (2\uppi)^d\sigma_{1,m}
\ds\int 
(-\Delta)^{-1/2}\sigma_2\phi_m\ud z=\delta_{0,m}.
\]
It yields, for $m\neq 0$,
$(\frac{m^4}4 - (k\cdot m)^2 - \gammaup |\sigma_{1,m}|^2m^2\Big)q_m=0$
and $
q_0  =-\frac{1}{2\gamma^2(2\uppi)^{2d}|\sigma_{1,0}|^2\kappa}$.
Therefore, we can set
\[
U_0=\mathscr L^{-1}Y_0=-\ds\frac{1}{2\gamma^2(2\uppi)^{2d}|\sigma_{1,0}|^2\kappa}\big(\mathbf 1,0,-2\gamma(2\uppi)^{d}
(-\Delta)^{-1/2}\sigma_2\sigma_{1,0},0\big),\]
solution of $\mathscr LU_0=Y_0$ such that $U_0\in (\mathrm{Ker}(\mathscr L))^\perp$.
We note that 
\begin{equation}\label{CalNo}
(U_0,Y_0)=-\ds\frac{1}{2\gamma^2(2\uppi)^{d}|\sigma_{1,0}|^2\kappa}<0.\end{equation}

\subsection{Reformulation of the eigenvalue problem, and counting theorem}

The aim of the section is to introduce several reformulations of the eigenvalue problem.
This will allow 
us to make use of general counting arguments, set up by \cite{ChouPel,PelBk,PelProc}.

\begin{proposition}\label{equiv1}
Let us set $\mathscr M=-\mathscr J\mathscr L\mathscr J$.
The coupled system
\begin{equation}\label{CouEig}
\mathscr MY=-\lambda X,\qquad \mathscr LX=\lambda Y,\end{equation}
admits a solution with $\lambda\neq 0$, $X\in D(\mathscr  L)\smallsetminus\{0\}$, $Y \in D(\mathscr J\mathscr 
 L\mathscr J)\smallsetminus\{0\}$ 
iff there exists two vectors $X_\pm \in D(\mathbb L)\smallsetminus\{0\}$ that satisfy $\mathbb LX_\pm=\pm\lambda X_\pm$.
\end{proposition}

\noindent
Let $\mathscr P$ stand for the orthogonal projection from $\mathscr V$ to
$(\mathrm{Ker}(\mathscr L))^\perp\subset \mathscr V$.

\begin{proposition}\label{equiv2}
Let us set $\mathbb A=\mathscr P\mathscr M\mathscr P$
and $\mathbb K=\mathscr P\mathscr L^{-1}\mathscr P$.
 Let us define the following Hilbert space
 \[\mathscr H=D(\mathscr M)\cap (\mathrm{Ker}(\mathscr L))^\perp\subset \mathscr V.\]
 The coupled system \eqref{CouEig} has a pair of non trivial  solutions $(\pm \lambda,X,\pm Y)$, with $\lambda\neq 0$ 
 iff the generalized eigenproblem
 \begin{equation}\label{Gene}
 \mathbb A W=\alpha \mathbb K W,\qquad W\in \mathscr H,\end{equation}
admits the  eigenvalue 
$\alpha=-\lambda^2\neq 0$, with 
at least two
linearly independent eigenfunctions.
\end{proposition}

Recall that the plane wave solution obtained Section \ref{PWsol} is spectrally stable, if the spectrum of $\mathbb L$ is contained in $i\R$.  In view of Propositions~\ref{equiv1} and~\ref{equiv2}, this happens if and only if all the eigenvalues of the generalized eigenproblem \eqref{Gene} are real and positive. In other words, the presence of spectrally unstable directions corresponds to the existence of negative eigenvalues or complex but non real eigenvalues of the generalized eigenproblem \eqref{Gene}. 

Our goal is then to count the eigenvalues $\alpha$ of the generalized eigenvalue problem \eqref{Gene}. In particular we define the following quantities: 
\begin{itemize}
\item $N^-_n$, the number of negative eigenvalues
\item $N^0_n$, the number of  eigenvalues zero
\item $N^+_n$, the number of positive eigenvalues
\end{itemize}
of \eqref{Gene}, counted with their algebraic multiplicity, the eigenvectors of which are associated to 
non-positive 
values of the the quadratic form $W\mapsto (\mathbb KW|W)=(\mathscr L^{-1}\mathscr PW|\mathscr PW)$.
 Moreover, let
$N_{C^+}$ be the number of  eigenvalues $\alpha\in \mathbb C$ with $\mathrm{Im}(\alpha)>0$.

As pointed out above, the eigenvalues counted by $N^-_n$ and $N_{C^+}$ correspond to cases of instabilities for the linearized problem \eqref{syst_lin_sw}. Note that to prove the spectral stability, it is enough to show that the generalized eigenproblem \eqref{Gene} does not have negative eigenvalues and $N_{C^+}=0$. Indeed, as a consequence of Propositions~\ref{equiv1} and~\ref{equiv2} and Lemma \ref{basic}, if $\alpha\in \C\smallsetminus\R$ is an eigenvalue of \eqref{Gene}, then $\bar\alpha$ is an eigenvalue too. Hence, if $N_{C^+}=0$, then the generalized eigenproblem \eqref{Gene} does not have solutions in $\C\smallsetminus \R$.

Finally, for using the counting argument introduce by Chugunova and Pelinovsky in \cite{ChouPel}, we need the following information on the essential spectrum of $\mathbb A$,
see \cite[Lemma~2-(H1') and Lemma~4]{PelProc}.

 \begin{lemma} Let 
$\mathscr M=-\mathscr J\mathscr L\mathscr J$ be defined on $\mathscr V$. 
 Then $\sigma_{\mathrm{ess}}(\mathscr M)=[0,+\infty)$. 
 Let  $\mathbb A=\mathscr P\mathscr M\mathscr P$ and $\mathbb K=\mathscr P\mathscr L^{-1}\mathscr P$ be defined on $\mathscr H$. Then $\sigma_{\mathrm{ess}}(\mathbb A)=[0,+\infty)$ and we can find $\delta_*, d_*>0$ such that   
  for any real number $0<\delta <\delta_*$, 
 $\mathbb A+\delta \mathbb K$ admits a bounded inverse and 
  we have
  $\sigma_{\mathrm{ess}}(\mathbb A+\delta \mathbb K)\subset [d_*\delta,+\infty)$.
\end{lemma}

\begin{Proof}
We check that 
\[
\mathscr J\mathscr L\mathscr JX=
\begin{pmatrix}
\ds\frac{\Delta_x q}{2} - k\cdot\nabla_x p
\\
\ds\frac{\Delta_x p}{2} +k\cdot\nabla_x q +2c\gamma \sigma_1\star \ds\int (-\Delta_z)^{-1/2}\sigma_2 (-\Delta_z)^{1/2}\pi
\ud z
\\
2c^2\Delta_z  \phi
\\
 2c^2\ds\frac{\Delta_z \pi}{2} + 2c\gamma \sigma_2\sigma_1\star p
\end{pmatrix}.\]
As a matter of fact, for any $\phi\in H^2(\mathbb R^n)$, the vector $X_e=(0,0,\phi,0)$ lies in $(\mathrm{Ker}(\mathscr L))^\perp$ and satisfies
\[
\mathscr J\mathscr L\mathscr JX_e=
\begin{pmatrix}
0
\\
0
\\
2c^2\Delta_z  \phi
\\
0
\end{pmatrix} \in (\mathrm{Ker}(\mathscr L))^\perp.
\]
Consequently $\mathscr MX_e=\mathbb A X_e=-\mathscr J\mathscr L\mathscr JX_e=(0,0,-2c^2\Delta_z  \phi,0)$.
It indicates that a Weyl sequence for $\mathbb A-\lambda \mathbb I$, $\lambda>0$, can be obtained by
adapting a Weyl sequence for $(-\Delta_z-\mu\mathbb I)$, $\mu>0$.
Let us consider a sequence of smooth functions $\zeta_\nu\in C^\infty_c(\mathbb R^n)$
such that $\mathrm{supp}(\zeta_\nu)\subset B(0, \nu+1)$, 
$\zeta_\nu(z)=1$ for $x\in B(0,\nu)$ and
$\|\nabla_z\zeta_\nu\|_{L^\infty(\mathbb R^n)}\leq C_0<\infty$, $\|D_z^2\zeta_\nu\|_{L^\infty(\mathbb R^n)}\leq C_0<\infty$, uniformly with respect to $\nu\in \mathbb N$.
We set $\phi_\nu(z)=\zeta_\nu(z)e^{i\xi\cdot z/(\sqrt 2 c)}$ for some $\xi\in \mathbb R^n$.
We get $$(-|\xi|^2 - 2c^2\Delta_z)\phi_\nu(z) =
{\color{black} -} e^{i\xi\cdot z/(\sqrt 2 c)}
\Big(\ds\frac{2i}{\sqrt 2 c}\xi\cdot \nabla_z\zeta_\nu+2c^2\Delta_z\zeta_\nu\Big)(z),
$$
which is thus bounded in $L^\infty(\mathbb R^n)$ and supported in $B(0,\nu+1)\smallsetminus B(0,\nu)$.
It follows that $\|(-|\xi|^2- 2c^2\Delta_z)\phi_\nu\|^2_{L^2(\mathbb R^n)}\lesssim \nu^{n-1}$, while 
$\|\phi_\nu\|^2_{L^2(\mathbb R^n)}\gtrsim \nu^n$. Accordingly, we obtain
$\frac{\|\phi_\nu\|^2_{L^2(\mathbb R^n)}}{\|(-|\xi|^2- 2c^2\Delta_z)\phi_\nu\|^2_{L^2(\mathbb R^n)}}
\gtrsim \nu\to \infty$ as $\nu\to \infty$.
Therefore, $\phi_\nu$ equally provides a Weyl sequence for 
$\mathscr M -|\xi|^2\mathbb I$ and 
$\mathbb A-|\xi|^2\mathbb I$, showing the inclusions
$ [0,\infty)\subset \sigma_{\mathrm{ess}}(\mathscr M)$ and 
$ [0,\infty)\subset \sigma_{\mathrm{ess}}(\mathbb A)$.

Next, let  $\lambda\notin [0,\infty)$. We suppose that we can find a Weyl sequence $(X_\nu)_{\nu\in \mathbb N}$ for $\mathscr M$, such that 
\[\begin{array}{lll}
\mathscr MX_\nu-\lambda X_\nu&=&
\begin{pmatrix}
-\lambda q_\nu -\ds\frac{\Delta_x q_\nu}{2} + k\cdot\nabla_x p_\nu
\\
-\lambda p_\nu -\ds\frac{\Delta_x p_\nu}{2} -k\cdot\nabla_x q_\nu - 2c\gamma \sigma_1\star \ds\int (-\Delta_z)^{-1/2}\sigma_2 (-\Delta_z)^{1/2}\pi_\nu
\ud z
\\
-\lambda \phi_\nu - 2c^2\Delta_z  \phi_\nu
\\
-\lambda \pi_\nu -  2c^2\Delta_z \pi_\nu - 2c\gamma \sigma_2\sigma_1\star p_\nu
\end{pmatrix}
\\
&=&
\begin{pmatrix}
q'_\nu\\p'_\nu\\\phi'_\nu\\\pi'_\nu\end{pmatrix}\xrightarrow [\nu\to \infty]{}0,\end{array}\]
with, moreover,  $\|X_\nu\|=1$ and  $X_\nu\rightharpoonup 0$ weakly in $\mathscr V$.
In particular, we can set \begin{equation}\label{hatphinu}
\widehat {\phi_\nu}(x,\xi)=\frac{\widehat {\phi'_\nu}(x,\xi)}{2c^2|\xi|^2- \lambda}.\end{equation}
It defines a sequence which tends to 0 strongly  $L^2(\mathbb T^d\times \mathbb R^n)$ since, writing $\lambda=a+ib\in \mathbb C\smallsetminus [0,\infty)$,  we get
$ |2c^2|\xi|^2- \lambda|^2
=|2c^2|\xi|^2- a|^2+b^2 $ which  is $\geq b^2>0$ when $\lambda \notin \mathbb R$, and, in case $b=0$, $\geq a^2>0$.
{\color{black}
Similarly, we can write 
\begin{equation}\label{hatpinu}
\widehat {\pi_\nu}(x,\xi)=\underbrace{\ds\frac{\widehat {\pi'_\nu}(x,\xi)}{2c^2|\xi|^2- \lambda}}_{=h_\nu(x,\xi)\in 
L^2(\mathbb T^d\times \mathbb R^n)}+
  \underbrace{\ds\frac{2c\gamma\widehat {\sigma_2}(\xi)}{2c^2|\xi|^2- \lambda}}_{\in L^2(\mathbb R^n)}\sigma_1\star p_\nu,\end{equation}}
where $h_\nu$ tends to 0 strongly  $L^2(\mathbb T^d\times \mathbb R^n)$.
We are led to the system
\begin{equation}\label{sysnu}\begin{array}{l}
\begin{pmatrix}
-\Big(\lambda  +\ds\frac{\Delta_x}{2} \Big)q_\nu + k\cdot \nabla_x p_\nu
\\
-k\cdot \nabla_x q_\nu    -\Big(\lambda  + \ds\frac{\Delta_x}{2}\Big)p_\nu -  4c^2 \gamma^2 
\ds\int 
\ds\frac{|\widehat {\sigma_2}|^2  }{(2\uppi)^n( 2c^2|\xi|^2-\lambda)}\ud \xi\times \Sigma\star p_\nu
\end{pmatrix}
\\[.4cm]
\hspace*{3cm}=
\begin{pmatrix}
q'_\nu
\\
p'_\nu -  2c \gamma\sigma_{1}\star\ds \int \ds\frac{\widehat {\sigma_2}(\xi)}{|\xi|} h_\nu(x,\xi)\frac{\ud \xi}
{(2\uppi)^n}
\end{pmatrix}
\xrightarrow[\nu\to \infty]{} 0.
\end{array}\end{equation}
Reasoning as in the proof of Proposition~\ref{SpecL}-1), we conclude that $X_\nu$ converges strongly to 0 in $\mathscr V$, a contradiction. 
Hence, $\lambda\in \mathbb C\smallsetminus [0,\infty)$ cannot belong to $\sigma_{\mathrm{ess}}(\mathscr M)$ and the identification $\sigma_{\mathrm{ess}}(\mathscr M)=[0,\infty)$ holds.
\\

Proposition~\ref{SpecL}-3) identifies $\mathrm{Ker}(\mathscr L)$.
Let us introduce  the mapping
\[
\widetilde {\mathscr P}: \begin{pmatrix} q\\ p
\end{pmatrix}\in L^2(\mathbb T^d)\times L^2(\mathbb T^d) \longmapsto 
\begin{pmatrix}
\ds\sum_{m\in \mathscr K_*,\ k\cdot m>0} (q_m-ip_m) e^{im\cdot x}
+\ds\sum_{m\in \mathscr K_*,\ k\cdot m<0} (q_m+ip_m) e^{im\cdot x}
 \\ p_0
 +i \ds\sum_{m\in \mathscr K_*,\ k\cdot m>0} (q_m-ip_m)
e^{im\cdot x}  
 -i \ds\sum_{m\in \mathscr K_*,\ k\cdot m<0} 
(q_m+ip_m)e^{im\cdot x}
\end{pmatrix}.
\]
Then, 
$$X=\begin{pmatrix} q\\ p\\ \phi\\\pi
\end{pmatrix}\longmapsto 
\begin{pmatrix} \widetilde{\mathscr P}\begin{pmatrix}q\\p\end{pmatrix}
\\0 \\ 0
\end{pmatrix}
$$
is the projection of $\mathscr V$ on $\mathrm{Ker}(\mathscr L)$.
Accordingly, we realize that 
 $\mathscr P$ does not modify the last two components of a vector $X=(q,p,\phi,\pi)\in \mathscr V$,
 and for $X\in (\mathrm{Ker}(\mathscr L))^\perp$, we have
 $p_0=0$, and $q_m=\pm i p_m$ for any $m\in\mathscr K_*$, depending on the sign of $k\cdot m$. 
 
Now, let $\lambda\in \mathbb C\smallsetminus[0,\infty)$ and suppose that we can exhibit  a Weyl sequence $(X_\nu)_{\nu\in\mathbb N}$ for $\mathbb A-\lambda \mathbb I$:
$X_\nu\in \mathscr H\subset (\mathrm{Ker}(\mathscr L))^\perp$, 
$\mathscr PX_\nu=X_\nu$,
$\|X_\nu\|=1$, 
$X_\nu\rightharpoonup 0$ in $\mathscr V$ and $\lim_{\nu\to \infty}\|(\mathbb A-\lambda \mathbb I)X_\nu\|
=0$.
We can apply the same reasoning as before for the last two components
of $(\mathbb A-\lambda \mathbb I)X_\nu$; it leads to \eqref{hatphinu} and \eqref{hatpinu}, where, using $\lambda\notin [0,\infty)$,  
 $\phi_\nu$ and $h_\nu$ converge strongly to 0 in $L^2(\mathbb T^d\times\mathbb R^n)$. 
We arrive at the following analog
to \eqref{sysnu}
\begin{equation}\label{sysnu2}\begin{array}{l}
(\mathbb I-\widetilde{\mathscr P})\begin{pmatrix}
-\Big(\lambda  +\ds\frac{\Delta_x}{2} \Big)q_\nu + k\cdot \nabla_x p_\nu
\\
-k\cdot \nabla_x q_\nu    -\Big(\lambda  + \ds\frac{\Delta_x}{2}\Big)p_\nu -  4c^2\gamma^2 
\ds\int 
\ds\frac{|\widehat {\sigma_2}|^2  }{(2\uppi)^n( 2c^2|\xi|^2-\lambda)}\ud \xi\times \Sigma\star p_\nu
\end{pmatrix}
\\[.4cm]
\hspace*{3cm}=
\begin{pmatrix}
q'_\nu
\\
p'_\nu\end{pmatrix}-
(\mathbb I-\widetilde{\mathscr P})
\begin{pmatrix}
0
\\
2c \gamma\sigma_{1}\star\ds \int \ds\frac{\widehat {\sigma_2}(\xi)}{|\xi|} h_\nu(x,\xi)\frac{\ud \xi}
{(2\uppi)^n}
\end{pmatrix}
\xrightarrow[\nu\to \infty]{} 0.
\end{array}\end{equation}
In order to derive from \eqref{sysnu2} an estimate in a positive Sobolev space as we did  in the proof of Proposition~\ref{SpecL}-1),
we should consider the 
Fourier coefficients arising from 
$-\frac{1}{2}\Delta_xq_\nu + k\cdot \nabla_x p_\nu$ and 
$-\frac{1}{2}\Delta_xp_\nu  -k\cdot \nabla_x q_\nu$, namely
$Q_m=\frac{m^2}{2}q_{\nu,m} + ik\cdot m p_{\nu,m}$ 
and 
$P_m=\frac{m^2}{2}p_{\nu,m} - ik\cdot m q_{\nu,m}$.
Only the coefficients belonging to $\mathscr K_*$ are affected by the action of $\widetilde{\mathscr P}$, which leads to 
$Q_m-(Q_m\mp iP_m)=\pm iP_m$ 
and $P_m\mp i(Q_m\mp iP_m)=\mp iQ_m$, according to the sign of $k\cdot m$.
However, we bear in mind that $q_m=\pm i p_m$ when $m\in 
\mathscr K_*$ with $\pm k\cdot m>0$.
Hence, for coefficients in $\mathscr K_*$, the contributions of the differential operators
reduces to $\pm i m^2 p_m =\pm m^2 q_m$ and $\mp i m^2 q_m=\pm m^2 p_m$, respectively.
Note also that for these coefficients there is no contributions coming from the convolution with $\sigma_1$ in \eqref{sysnu2} since $\sigma_{1,m}=0$ for $m\in \mathscr K_*$.
Therefore, reasoning as in the proof of Proposition~\ref{SpecL}-1) for coefficients $m\in \mathbb Z^d\smallsetminus \mathscr K_*$, 
we can obtain a uniform bound on $\sum_{m\in \mathbb Z^d} 
|m|^4(|q_{\nu,m}|^2+|p_{\nu,m}|^2)$, which  provides a uniform $H^2$ bound on $q_\nu$ and $p_\nu$, leading eventually to a contradiction.
We conclude that $\sigma_{\mathrm{ess}}(\mathbb A)=[0,\infty)$.
\\

Let $\delta>0$ and consider the shifted operator $\mathbb A+\delta \mathbb K$.
As a consequence of Lemma~\ref{l:nmn}, we will see that
 $\mathrm{Ker}(\mathbb A+\delta \mathbb K)=\{0\}$ for any $\delta>0$: 0 is not an eigenvalue for $\mathbb A+\delta \mathbb K$; let us  justify it does not belong to the  essential spectrum neither. To this end, we need to 
 detail the expression of the operator $\mathbb K$. Given $X\in \mathscr H$, we wish to find $X'\in \mathscr H$ satisfying
\[
\mathscr LX'=
  \begin{pmatrix}
 -  \ds\frac{ 1}{2}\Delta_x q'+    k\cdot \nabla_x p' +
 \gamma \sigma_1\star\left(
  \ds\int_{\mathbb R^n}(-\Delta)^{-1/2}\sigma_2\phi' \ud z\right) 
\\
-\ds\dfrac12\Delta_x p'-k\cdot \nabla_x q'
 \\
\ds\frac12 \phi'
+\gamma    
(-\Delta)^{-1/2}\sigma_2\sigma_1\star   q'\\
 \ds\frac{ \pi'}{2}
\end{pmatrix}
=
X.
\]
We infer  $\pi'=2\pi$ and the relation $\phi'=2\phi -2\gamma(-\Delta_z)^{-1/2}\sigma_2 \sigma_1\star q'$.
In turn, the Fourier coefficients of $q',p'$ are required to satisfy
\[
\begin{pmatrix}
m^2/2-2\gamma^2\kappa(2\uppi)^{2d}|\sigma_{1,m}|^2 & ik\cdot m
\\
-ik\cdot m & m^2/2
\end{pmatrix}
\begin{pmatrix}q'_m
\\
p'_m
\end{pmatrix}
=
\begin{pmatrix}
q_m -2\gamma(2\uppi)^d\sigma_{1,m}\ds\int(-\Delta)^{-1/2}\sigma_2\phi_m\ud z
\\
p_m
\end{pmatrix}.
\]
When $m\neq 0$, $m\notin \mathscr K_*$, the matrix of this system 
has  its determinant equal to $$
\mathrm{det}=\frac{m^4}{4}\big(1-4\gamma^2 \kappa(2\uppi)^{2d}\frac{|\sigma_{1,m}|^2}{m^2}\big)-(k\cdot m)^2.$$ Owing to \eqref{small}, since $(k\cdot m)^2$ takes values in $\mathbb N$, it does not vanish
and we obtain $q'_m,p'_m$ by solving the system
\[\begin{array}{l}
q'_m
=\ds\frac{1}{\mathrm{det}}\left(\ds\frac{m^2}{2} 
\Big(q_m -2\gamma(2\uppi)^d\sigma_{1,m}\ds\int(-\Delta)^{-1/2}\sigma_2\phi_m\ud z\Big)
-ik\cdot m p_m
\right),
\\
p'_m
=\ds\frac{1}{\mathrm{det}}\left(\ds+ik\cdot m 
\Big(q_m -2\gamma(2\uppi)^d\sigma_{1,m}\ds\int(-\Delta)^{-1/2}\sigma_2\phi_m\ud z\Big)
+\Big(\ds\frac{m^2}{2}-2\gamma^2\kappa(2\uppi)^{2d}|\sigma_{1,m}|^2\Big) p_m
\right).
\end{array}\]
If $m\in \mathscr K_*$ we find a solution in $(\mathrm{Ker}(\mathscr L))^\perp$ by setting
$p'_m=\frac{p_m}{m^2}$, $q'_m=\pm i p'_m$, according to the sign of $k\cdot m$;
if $m=0$, we set $p'_0=0$ and $q_0'=\frac{1}{2\gamma^2\kappa(2\uppi)^{2d}|\sigma_{1,0}|^2}
\big(q_0 -2\gamma(2\uppi)^d\sigma_{1,0}\int(-\Delta)^{-1/2}\sigma_2\phi_0\ud z\big)$.
This defines $X'=\mathbb KX$.

Therefore, the last two components of $(\mathbb A+\delta \mathbb K-\lambda \mathbb I)X$ read
\[
\begin{array}{l}
(2\delta-\lambda ) \phi -2c^2 \Delta_z \phi -2\delta \gamma (-\Delta)^{-1/2}\sigma_2 \sigma_1\star q',
\\
%\Big(\ds\frac{\delta}{2c^2}-\lambda \Big)
(2\delta -\lambda) \pi -\ds\frac12\Delta_z \pi
- \gamma\sigma_2 \sigma_1\star p'.
\end{array}\]
Hence, when $\lambda$ does not belong to  $[2\delta,\infty)$, 
% d_*,\infty)$, 
%with $d_*=\min(2,\frac{1}{2c^2})$, 
we can repeat the analysis performed above
to establish that $\lambda \notin \sigma_{\mathrm{ess}}(\mathbb A+\delta \mathbb K)$.
In particular the essential spectrum  of $\mathbb A$ has been shifted away from 0.
\end{Proof}

We are now able to apply the results of Chugunova and Pelinovsky \cite{ChouPel} (see also \cite{PelProc}), to obtain the following. 

\begin{theo}\cite[Theorem~1]{ChouPel}\label{ThCP}
Let $\mathscr L$ be defined by \eqref{lin_self_adj_opSW}.  Suppose \eqref{small}.
With the 
operators $\mathscr M, \mathbb A, \mathbb K$ defined as in Propositions~\ref{equiv1}-\ref{equiv2},
the following identity holds
\[N^-_n+N^0_n+N^+_n+N_{C^+}=n(\mathscr L).\]
\end{theo}

\noindent
Let us now detail the proof of Proposition~\ref{equiv1} and \ref{equiv2}, adapted from \cite[Prop.~1 \& Prop.~3]{PelProc}.
\\

\begin{ProofOf}{Propositions~\ref{equiv1} and ~\ref{equiv2}}
The goal is to establish connections between the following three problems:
\begin{itemize}
\item  [\namedlabel{Ev}{\textbf{(Ev)}}] the eigenvalue  problem $\mathbb LX=\lambda X$, with $\mathbb L=\mathscr J\mathscr L$,
\item  [\namedlabel{Co}{\textbf{(Co)}}] the coupled problem $\mathscr L X=\lambda Y$, $\mathscr M Y =-\lambda X$, with $\mathscr M=-\mathscr J\mathscr L\mathscr J$,
\item [\namedlabel{GEv}{\textbf{(GEv)}}] the generalized eigenvalue  problem $\mathbb A W=\alpha \mathbb KW$, with $\mathbb A=\mathscr P \mathscr M\mathscr P$, $\mathbb K=\mathscr P \mathscr L^{-1}\mathscr P$,  the projection $\mathscr P$ on $(\mathrm{Ker}(\mathscr L))^\perp$, and $W\in \mathscr H=D(\mathscr M)\cap  (\mathrm{Ker}(\mathscr L))^\perp$.
\end{itemize}
The proof of Propositions~\ref{equiv1} and ~\ref{equiv2} follows from the following sequence
of arguments.
\begin{enumerate}[label=(\roman*)]
\item By Lemma~\ref{basic}, we already know that if there exists a solution $(\lambda,X_+)$ of \ref{Ev}, with $\lambda \neq 0$ and $X_+\neq 0$, then, there exists $X_-\neq 0$, such that  $(-\lambda,X_-)$ satisfies \ref{Ev}.
Being eigenvectors associated to distinct eigenvalues, $X_+$ and $X_-$ are linearly independent.
Note that only this part of the proof uses the specific structure of the operator $\mathbb L$.
\item From these eigenpairs for $\mathbb L$, we set 
\[X=\ds\frac{X_++X_-}{2}, \qquad
Y=\widetilde{\mathscr J}\left(\ds\frac{X_+-X_-}{2}\right).\]
Since  $X_+$ and $X_-$ are linearly independent, we have $X\neq 0$, $Y\neq 0$.
Moreover, $X=\frac{X_++X_-}{2}$ and $\mathscr JY=\frac{X_+-X_-}{2}$ are linearly independent.
We get 
\[\begin{array}{l}
\mathscr L X=\widetilde{\mathscr J}\mathbb LX=\widetilde{\mathscr J}\left(\ds\frac{\lambda}{2}(X_+-X_-)\right)=\lambda Y,
\\
\mathscr MY=-\mathscr J\mathscr L \left(\ds\frac{X_+-X_-}{2}\right)=-\mathbb L\left(\ds\frac{X_+-X_-}{2}\right)
=-\ds\frac{\lambda}{2}(X_++X_-)=-\lambda X,
\end{array}\]
so that $(\lambda, X, Y)$ satisfies \ref{Co}.
\item If  $(\lambda, X, Y)$ is a solution \ref{Co}, then  $(-\lambda, X, -Y)$ satisfies \ref{Co} too.
\item  Let $(\lambda, X, Y)$ be a solution \ref{Co}. Set 
\[
X'=\mathscr JY,\qquad Y'=\widetilde{\mathscr J}X.\]
We observe that 
\[\begin{array}{l}
\mathscr M Y'=-\mathscr J\mathscr L \mathscr J \widetilde{\mathscr J} X
=-\mathscr J\mathscr L X=-\mathscr J(\lambda Y)=-\lambda X',
\\
\mathscr L X'=\mathscr L \mathscr JY
=\widetilde{\mathscr J}\mathscr J\mathscr L \mathscr JY 
=-\widetilde{\mathscr J}\mathscr MY=\lambda \widetilde{\mathscr J}X=\lambda Y',
\end{array}\]
which means that $(\lambda, \mathscr JY,\widetilde{\mathscr J}X)$ is a solution of \ref{Co}.
Moreover,  if $X$ and $ \mathscr JY$ are linearly independent,  $Y$ and $\widetilde{\mathscr J}X$ are linearly independent too.
\item Let $(\lambda, X, Y)$ be a solution \ref{Co}, with $X\neq 0$.
We get
\[\begin{array}{lll}
\mathbb L(X\pm \mathscr J Y)&=&
\mathscr J\mathscr L X\pm \mathscr J  \mathscr L\mathscr  JY
=\mathscr J\mathscr LX\mp \mathscr MY
\\
&=&\mathscr J(\lambda Y)\pm \lambda X=\pm \lambda (X\pm \mathscr  JY),
\end{array}\]
so that $(\pm \lambda , X\pm \mathscr JY)$ satisfy \ref{Ev}.
In the situation where  $X$ and $\mathscr JY$ are linearly independent, 
we have $X\pm\mathscr JY\neq 0$ and 
$(\pm \lambda,X\pm\mathscr JY)$ are eigenpairs for $\mathbb L$.
Otherwise, one of the vectors 
$X\pm\mathscr JY$ might vanish.
Nevertheless, since only one of these two vectors can be 0, 
we still obtain
an eigenvector $X_\pm\neq 0$ of $\mathbb L$, associated to either $\pm \lambda$.
Coming back to i), we conclude that $\mp\lambda$ is an eigenvalue too.
\end{enumerate}

Items i) to v) justify the equivalence stated in Proposition~\ref{equiv1}.

\begin{enumerate}[label=(\roman*),resume]
\item Let $(\lambda, X, Y)$ be a solution \ref{Co}.
From $\mathscr LX=\lambda Y$, we infer $Y\in \mathrm{Ran}(\mathscr L)\subset
(\mathrm{Ker}(\mathscr L))^\perp$ so that $\mathscr PY=Y$.
The relation thus recasts as
\[
X=\lambda \mathscr P\mathscr L^{-1}\mathscr P Y+\tilde Y,\qquad
\tilde Y\in \mathrm{Ker}(\mathscr L),\qquad \mathscr P\tilde Y=0.\]
(Here, $\mathscr P\mathscr L^{-1}\mathscr PY$ stands for the unique solution of $\mathscr L Z=Y$ which lies  in 
$(\mathrm{Ker}(\mathscr L))^\perp$.)
We obtain 
\[\begin{array}{lll}
\mathscr P\mathscr MY&=&\mathscr P\mathscr (-\lambda X)=
-\lambda \mathscr P(\lambda \mathscr P\mathscr L^{-1}\mathscr P Y+\tilde Y)
\\&=&-\lambda^2  \mathscr P\mathscr L^{-1}\mathscr P Y=-\lambda ^2 \mathbb KY
=\mathscr P\mathscr M\mathscr PY=\mathbb AY,\end{array}\]
so that $(-\lambda^2,Y)$ satisfies \ref{GEv}.
Going back to iv), we  know that $(-\lambda^2,\widetilde {\mathscr J}X)$ is equally a solution to \ref{GEv}.
If $X$ and $\mathscr JY $ are linearly independent, we obtain this way  two linearly independent
vectors, $Y$ and $\widetilde {\mathscr J}X$, solutions of \ref{GEv} with $\alpha=-\lambda^2$.

\item Let $(\alpha,W)$ satisfy \ref{GEv}, with $\alpha\neq 0$, $W\neq 0$.
We set 
$X=\frac{-\mathscr M W}{\sqrt{-\alpha}}$.
We have 
$$\widetilde {\mathscr J}X=-\ds\frac{1}{\sqrt{-\alpha}}\widetilde {\mathscr J}\mathscr M W
=\ds\frac{1}{\sqrt{-\alpha}}\widetilde {\mathscr J}\mathscr J\mathscr L\mathscr J W
=\ds\frac{1}{\sqrt{-\alpha}}\mathscr L\mathscr J W$$
which lies in $\mathrm{Ran}(\mathscr L)\subset  (\mathrm{Ker}(\mathscr L))^\perp$.
Thus, using $\mathscr P 
\widetilde {\mathscr J}X= \widetilde {\mathscr J}X$, we compute 
\[
\mathbb K\widetilde {\mathscr J}X=\mathscr P\mathscr L^{-1}\mathscr P 
\widetilde {\mathscr J}X= 
\mathscr P\mathscr L^{-1}
\widetilde {\mathscr J}X=
\ds\frac{1}{\sqrt{-\alpha}}\mathscr P\mathscr L^{-1}\mathscr L\mathscr J W
=\ds\frac{1}{\sqrt{-\alpha}}\mathscr P\mathscr J W.
\]
Next, we observe that
\[
\mathbb A\widetilde {\mathscr J}X=
\mathscr P\mathscr M\mathscr P 
\widetilde {\mathscr J}X
=-
\mathscr P\mathscr J\mathscr L\mathscr J 
\widetilde {\mathscr J}X=-
\mathscr P\mathscr J\mathscr L X
=
\ds\frac{1}{\sqrt{-\alpha}} \mathscr P\mathscr J\mathscr L  \mathscr M W
.
\]
However, we can use $\mathscr P W=W$ (since $W\in \mathscr H\subset (\mathrm{Ker}(\mathscr L))^\perp$) and the fact that, for any vector $Z$,
$\mathscr L Z=\mathscr L (\mathbb I-\mathscr P)Z+\mathscr L\mathscr PZ
=0+\mathscr L\mathscr PZ$, which yields
\[\begin{array}{lll}
\mathbb A\widetilde {\mathscr J}X
&=&
\ds\frac{1}{\sqrt{-\alpha}} \mathscr P\mathscr J\mathscr L \mathscr P \mathscr M \mathscr PW
=\ds\frac{1}{\sqrt{-\alpha}} \mathscr P\mathscr J\mathscr L\mathbb AW
=
-\sqrt{-\alpha} \mathscr P\mathscr J\mathscr L\mathbb K W
\\&=&-\sqrt{-\alpha} \mathscr P\mathscr J\mathscr L\mathscr P\mathscr L^{-1}\mathscr P W
=-\sqrt{-\alpha} \mathscr P\mathscr J\mathscr L\mathscr L^{-1} W
=-\sqrt{-\alpha} \mathscr P\mathscr JW.
\end{array}\]
We conclude that $\mathbb A\widetilde {\mathscr J}X=\alpha\mathbb K \widetilde {\mathscr J}X$:
$(\alpha,\widetilde {\mathscr J} X)$ satisfies \ref{GEv}. 
\item 
Let $(\alpha,W)$ satisfy \ref{GEv}, with $\alpha\neq 0$, $W\neq 0$.
We have 
\[\mathscr P(\mathscr M \mathscr PW-\alpha \mathscr L^{-1}\mathscr PW)=0
\] and thus 
\[\mathscr M \mathscr PW-\alpha \mathscr L^{-1}\mathscr PW=\tilde Y\in \mathrm{Ker}(\mathscr L).\]
Let us set
\[
Y=\mathscr PW\in (\mathrm{Ker}(\mathscr L))^\perp,\qquad
X=-\ds\frac{\mathscr M \mathscr PW}{\sqrt{-\alpha}}=
\ds\frac {-1}{\sqrt{-\alpha}}(\tilde Y +\alpha \mathscr L^{-1}\mathscr PW),\]
so that
\[
\mathscr LX=\sqrt{-\alpha}\mathscr PW=\sqrt{-\alpha}Y,\qquad
\mathscr M Y=\mathscr M\mathscr P W=-\sqrt{-\alpha}X.\]
{\color{black} (Incidentally, since $W$ is assumed to belong to $\mathscr H$, we have $W=\mathscr PW=Y$.)}
Therefore $(\sqrt{-\alpha},X,Y)$ satisfies \ref{Co}.
By v), $(\pm \sqrt{-\alpha},X\pm \mathscr JY)$ satisfy \ref{Ev}, and at least one of the vectors $
X\pm \mathscr JY$ does not vanish; using i), we thus obtain 
eigenpairs $(\pm \sqrt{-\alpha},X_\pm)$ of $\mathbb L$. With ii), we construct 
solutions of \ref{Co} under the form $\big( \sqrt{-\alpha},\frac{X_++X_-}{2},
\widetilde {\mathscr J}\big(\frac{X_+-X_-}{2}\big)\big)$, which, owing to iv) and vi), 
provide the linearly independent solutions $\big(\alpha,\widetilde {\mathscr J}\big(\frac{X_+\pm X_-}{2}\big)\big)$ of \ref{GEv}.
The dimension of the linear space 
of solutions of \ref{GEv} is at least 2.

At least one of these vectors $X_\pm$ is given by the formula
$$\tilde X_\pm=-\ds\frac{\mathscr MW}{\sqrt{-\alpha}}\pm \mathscr JW.$$
By the way, we indeed note that
$\mathbb A W=\alpha\mathbb K W$, with $W\in \mathscr H$, 
can be cast as 
$\mathscr L \mathscr J\mathscr L\mathscr J W=-\alpha W$
{\color{black} since it means
$$
(\mathbb A -\alpha\mathbb K )W=\mathscr P(\mathscr M-\alpha \mathscr L^{-1})
\underbrace{\mathscr P W}_{=W\in \mathscr H}=
-\mathscr P(\mathscr J\mathscr L\mathscr J+\alpha \mathscr L^{-1})W=0
$$
so that  $(\mathscr J\mathscr L\mathscr J+\alpha \mathscr L^{-1})W\in \mathrm{Ker}(\mathscr L)$.
It follows that}
% (see Lemma~\ref{idKernel} below) so that
\[\begin{array}{lll}
\mathbb L\Big(-\ds\frac{\mathscr MW}{\sqrt{-\alpha}}\pm \mathscr JW\Big)
&=&\ds\frac{1}{\sqrt{-\alpha}}\mathscr J(\mathscr L \mathscr J \mathscr L\mathscr J W)
\pm \mathscr J\mathscr L\mathscr J W
\\
&=&\sqrt{-\alpha} \mathscr J W \mp \mathscr M W
=\pm \sqrt{-\alpha} \Big(-\ds\frac{\mathscr MW}{\sqrt{-\alpha}}\pm \mathscr JW\Big).
\end{array}\]
 With these manipulations we have checked  that $(\pm  \sqrt{-\alpha},\tilde X_\pm)
 $ satisfy \ref{Ev}.
 If  both vectors $\tilde X_\pm $ are non zero, we get $X_\pm=\tilde X_\pm$
 and we recover $W=\widetilde {\mathscr J}\big(\frac{X_+-X_-}{2}\big)$.
 If $\tilde X_\pm=0$, then, we get $\tilde X_\mp=\mp \mathscr JW\neq 0$,
 and we directly obtain $X_\mp=\tilde X_\mp$, $W=\mp \widetilde {\mathscr J} X_\mp$.
 In any cases, $W$ lies in the space spanned by $X_+$ and $X_-$  and 
 the dimension of the space of solutions of \ref{GEv} is even.

\end{enumerate}

This ends the proof of Proposition~\ref{equiv1} and \ref{equiv2}.
\end{ProofOf}

\subsection{Spectral instability}

We are going to compute 
the terms arising in Theorem~\ref{ThCP}. Eventually, it will allow us 
to identify the possible unstable modes.
In what follows, we find convenient to work with the operator $\mathscr M-\alpha \mathscr L^{-1}$
instead of $\mathscr P(\mathscr M-\alpha \mathscr L^{-1})\mathscr P=\mathbb A-\alpha \mathbb K$, owing to to the following claim.

{\color{black} 
\begin{lemma}\label{idKernel}
Let $\alpha\neq 0$ and $X\in \mathscr H$. The following two problems are equivalent:
\begin{itemize}
\item [\textcircled{1}] $X\in \mathrm{Ker}(\mathbb A-\alpha \mathbb K)$,
\item [\textcircled{2}] there exists $\tilde X\in \mathscr V$ such that $\mathscr MX=\alpha \tilde X$ and $\mathscr L\tilde X=X$.
\end{itemize}
%In the space $\mathscr H=D(\mathscr M)\cap (\mathrm{Ker}(\mathscr L))^\perp$, the 
%two subspaces $$ and 
%$\mathrm{Ker}(\mathscr M-\alpha \mathscr  L^{-1}) $ coincide.
\end{lemma}

\begin{Proof}
Suppose \textcircled{1}.
Since $X=\mathscr PX\in \mathscr H$, it means $\mathscr P(\mathscr M-\alpha \mathscr L^{-1})X=0$, that is
$(\mathscr M-\alpha \mathscr L^{-1})X=Z\in \mathrm{Ker}( \mathscr L)$.
Since $\alpha\neq 0$, we can set $\tilde X=\frac{\mathscr MX}{\alpha}\in \mathscr V$.
It satisfies 
$\mathscr L\tilde X=\frac1\alpha\mathscr L(Z+\alpha \mathscr L^{-1}X)=X$, and \textcircled{2} holds.

Conversely, suppose \textcircled{2}.
We bear in mind that the pseudo-inverse $\mathscr L^{-1}$ is defined as an application from $(\mathrm {Ker}(\mathscr L))^\perp$ to itself, hence we can decompose $\tilde  X=\mathscr L^{-1}X+Z$, with $Z\in \mathrm {Ker}(\mathscr L)$.
Therefore, we get 
$\mathscr MX-\alpha \tilde X=(\mathscr M-\alpha \mathscr L^{-1})X-\alpha Z=0$.
In other words, 
$(\mathscr M-\alpha \mathscr L^{-1})X=\alpha Z\in \mathrm {Ker}(\mathscr L)$ which implies, since $X=\mathscr PX\in \mathscr H$,
$(\mathbb A-\alpha \mathbb K)X=\mathscr P(\mathscr M-\alpha \mathscr L^{-1})X=0$: \textcircled{1} is satisfied.
%Let  $X\in \mathscr H$ satisfy $ \mathscr MX=\alpha \mathscr  L^{-1}X$. Then, we have
%$X=\mathscr PX$ and, thus,  $
%(\mathbb A-\alpha \mathbb K)X=\mathscr P(\mathscr M-\alpha \mathscr L^{-1})\mathscr PX=
%\mathscr P(\mathscr MX-\alpha \mathscr L^{-1}X)=0$, showing the inclusion 
%$\mathrm{Ker}(\mathscr M-\alpha \mathscr  L^{-1})\cap \mathscr H \subset \mathrm{Ker}(\mathbb A-\alpha \mathbb K) $.
%
%Conversely,
%the equation $(\mathbb A-\alpha \mathbb K)X=0$, with $X=\mathscr PX\in (\mathrm{Ker}(\mathscr L))^\perp$ means that
%$(\mathscr M-\alpha \mathscr L^{-1})X=Y\in \mathrm{Ker}(\mathscr L)$.
%Applying $\mathscr L$ then yields
%$\mathscr L\mathscr MX=\alpha X$. Since both terms  of this relation lie in $(\mathrm{Ker}(\mathscr L))^\perp$, it is legitimate to apply $\mathscr L^{-1}$, showing that 
%$\mathscr MX=\alpha \mathscr L^{-1}X$: we have shown
%$ \mathrm{Ker}(\mathbb A-\alpha \mathbb K)\cap \mathscr H
%\subset 
%\mathrm{Ker}(\mathscr M-\alpha \mathscr  L^{-1})
%$.
\end{Proof}
}

\noindent
Therefore, we shall consider 
%the solutions of the generalized eigenvalue problem
%$\mathscr M X=\alpha \mathscr L^{-1}X$, with $X\in \mathscr H$. 
%We rewrite the equation by introducing an 
auxiliary problem: %unknown:
$$\mathscr M X=\alpha \tilde X,\qquad \mathscr L \tilde X=X.$$

\begin{lemma} Suppose \eqref{small}. $N^0_n=1$.\end{lemma}

\begin{Proof}
We are interested in the solutions of
\[\begin{array}{l}
-\ds\frac12\Delta_x q+k\cdot\nabla_x p=0,
\\[.3cm]
-\ds\frac12\Delta_x p-k\cdot\nabla_x q- 2c \gamma\sigma_1\star \ds\int \sigma_2\pi\ud z=0,
\\[.3cm]
-2c^2\Delta_z \phi=0, 
\\[.3cm]
- 2c^2\Delta_z\pi - 2c \gamma\sigma_2\sigma_1\star p=0.
\end{array}\]
We infer $\phi(x,z)=0$ and $\widehat
\pi(x,\xi)= \frac\gamma c\frac{\widehat{\sigma_2}(\xi)}{|\xi|^2}\sigma_1\star p(x)$, and, next,
\[
-\ds\frac12\Delta_x q+k\cdot\nabla_x p=0,\qquad
-\ds\frac12\Delta_x p-k\cdot\nabla_x q-2\gamma^2\kappa \Sigma\star p=0
\]
with $\Sigma=\sigma_1\star \sigma_1$.
In terms of Fourier coefficients, it becomes
\[
\ds\frac{m^2}2 q_m+ik\cdot m p_m=0,\qquad
\ds\frac{m^2}2 p_m-ik\cdot m q_m-2(2\uppi)^{2d}\gamma^2\kappa |\sigma_{1,m}|^2 p_m=0
.\]
For $m=0$, we get $p_0=0$ and
we find the eigenfunction $(\mathbf 1,0,0,0)=Y_0={\color{black} +\mathscr JX_0}
$ with $X_0=(0,\mathbf 1,0,0)
\in\mathrm{Ker}(\mathscr L)$.

For $m\neq 0$ with $\sigma_{1,m}\neq 0$, we get
$$
m^4
-4(k\cdot m)^2=\underbrace{2(2\uppi)^{2d}\gamma^2\kappa |\sigma_{1,m}|^2}_{\in (0,1)} m^2.
$$
which cannot hold (see the proof of Proposition~\ref{SpecL} for more details).

For $m\neq 0$ with $\sigma_{1,m}= 0$, we get
$M_m\begin{pmatrix}q_m \\ p_m\end{pmatrix}=0$
with $M_m$ defined in \eqref{defmatM}.
As far as
$m^4-4(k\cdot m)^2\neq 0$, $M_m$ is invertible and the only solution is $p_m=0=q_m$.
If $m^4-4(k\cdot m)^2=0$, we find the eigenfunctions $(e^{ik\cdot m},\pm ie^{ik\cdot m},0,0)$.
These functions belong to $\mathrm{Ker}(\mathscr L)$, and thus do not lie in the working space $\mathscr H$.
\\

We conclude that $\mathrm{Ker}(\mathscr M)=\mathrm{span}_{\mathbb R}\{Y_0\}$.
Moreover, this vector $Y_0$ does not belong to $\mathrm{Ran}(\mathscr M)$ so that 
the algebraic multiplicity of the eigenvalue 0 is 1.
Finally, bearing in mind \eqref{CalNo}, which can be recast as
$(\mathbb K Y_0 | Y_0)<0$,  we arrive at
$N^0_n=1$. 
\end{Proof}

\begin{lemma}\label{l:nmn}
  Suppose \eqref{small}. The generalized eigenproblem \eqref{Gene} does not admit negative eigenvalues. In particular, $N^-_n=0$.
\end{lemma}

\begin{Proof}
Let $\alpha<0$, $X=(q,p,\phi,\pi)$ and $\tilde X=(\tilde q,\tilde p,\tilde \phi,\tilde\pi)$ satisfy
\begin{equation}\label{eqal1}\begin{array}{l}
-\ds\frac12\Delta_x q +k\cdot\nabla_x p =\alpha \tilde q,
\\ [.3cm]
-\ds\frac12\Delta_x p -k\cdot\nabla_x q -2c\gamma\sigma_1\star\ds\int \sigma_2\pi\ud z=\alpha \tilde p,
\\[.3cm]
-2c^2\Delta_z \phi=\alpha\tilde \phi,
\\[.3cm]
- 2c^2\Delta_z\pi-2c\gamma\sigma_2\sigma_1\star p=\alpha\tilde\pi,
\end{array}\end{equation}
where
\begin{equation}\label{eqal2}
\begin{array}{l}
q=-\ds\frac12\Delta_x \tilde q +k\cdot\nabla_x  \tilde p +\gamma\sigma_1\star\ds\int(-\Delta_z)^{-1/2}\sigma_2 
 \tilde \phi\ud z,
\\ [.3cm]
p=-\ds\frac12\Delta_x\tilde p -k\cdot\nabla_x \tilde q 
,
\\[.3cm]
\phi=\ds\frac12\tilde \phi+\gamma(-\Delta_z)^{-1/2}\sigma_2\sigma_1\star\tilde q,
\\[.3cm]
\pi= \ds\frac{\tilde\pi}{2}.
\end{array}\end{equation}
This leads to solve an elliptic equation for $\pi$
\[
\Big(\ds\frac{|\alpha|}{c^2}-\Delta_z\Big)\pi= \frac\gamma c\sigma_2\sigma_1\star p.\]
In other words, we get, by means of Fourier transform
\[
\widehat \pi(x,\xi)=\ds\frac{\gamma}{c}\sigma_1\star p(x)\times  \ds\frac{\widehat {\sigma_2}(\xi)}{|\xi|^2+|\alpha|/c^2}.\]
On the same token, we obtain 
\[
\Big(\ds\frac{|\alpha|}{c^2}-\Delta_z\Big)\tilde\phi=-2\gamma(-\Delta_z)^{1/2}\sigma_2\sigma_1\star\tilde q,\]
which yields
\[
\widehat {\tilde \phi}(x,\xi)=-2\gamma\sigma_1\star \tilde q(x)\times 
 \ds\frac{|\xi|\widehat{\sigma_2}(\xi)}{|\xi|^2+|\alpha|/c^2}.\]
For $\lambda>0$, 
we introduce the symbol
\[
0\leq \kappa_\lambda=\ds\int \ds\frac{|\widehat {\sigma_2}(\xi)|^2}{|\xi|^2+\lambda}\leq \kappa.\]
It turns out that 
\[\begin{array}{lll}
-\ds\frac12\Delta_x q +k\cdot\nabla_x p =\alpha \tilde q,
\\ [.3cm]
-\ds\frac12\Delta_x p -k\cdot\nabla_x q -2\gamma^2\kappa_{|\alpha|/c^2}\Sigma\star p=\alpha \tilde p,
\end{array}\]
with
\[\begin{array}{lll}
q=-\ds\frac12\Delta_x \tilde q +k\cdot\nabla_x  \tilde p -2\gamma^2\kappa_{|\alpha|/c^2}
\Sigma\star\tilde q,
\\ [.3cm]
p=-\ds\frac12\Delta_x\tilde p -k\cdot\nabla_x \tilde q. 
\end{array}\]
For the Fourier coefficients, it casts as 
\[\begin{array}{lll}
\ds\frac{m^2}{2} q_m +ik\cdot m p_m =\alpha \tilde q_m,
\\ [.3cm]
\ds\frac{m^2}{2} p_m  -ik\cdot m q_m -2\gamma^2\kappa_{|\alpha|/c^2}(2\uppi)^{2d}|\sigma_{1,m}|^2 p_m=\alpha \tilde p_m,
\end{array}\]
with
\[\begin{array}{lll}
q_m=\ds\frac{m^2}{2}  \tilde q_m +ik\cdot m  \tilde p_m -2\gamma^2\kappa_{|\alpha|/c^2}
(2\uppi)^{2d}|\sigma_{1,m}|^2 \tilde q_m,
\\ [.3cm]
p_m=\ds\frac{m^2}{2} \tilde p_m -ik\cdot m \tilde q_m. 
\end{array}\]
We are going to see that these equations do not have non trivial solutions with $\alpha<0$:
\begin{itemize}
\item 
If $m=0$, we get $p_0=0$, 
$\tilde q_0=0$, and, consequently, $\tilde p_0=0$, $q_0=0$.
Hence, for $\alpha<0$, we cannot find an eigenvector 
with a non trivial 0-mode. 
\item 
If $m\neq 0$ and $\sigma_{1,m}= 0$, we see that $(q_m,p_m)$ and $(\tilde q_m,\tilde p_m) $ 
are related by 
\begin{equation}\label{systil}
M_m\begin{pmatrix} q_m\\ p_m\end{pmatrix}=\alpha \begin{pmatrix} \tilde q_m\\ \tilde p_m\end{pmatrix},\qquad
\begin{pmatrix} q_m\\ p_m\end{pmatrix}=M_m\begin{pmatrix} \tilde q_m\\ \tilde p_m\end{pmatrix}.\end{equation}
It means that $\alpha$ is an eigenvalue of 
\[
M_m^2=\begin{pmatrix}
\frac{m^4}{4}+(k\cdot m)^2 & im^2 k\cdot m
\\
- im^2 k\cdot m& \frac{m^4}{4}+(k\cdot m)^2
\end{pmatrix}.
\]
The roots of the characteristic polynomial  of $M_m^2$ are $(\frac{m^2}{2}\pm k\cdot m)^2\geq 0$, which contradicts the assumption $\alpha<0$.
\item 
For the case where  $m\neq 0$ and $\sigma_{1,m}\neq 0$, we introduce the shorthand notation $a_m=2\gamma^2(2\uppi)^{2d}|\sigma_{1,m}|^2 \kappa_{|\alpha|/c^2}$, bearing in mind that $0<a_m<\frac{m^2}{2}$ by virtue of the smallness condition  \eqref{small}.
We are led to the systems
\[
\left(M_m-\begin{pmatrix} 0 & 0 
\\
0 & a_m\end{pmatrix}
\right)\begin{pmatrix} q_m\\ p_m\end{pmatrix}=\alpha \begin{pmatrix} \tilde q_m\\ \tilde p_m\end{pmatrix},\qquad
\begin{pmatrix} q_m\\ p_m\end{pmatrix}=\left(M_m
-\begin{pmatrix} a_m & 0 
\\
0 & 0\end{pmatrix}
\right)\begin{pmatrix} \tilde q_m\\ \tilde p_m\end{pmatrix},\]
which imply that $\alpha$ is an eigenvalue of the matrix
\[
\left(M_m-\begin{pmatrix} 0 & 0 
\\
0 & a_m\end{pmatrix}
\right)\left(M_m
-\begin{pmatrix} a_m & 0 
\\
0 & 0\end{pmatrix}
\right).
\]
However the eigenvalues of this matrix read
$\big(\sqrt{\frac{m^2}{2}(\frac{m^2}{2}-a_m)}\pm (k\cdot m)^2\big)^2\geq 0$, contradicting that $\alpha$ is negative.\end{itemize}
\end{Proof}

\begin{lemma}
Suppose \eqref{small}. $N^+_n=\#\{m\in \mathbb Z^d\smallsetminus\{0\}, \sigma_{1,m}=0,\textrm{ and } m^4-4(k\cdot m)^2<0\}$.\end{lemma}

\begin{Proof}
We should consider the system of equations \eqref{eqal1}-\eqref{eqal2}, now with $\alpha>0$.
For Fourier coefficients, it casts as 
\[\begin{array}{l}
\ds\frac{m^2}2 q_m +ik\cdot m p_m =\alpha \tilde q_m,
\\ [.3cm]
\ds\frac{m^2}2 p_m -ik\cdot m q_m - 2c\gamma(2\uppi)^d\sigma_{1,m}\ds\int \sigma_2\pi_m\ud z=\alpha \tilde p_m,
\\[.3cm]
-2c^2\Delta_z \phi_m=\alpha\tilde \phi_m,
\\[.3cm]
- 2c^2\Delta_z\pi_m-2c\gamma(2\uppi)^d\sigma_{1,m}\sigma_2 p_m=\alpha\tilde\pi_m,
\end{array}\]
where
\[
\begin{array}{l}
q_m=\ds\frac{m^2}2\tilde q_m +i k\cdot m  \tilde p_m +\gamma(2\uppi)^d\sigma_{1,m}\ds\int(-\Delta_z)^{-1/2}\sigma_2 
 \tilde \phi_m\ud z,
\\ [.3cm]
p_m=\ds\frac{m^2}2\tilde p_m -ik\cdot m \tilde q _m
,
\\[.3cm]
\phi_m=\ds\frac12\tilde \phi_m+\gamma(2\uppi)^d(-\Delta_z)^{-1/2}\sigma_2\sigma_{1,m}\tilde q_m,
\\[.3cm]
\pi_m= \ds\frac{\tilde\pi_m}{2}.
\end{array}\]
\begin{itemize}
\item For $m=0$, we obtain $p_0=0$, $\tilde q_0=0$.
Hence $\pi_0$ satisfies $(-\alpha/c^2-\Delta_z)\pi_0=0$. Here, $+\alpha/c^2$ lies in the essential spectrum of $-\Delta_z$ and the only solution in $L^2$ of this equation is $\pi_0=0$. In turn, this implies $\tilde p_0=0$, $(-\alpha/c^2-\Delta_z)\phi_0=0$, and thus 
$\phi_0=0$, $q_0=0$.
Hence, for $\alpha>0$, we cannot find an eigenvector 
with a non trivial 0-mode. 

\item When $m\neq 0$ and $\sigma_{1,m}=0$, we are led to 
$(-\alpha/c^2-\Delta)\phi_m=0$, $(-\alpha/c^2-\Delta)\pi_m=0$
that imply $\phi_m=0$, $\pi_m=0$.
In turn, we get 
\eqref{systil} for $q_m,p_m,\tilde q_m,\tilde p_m$. This holds iff $\alpha$ is an eigenvalue 
of $M_m^2$.
If $m^4\neq 4(k\cdot m)^2$, we find two  eigenvalues 
$\alpha_{m,\pm}=(\frac{m^2}{2}\pm k\cdot m)^2>0$, with associated eigenvectors 
$X_{m,\pm}=(e^{im\cdot x},\mp ie^{im\cdot x},0,0)$, respectively.
To decide whether these modes should be counted, we need to evaluate the sign of
$(\mathscr L^{-1}X_{m,\pm}|X_{m,\pm})$.
We start by solving $\mathscr LX'_{m,\pm}=X_{m,\pm}$. It yields
$\frac{\phi'_{m,\pm}}{2}=0$, $\frac{\pi'_{m,\pm}}{2}=0$ and 
\[M_m\begin{pmatrix}q'_{m,\pm}\\p'_{m,\pm}\end{pmatrix}
=\begin{pmatrix}1\\\mp i \end{pmatrix}.\]
We obtain
\[q'_{m,\pm}=
\ds\frac{2}{m^2\pm 2k\cdot m},\quad
\pi'_{m,\pm}=\ds\frac{\mp 2i}{m^2\pm 2k\cdot m},\]
so that 
\[\begin{array}{lll}
(\mathscr L^{-1}X_{m,\pm}|X_{m,\pm})
&=&\ds\frac{2}{m^2\pm 2k\cdot m}
\left(\ds\int_{\mathbb T^d}
e^{im\cdot x}e^{-im\cdot x}
\ud x
+\ds\int_{\mathbb T^d}
(\mp i) e^{im\cdot x}{\pm i}e^{-im\cdot x}
\ud x\right)
\\&=&\ds\frac{4(2\uppi)^d}{m^2\pm 2k\cdot m},\end{array}\]
the sign of which is determined by the sign of $m^2\pm 2k\cdot m$.
We count only the situation where these quantities are negative; reproducing a discussion made in the proof of Proposition~\ref{SpecL}, we conclude that 
\[
N^+_n\geq \#\{m\in \mathbb Z^d\smallsetminus\{0\},\ \sigma_{1,m}=0 \textrm { and } 
m^4-4(k\cdot m)^2<0\}.\]

When $m^4-4 (k\cdot m)^2=0$, the eigenvalues of $M_n^2$ are $0$ and $m^4$, 
and we just have to consider the positive eigenvalue $\alpha=m^4$, associated to the eigenvector $X_m=(  e^{im\cdot x},\pm ie^{im\cdot x},0,0)$ (depending whether $\frac{m^2}{2}=\mp k\cdot m$).
The equation $\mathscr LY_m=X_m$ \
has infinitely many solutions of the form $(2/m^2 e^{im\cdot x}, 0, 0 , 0) +z(\pm i e^{im\cdot x},e^{im\cdot x},0,0)$, with $z\in \mathbb C$. We deduce that
$(\mathscr L^{-1}X_m|X_m)=\frac{2(2\uppi)^d}{m^2}>0$. Thus these modes do not affect the counting.
\item When $m\neq 0$ and $\sigma_{1,m}\neq 0$, 
we are led to the relations $(-\alpha/c^2-\Delta_z)\pi_m=\frac\gamma c\sigma_2(2\uppi)^d\sigma_{1,m}p_m$,
 $(-\alpha/c^2-\Delta_z)\tilde \phi_m=-2(-\Delta_z)^{1/2}\sigma_2\gamma(2\uppi)^d\sigma_{1,m}\tilde q_m$.
The only solutions with square integrability on $\mathbb R^n$ are
$\pi_m=0$, $\tilde \phi_m=0$, $p_m=0$, $\tilde q_m=0$.
This can be seen by means of Fourier transform: $(-\alpha/c^2-\Delta_z)\phi=\sigma$
amounts to $\widehat \phi(\xi)=\frac{\widehat \sigma(\xi)}{|\xi|^2-\alpha/c^2}$;
due to \ref{H4} this function has a singularity which cannot be square-integrable.
In turn, this equally implies $\phi_m=0$ and $\tilde \pi_m=0$.
Hence, we arrive at $\frac{m^2}{2}q_m=0$ and 
$-ik\cdot m q_m=\alpha \tilde p_m$, together with $q_m= i k\cdot m\tilde p_m$ and $\frac{m^2}{2}\tilde p_m=0$.
We conclude that $\alpha>0$ cannot be an eigenvalue associated to a $m$-mode such that $m\neq 0$ and $\sigma_{1,m}\neq 0$.
\end{itemize}

\end{Proof}

We can now make use of Theorem~\ref{ThCP}, together with
Proposition~\ref{SpecL}. This leads to 
\begin{align*}
  0+1 
+\#\{m\in \mathbb Z^d\smallsetminus\{0\}, \sigma_{1,m}=0,\textrm{ and } m^4-4(k\cdot m)^2<0\}
+N_{C^+}= N^-_n+N^0_n+N^+_n+N_{C^+}\\
= n(\mathscr L)= 1+\#\{m\in \mathbb Z^d\smallsetminus\{0\}, m^4-4(k\cdot m)^2<0 \text{ and } 
\sigma_{1,m}=0\}\\
+
\#\{m\in \mathbb Z^d\smallsetminus\{0\}, m^4-4(k\cdot m)^2\leq 0 \text{ and } 
\sigma_{1,m}\neq  0\}
\end{align*}
so that
\[N_{C^+}
=\#\{m\in \mathbb Z^d\smallsetminus\{0\}, m^4-4(k\cdot m)^2\leq 0 \text{ and } 
\sigma_{1,m}\neq 0\}.\]
 Since 
the eigenvalue problem \eqref{Gene} does not have negative (real) eigenvalues, this is the only source of instabilities.

As a matter of fact, when $k=0$, we obtain $N_{C^+}=0$, 
which yields the following statement,  (hopefully!) consistent with Lemma \ref{mod0SW}
and Proposition~\ref{LinStabk0SW}.

\begin{coro}\label{spectralstabSWk0} Let $k=0$ and $\omega\in \mathbb R$ such that the dispersion relation \eqref{eq:dispersion} is satisfied. Suppose \eqref{small} holds. Then the plane wave solution $(e^{i\omega t}\mathbf 1(x), -\gamma\Gamma(z)\avex{\sigma}, 0)$ is spectrally stable.
\end{coro}

In contrast to what happens for the Hartree equation, for which the eigenvalues are purely imaginary, see Lemma~\ref{SpecStabH}, we can find unstable modes, despite the smallness condition \eqref{small}. Let us consider the following two examples in dimension $d=1$, with $k\in \mathbb Z\smallsetminus\{0\}$.
 
 \begin{example}
 Suppose $\sigma_{1,0}\neq 0$, and  $\sigma_{1,1}\neq 0$.
 Then, the set $\{m\in \mathbb Z\smallsetminus\{0\},\ m^4-4k^2m^2\leq 0\textrm{ and } 
 \sigma_{1,m}\neq 0\}$ contains $\{-1,+1\}$ (since $4k^2\geq 1$).
 Let $k\in \mathbb Z\smallsetminus\{0\}$ and $\omega\in \mathbb R$ such that the dispersion relation \eqref{eq:dispersion} is satisfied. Then the plane wave solution $(e^{i\omega t}e^{ikx}, -\gamma\Gamma(z)\avex{\sigma_1}, 0)$ is spectrally unstable.
\end{example}

 \begin{example}
 Let $m_*\in  \mathbb Z\smallsetminus\{0\}$ be the first Fourier mode
 such that $\sigma_{1,m_*}\neq 0$. Let $k\in \mathbb Z$ and $\omega\in\mathbb R$ such that the dispersion relation \eqref{eq:dispersion} is satisfied.
 Then, for all $k\in \mathbb Z$ such that $4k^2<m_*^2$,
 the plane wave solution $(e^{i\omega t}e^{ikx}, -\gamma\Gamma(z)\avex{\sigma}, 0)$ is spectrally stable, while for all $k\in \mathbb Z$ such that $4k^2\geq m_*^2$, the plane wave solution $(e^{i\omega t}e^{ikx}, -\gamma\Gamma(z)\avex{\sigma_1}, 0)$ is spectrally unstable.
\end{example}

 In general, if $k\in \Z^d\smallsetminus\{0\}$, the set  $\{m\in \mathbb Z^d\smallsetminus\{0\}, m^4-4(k\cdot m)^2\leq 0 \text{ and } 
\sigma_{1,m}\neq 0\}$ contains $-k$ and $k$ provided $\sigma_{1,k}\neq 0$. Hence, we have the following result. 
\begin{coro}\label{spectralinstabSWk} Let $k\in \Z^d\smallsetminus\{0\}$ and $\omega\in\mathbb R$ such that the dispersion relation \eqref{eq:dispersion} is satisfied. Suppose \eqref{small} holds and $\sigma_{1,m}\neq 0$ for all $m\in \mathbb Z^d\smallsetminus\{0\}$. Then the plane wave solution $(e^{i(\omega t+k\cdot x)},-\gamma \Gamma(z)\avex{\sigma_1},0)$ is spectrally unstable.
\end{coro}

%\begin{rmk}[Orbital instability] Given Corollary \ref{spectralinstabSWk}, it is natural to ask whether in this case the plane wave solution $(e^{i(\omega t+k\cdot x)},-\gamma \Gamma(z)\avex{\sigma_1},0)$ is orbitally unstable.
%
%Note that, if $\sigma_{1,m}\neq 0$ for all $m\in \mathbb Z^d\smallsetminus\{0\}$, we deduce from Proposition \ref{SpecL} that $n(\mathcal L)\ge 3$. As a consequence, the arguments used in \cite{GSS} to prove the orbital instability (see also \cite{Maeda-10,Ohta-11}) do not apply. It seems then necessary to work directly with the propagator generated by the linearized operator as in \cite{GSS2,GeoOht-10}. In particular, one has to establish Strichartz type estimates for the propagator of $\mathbb L$ (a task we leave for future work).
%\end{rmk}

\subsection{Orbital instability}

 Given Corollary \ref{spectralinstabSWk}, it is natural to ask whether or not the plane wave solution with $k\neq 0$ is orbitally unstable in this case.

\begin{theo}
\label{Th:unsta}
Let $k\in \Z^d\smallsetminus\{0\}$ and $\omega\in\mathbb R$ such that the dispersion relation \eqref{eq:dispersion} is satisfied. Suppose \eqref{small} holds and $\sigma_{1,m}\neq 0$ for all $m\in \mathbb Z^d\smallsetminus\{0\}$. Then the plane wave solution $(e^{i(\omega t+k\cdot x)},-\gamma \Gamma(z)\avex{\sigma_1},0)$ is orbitally unstable.
\end{theo}

Note that, if $\sigma_{1,m}\neq 0$ for all $m\in \mathbb Z^d\smallsetminus\{0\}$, we deduce from Proposition \ref{SpecL} that $n(\mathcal L)\ge 3$. As a consequence, the arguments used in \cite{GSS} to prove the orbital instability (see also \cite{Maeda-10,Ohta-11}) do not apply. It seems then necessary to work directly with the propagator generated by the linearized operator as in \cite{CCO,GeoOht-10,GSS2}. 
These arguments are of different nature: the former relies 
 on specific spectral properties of the self-adjoint operator $\mathscr L$, the latter  uses the existence of at least an eigenvalue of the linearized operator $\mathbb L$ with positive real part, a fact which has been just justified by the counting
argument. 

We go back to the non linear problem \eqref{sw_k}.
More precisely, we write $u(t,x)=e^{i\omega t}(\mathbf 1+\tilde u(t,x))$ and 
$\Psi(t,x,z)=-\gamma \avex{\sigma_1} \Gamma(z) +\tilde \Psi(t,x,z)$, where the perturbation $(\tilde u,\tilde \Psi)$ now satisfies
\begin{equation}\label{etoile}
\begin{array}{l}
i\partial_ t\tilde u +\ds\frac{\Delta_x \tilde u}{2}+ik\cdot \nabla_x \tilde u =
\gamma\sigma_1\star\ds\int_{\mathbb R^n}\sigma_2\tilde \Psi\ud z
+\left(\gamma  \sigma_1\star\ds\int_{\mathbb R^n}\sigma_2\tilde \Psi\ud z\right)\tilde u,
\\
\ds\frac{1}{c^2}\partial^2_{tt}\tilde \Psi-\Delta _z  \tilde \Psi
=-2\gamma \sigma_2\sigma_1\star\mathrm{Re}(\tilde u)
-\gamma\sigma_2\sigma_1\star |\tilde u|^2. 
\end{array}\end{equation}
Showing that the plane wave solution is orbitally instable   is then equivalent to prove that the solution 
$(0,0)$ of \eqref{etoile} is orbitally instable. 
By setting $\tilde \Psi=(-\Delta)^{-1/2}\phi$ and $\pi=-\frac{(-\Delta)^{-1/2}\partial_t \phi}{c}$ as before, 
we obtain that  \eqref{etoile}
can be expressed as a perturbation from the linearized equation 
\begin{equation}\label{perteq}\partial_t X=\mathbb LX+F(X).\end{equation}
The strategy consists in showing that we can exhibit initial data, as small as we wish, such that the 
solution exits a certain ball in finite time.
The exit time is related to the logarithm of the inverse of the size of the initial perturbation  (the smaller the initial
data, the larger the exit time).
In \eqref{perteq}, the non linear reminder is given by 
\[
F(X)=
\begin{pmatrix}
-\gamma p \sigma_1\star \ds\int_{\mathbb R^n} (-\Delta)^{-1/2}\sigma_2\phi\ud z
\\
\gamma q \sigma_1\star \ds\int_{\mathbb R^n} (-\Delta)^{-1/2}\sigma_2\phi\ud z
\\
0
\\
%\ds\frac
\gamma c\sigma_2 \sigma_1\star (|q|^2+|p|^2)
%\ds\int_{\mathbb T^d} \sigma_1(x-y) (q^2(y)+p^2(y))\ud y
\end{pmatrix},
\]
and $\mathbb L:D(\mathbb L)\subset \mathscr V\rightarrow \mathscr V$  is the linear operator defined in \eqref{defopLd}.

\begin{lemma}\label{EstNL}
We can find a constant $C_F$ such that, for any $X$, there holds 
$\|F(X) \|_{\mathscr V}\leq C_F\|X\|_{\mathscr V}^2$.
\end{lemma}

\begin{Proof}
For the first two components of $F(X)$, it suffices to obtain a uniform estimate on the potential
\[\begin{array}{l}
\Big|
\sigma_1\star \ds\int_{\mathbb R^n} (-\Delta)^{-1/2}\sigma_2\phi\ud z
\Big|
=
\left|
\ds\int_{\mathbb T^d}
\sqrt{\sigma_1(x-y)}\sqrt{\sigma_1(x-y)} \ds\int_{\mathbb R^n} (-\Delta)^{-1/2}\sigma_2(z)\phi(y,z)\ud z\ud y
\right|
\\
[.3cm]
\qquad\leq 
\left(\ds\int_{\mathbb T^d}\sigma_1(y)\ud y\right)^{1/2}
\left(\ds\int_{\mathbb T^d}\sigma_1(x-y)
\Big| \ds\int_{\mathbb R^n} (-\Delta)^{-1/2}\sigma_2(z)\phi(y,z)\ud z\Big|^2
\ud y
\right)^{1/2}
\\
[.3cm]
\qquad\leq \sqrt{\avex{\sigma_1}}
\left(\ds\int_{\mathbb T^d}\sigma_1(x-y)
 \ds\int_{\mathbb R^n} \ds\frac{\widehat \sigma_2(\xi)}{|\xi|^2}\ud \xi
  \ds\int_{\mathbb R^n} |\phi(y,z)|^2\ud z
\ud y
\right)^{1/2}\\
[.3cm]
\qquad\leq \sqrt{\avex{\sigma_1}}
\sqrt{\kappa\|\sigma_1\|_{L^\infty(\mathbb T^d)}}
\left(\ds\iint_{\mathbb T^d\times \mathbb R^n} |\phi(y,z)|^2\ud z\ud y\right)^{1/2}
.
\end{array}\]
It implies that the $L^2$ norm of the  first 
component of $F(X)$ is dominated by 
\[
\gamma\sqrt{\avex{\sigma_1}}
\sqrt{\kappa\|\sigma_1\|_{L^\infty(\mathbb T^d)}} \|p\|_{L^2(\mathbb T^d)}\|\phi\|_{L^2(\mathbb T^d\times \mathbb R^n)},\]
and a similar estimate holds for the second component.
Finally, for the forth component of $F(X)$, we get, with $|u|^2=|q|^2+
|p|^2$,
\[\begin{array}{lll}
\ds\iint_{\mathbb T^d\times \mathbb R^n} 
|\sigma_2(z)|^2\ | \sigma_1\star |u|^2(x)|^2
\ud z\ud x
&\leq & 
\|\sigma_2\|_{L^2(\mathbb R^n)}^2
\ds\int_{\mathbb T^d}
\left|
\ds\int_{\mathbb T^d}
\sigma_1(x-y)|u|(y)\times |u|(y)\ud y
\right|^2\ud x
\\
[.3cm]
&\leq & 
\|\sigma_2\|_{L^2(\mathbb R^n)}^2
\ds\int_{\mathbb T^d}
\ds\int_{\mathbb T^d}|\sigma_1|^2(x-y)|u|^2(y)\ud y
\ds\int_{\mathbb T^d}|u|^2(y)\ud y\ud x\\
[.3cm]
&\leq & 
\|\sigma_2\|_{L^2(\mathbb R^n)}^2
\|\sigma_1\|_{L^2(\mathbb T^d)}^2
\left(\ds\int_{\mathbb T^d}|u|^2(y)\ud y\right)^2
\end{array}
\]
Hence the $L^2$ norm of the last component of $F(X)$ is dominated by 
\[
%\ds\frac
\gamma c \|\sigma_2\|_{L^2(\mathbb R^n)}\|\sigma_1\|_{L^2(\mathbb T^d)}(
\|q\|_{L^2(\mathbb T^d)}^2 + \|p\|_{L^2(\mathbb T^d)}^2)
\]
\end{Proof}

Next, we are going to use the Duhamel formula
\begin{equation}\label{Duhamel}
X(t)=e^{\mathbb L t}X(0)+\ds\int_0^t e^{\mathbb L (t-s)}F(X(s))\ud s.\end{equation}
The definition of the operator semi-group $\{e^{\mathbb L t}, \ t\geq 0\}$
follows from the application of Lumer-Phillips' theorem \cite[Th.~12.22]{Ren} by combining  the basic estimate 
\[\begin{array}{lll}
|\langle \mathbb LX|X\rangle|
&=&\left|
-\gamma \ds\int_{\mathbb T^d} p\sigma_1\star \Big(
\ds\int_{\mathbb R^n}(-\Delta)^{-1/2}\sigma_2\phi\ud z\Big)\ud x
+2c\gamma \ds\iint_{\mathbb T^d\times \mathbb R^n} \sigma_2\pi \sigma_1\star q\ud z\ud x
\right|
\\
&\leq& \gamma  \avex{\sigma_1}\Big(\sqrt \kappa 
+2c\sqrt { 
\|\sigma_2\|_{L^\infty(\mathbb R^n)}\|\sigma_2\|_{L^1(\mathbb R^n)}}
\Big) \|X\|_{\mathscr V}^2,
\end{array}\]
together with the following claim.
\begin{lemma}\label{LIsOnto}
%which allows us to apply the Hille-Yosida theorem \cite[Chap. VII]{Brez}.
There exists $\lambda_*>0$ such that for any real $\lambda\geq  \lambda_*$, the operator $\lambda -\mathbb L$ is onto.
\end{lemma}

\begin{Proof}
We try to solve the system
\[\begin{array}{l}
\lambda q +\ds\frac{\Delta_x p}{2}+k\cdot\nabla_x q=q',
\\
\lambda p -\ds\frac{\Delta_x q}{2}+k\cdot\nabla_x p+\gamma\sigma_1\star\ds\int_{\mathbb R^n}(\Delta_z)^{-1/2}\sigma_2\phi\ud z=p',
\\
\lambda \phi+c(-\Delta)^{1/2}\pi=\phi',
\\
\lambda \pi -c(-\Delta)^{1/2}\phi-2c \gamma\sigma_2\sigma_1\star q=\pi',
\end{array}\]
with $\lambda\in \mathbb R\setminus\{0\}$.
By using the Fourier transform, the last two equations become
\[\widehat \pi=\ds\frac{-\lambda \widehat\phi+\widehat\phi'}{c|\xi|},
\qquad
\lambda \widehat \pi -c|\xi| \widehat\phi-2c \gamma\widehat{\sigma_2}\sigma_1\star q=\widehat\pi',
\]
which yields
\[
\widehat\phi(x,\xi)=\ds\frac{
\lambda\widehat\phi'(x,\xi)/c^2-|\xi|\widehat \pi'(\xi)/c -2\gamma |\xi|\widehat{\sigma_2}(\xi)\sigma_1\star q(x)
}{\lambda^2/c^2+|\xi|^2}.
\]
%Note that this formula is meaningful provided  $\lambda\notin i\mathbb R$. 
Let us introduce the quantity
\[\mu\in \mathbb R \longmapsto \kappa_\mu=\ds\int_{\mathbb R^n} \ds\frac{|\widehat \sigma_2(\xi)|^2}{\mu^2+|\xi|^2}\ds\frac{\ud \xi}{(2\uppi)^n}.\]
The function $\mu\mapsto \kappa_\mu$ is non increasing on $[0,\infty)$, and the inequality $0\leq \kappa_\mu\leq \kappa$ holds for any $\mu \in \mathbb R$.
Reasoning by means of Fourier coefficients we are led to
\[\begin{array}{l}
\begin{pmatrix}
\lambda +ik\cdot m & -m^2/2
\\
m^2/2-2\gamma^2(2\uppi)^{2d}|\sigma_{1,m}|^2\kappa_{\lambda^2/c^2}& \lambda +ik\cdot m
\end{pmatrix}
\begin{pmatrix}
q_m\\
p_m
\end{pmatrix}=\begin{pmatrix}q_m'
\\
S_m
\end{pmatrix}
\end{array}\]
with 
\[
S_m=
p'_m-\gamma(2\uppi)^d\sigma_{1,m} \ds\int_{\mathbb R^n} 
\ds\frac{\widehat \sigma_2(\xi)}{|\xi|}\ds\frac{\lambda \widehat \phi'_m(\xi)/c^2-|\xi|\widehat \pi'_m(\xi)/c}{\lambda^2/c^2+|\xi|^2}\ds\frac{\ud \xi}{(2\uppi)^n}
\]
Since $\lambda^2/c^2+|\xi|^2\geq \lambda^2/c^2$, we observe that the $\ell^2$ norm of the right hand side $S_m$ is dominated by
\[
 \|p'\|_{L^2(\mathbb T^d)}+\gamma \avex{\sigma_{1}}
\left(
\ds\frac{\sqrt\kappa}{|\lambda| }\|\phi'\|_{L^2(\mathbb T^d\times\mathbb R^n)}
+\ds\frac{c}{|\lambda|^2}\|\sigma_2\|_{L^2(\mathbb R^n)} \|\pi'\|_{L^2(\mathbb T^d\times\mathbb R^n)}\right)
.\]
We obtain $\lambda q_0=q'_0$, $\lambda p_0=S_0+2\gamma^2 \avex{\sigma_1}^2\kappa_{\lambda^2/c^2}q_0$
and, for $m\neq 0$, \[\begin{array}{l}
\underbrace{\Big((\lambda + ik\cdot m)^2+\ds\frac{|m|^4}{4}
\Big(1-4\gamma^2(2\uppi)^{2d}\ds\frac{|\sigma_{1,m}|^2}{m^2}\kappa_{\lambda^2/c^2}\Big)
\Big)}_{=R_m(\lambda)}q_m
=(\lambda + ik\cdot m)q'_m+\ds\frac{m^2}{2}S_m,
\\
p_m=\ds\frac{2}{m^2}(\big({\lambda+ik\cdot m} )q_m-q'_m\big).
\end{array}\]
By virtue  of \eqref{small}, $1-4\gamma^2(2\uppi)^{2d}\frac{|\sigma_{1,m}|^2}{m^2}\kappa_{\lambda^2/c^2}
\geq 1-4\gamma^2\avex {\sigma_{1}}^2\kappa >
0$, 
so that the coefficient $R_m(\lambda)$ does not vanish:
either its imaginary part $\lambda k\cdot m\neq 0$, or 
when  $k\cdot m=0$, its real part $\lambda^2 +\frac{m^4}{4}(1-4\gamma^2(2\uppi)^{2d}\frac{|\sigma_{1,m}|^2}{m^2}\kappa_{\lambda^2/c^2})$ is bounded from below by a positive quantity.
It remains to check that 
\[
q_m=\ds\frac{(\lambda + ik\cdot m)q'_m+\frac{m^2}{2}S_m}{R_m(\lambda)}
\]
defines a square-summable sequence, at least when $\lambda$ is large enough.
To this end, for $m\neq 0$, we evaluate
\[\begin{array}{lll}
|R_m(\lambda)|^2
&=&
\Big|2i\lambda k\cdot m+
\lambda^2 - (k\cdot m)^2+\ds\frac{|m|^4}{4}
\Big(1-4\gamma^2(2\uppi)^{2d}\ds\frac{|\sigma_{1,m}|^2}{m^2}\kappa_{\lambda^2/c^2}\Big) 
\Big|^2
\\[.4cm]
&=&
4\lambda^2 (k\cdot m)^2
+\big(\lambda^2 - (k\cdot m)^2\big)^2 
+\Big(\ds\frac{|m|^4}{4}
\Big(1-4\gamma^2(2\uppi)^{2d}\ds\frac{|\sigma_{1,m}|^2}{m^2}\kappa_{\lambda^2/c^2}\Big)\Big)^2
\\[.4cm]
&&
\hspace*{2cm}+\big(\lambda^2 - (k\cdot m)^2\big)\ds\frac{|m|^4}{2}
\Big(1-4\gamma^2(2\uppi)^{2d}\ds\frac{|\sigma_{1,m}|^2}{m^2}\kappa_{\lambda^2/c^2}\Big)
\\[.4cm]
&=&
\big(\lambda^2 + (k\cdot m)^2\big)^2 
+\Big(\ds\frac{|m|^4}{4}
\Big(1-4\gamma^2(2\uppi)^{2d}\ds\frac{|\sigma_{1,m}|^2}{m^2}\kappa_{\lambda^2/c^2}\Big)\Big)^2
\\[.4cm]&&
\hspace*{2cm}+
\big(\lambda^2 - (k\cdot m)^2\big)\ds\frac{|m|^4}{2}
\Big(1-4\gamma^2(2\uppi)^{2d}\ds\frac{|\sigma_{1,m}|^2}{m^2}\kappa_{\lambda^2/c^2}\Big)
\\[.4cm]
&\geq 
&\Big(\ds\frac{|m|^4}{4}
\Big(1-4\gamma^2(2\uppi)^{2d}\ds\frac{|\sigma_{1,m}|^2}{m^2}\kappa_{\lambda^2/c^2}\Big)\Big)^2
\\[.4cm]&&
\hspace*{2cm}+
\big(\lambda^2 - k^2 m^2\big)\ds\frac{|m|^4}{2}
\Big(1-4\gamma^2(2\uppi)^{2d}\ds\frac{|\sigma_{1,m}|^2}{m^2}\kappa_{\lambda^2/c^2}\Big)
.
\end{array}\]
Let $\delta>0$, that will be made precise later on.
We split the last term depending whether $k^2\geq \delta m^2$ or $k^2<\delta m^2$:
\[\begin{array}{l}
\big(\lambda^2 - (k\cdot m)^2\big)\ds\frac{|m|^4}{2}
\Big(1-4\gamma^2(2\uppi)^{2d}\ds\frac{|\sigma_{1,m}|^2}{m^2}\kappa_{\lambda^2/c^2}\Big)
\mathbf 1_{k^2\geq \delta m^2}
\\[.4cm]
\qquad\geq 
\big(\lambda^2 - k^4/\delta\big)\ds\frac{|m|^4}{2}
\Big(1-4\gamma^2(2\uppi)^{2d}\ds\frac{|\sigma_{1,m}|^2}{m^2}\kappa_{\lambda^2/c^2}\Big)
\mathbf 1_{k^2\geq \delta m^2}
\end{array}\]
and
\[\begin{array}{l}
\big(\lambda^2 - (k\cdot m)^2\big)\ds\frac{|m|^4}{2}
\Big(1-4\gamma^2(2\uppi)^{2d}\ds\frac{|\sigma_{1,m}|^2}{m^2}\kappa_{\lambda^2/c^2}\Big)
\mathbf 1_{k^2< \delta m^2}
\\[.4cm]
\qquad\geq 
\big(\lambda^2 - \delta m^4\big)\ds\frac{|m|^4}{2}
\Big(1-4\gamma^2(2\uppi)^{2d}\ds\frac{|\sigma_{1,m}|^2}{m^2}\kappa_{\lambda^2/c^2}\Big)
\mathbf 1_{k^2< \delta m^2}
\\[.4cm]
\qquad\geq 
 - \delta\ds\frac{m^8}{2}
\Big(1-4\gamma^2(2\uppi)^{2d}\ds\frac{|\sigma_{1,m}|^2}{m^2}\kappa_{\lambda^2/c^2}\Big)
\mathbf 1_{k^2< \delta m^2}.
\end{array}\]
When $\lambda\geq \lambda_*=k^2/\sqrt\delta$, we can get rid of the first term in the evaluation of $|R_m(\lambda)|^2$ and we arrive at 
\[\begin{array}{lll}
|R_m(\lambda)|^2&\geq & 
\ds\frac{m^8}{16}\Big(1-4\gamma^2(2\uppi)^{2d}\ds\frac{|\sigma_{1,m}|^2}{m^2}\kappa_{\lambda^2/c^2}\Big)
\Big\{
\mathbf 1_{k^2\geq \delta m^2}\big(1-4\gamma^2\avex{\sigma_1}^2\kappa\big)
\\
[.4cm]
&&\hspace*{5cm}
+
\mathbf 1_{k^2< \delta m^2}\Big(
\big(1-4\gamma^2\avex{\sigma_1}^2\kappa\big)
-8\delta
\Big)
\Big\}.\end{array}\]
We choose $\delta$ so that the last term contributes positively, for instance
$\delta =\frac{1-4\gamma^2\avex{\sigma_1}^2\kappa}{16}$.
Having defined this way $\delta$, and thus $\lambda_*$, we exhibit $c_*>0$ such that
$|R_m(\lambda)|^2\geq c_* m^8$.
Combined to the $\ell^2$  estimate on $S_m$, this allows us to conclude that $\|X\|_{\mathscr V}
=\|(\lambda-\mathbb L)^{-1}X'\|_{\mathscr V}\leq M\|X'\|_{\mathscr V}$ holds for a certain constant $M$, when $\lambda\geq \lambda_*$.
\end{Proof}

Moreover, 
a continuity estimate holds: we can find $\Lambda>0$ 
such that  for any $t\geq 0$, $\|e^{\mathbb Lt}\|_{\mathscr L(\mathscr V)} \leq e^{\Lambda t}$.
Let us also introduce
\[K_0=\sup\big\{\|e^{\mathbb Lt}\|_{\mathscr L(\mathscr V)}, 0\leq t\leq 1\big\}.\]

The proof of instability slightly simplifies when $\sigma(e^{\mathbb L})=e^{\sigma(\mathbb L)}$, see  \cite{SRN3}, and the references therein, for a  situation where this equality is fulfilled.
 According to 
Gearhart-Greiner-Herbst-Pr\"uss' theorem, see 
\cite[Prop.~1]{Pruss} and the formulation proposed in \cite[Section~2]{Gesz}), 
such identification holds provided  the resolvent $(\lambda-\mathbb L)^{-1}$ satisfies  a uniform estimate as $\mathrm{Im}(\lambda)\to \pm \infty $ with $\mathrm{Re}(\lambda) \neq 0$ fixed, which is far from obvious.
Nevertheless, the arguments of \cite{ShSt}
only relies on the trivial embedding
$e^{\sigma(\mathbb L)}\subset \sigma(e^{\mathbb L})$.

We are concerned with the case where  spectral instability  holds, which means that  $\mathbb L$ has eigenvalues with positive real value.
 There is only a finite number of such eigenvalues (as indicated by the counting argument).
 In turn, the spectral radius of 
 $e^{\mathbb L}$ is larger than 1.  Let  $\lambda_*=a_*+ib_*$ with $a_*>0$, be  such that $e^{\lambda_*}$ lies in the boundary of 
 $\sigma(e^{\mathbb L})$:
 \[
 |e^{\lambda_*}|=e^{a_*}=\max\big\{|\mu|,\ \mu \in \sigma(e^{\mathbb L})\big\}.\]
 
%  \[
% \lambda_*=a_*+ib_*\in \sigma(\mathbb L),\quad
% a_*=\sup\big\{\mathrm{Re}(\lambda),\ \lambda\in \sigma(\mathbb L)\big\}>0.\]
%  Of course, for any $t\geq 0$, we have 
% $|e^{\lambda_* t}|=e^{a_*t}$ and the spectral radius of $e^{\mathbb L}$ is $e^{a_*}>1$ \cite{GeoOht-10}.
%Finally, we need the following claim.

\begin{lemma}\cite[Lemma~2 \& Lemma~3]{ShSt}\label{LShSt}
The following assertions hold:
\begin{enumerate}
\item For any $\gamma>0$ and any $m\in \mathbb N\setminus\{0\}$, there exists $Y_*\in \mathscr V$ such that $\|Y_*\|_{\mathscr V}=1$ and $\|(e^{m\mathbb L}-e^{m\lambda_*})Y_*\|_{\mathscr V}
\leq \gamma  $;
\item For any $0\leq t\leq m$, we have
$\|e^{t\mathbb L}Y_*\|_{\mathscr V}\leq 2K_0e^{a_*t}$;
\item 
There exists a constant $K_1$, such that for any $t\geq 0$, there holds
$e^{a_*t}\leq \|e^{t\mathbb L}\|_{\mathscr L(\mathscr V)}\leq K_1e^{3a_*t/2}$.
\end{enumerate}\end{lemma}

Let us define $\epsilon$ such that 
\[\ds\frac{4  K_1(2K_0+C_F)^2e^{a_*} }{a_*}\epsilon<1,
\qquad \ds\frac{8 K_1C_F(2K_0+C_F)^2 e^{2a_*} }{a_*}\epsilon<1.
\]
Then, pick  $\delta >0$ as small as we wish and set
$$T_\delta=\ds\frac{1}{a_*}\ln\Big(\ds\frac\epsilon \delta\Big),
\qquad m_\delta=\lfloor T_\delta \rfloor +1.
$$
Let $Y_*$ be a normalized 
vector as defined by Lemma~\ref{LShSt}-1 with $\gamma=\frac\epsilon {2\delta}$ and $m=m_\delta$.
%eigenvector of $\mathbb L$ associated to $\lambda_*$:
%\[\mathbb LY_*=\lambda_*Y_*,\qquad \|Y_*\|_{\mathscr V}=1.\]
The initial data
\[
X\big|_{t=0}=\delta Y_*,\]
 has thus an arbitrarily small norm.
Now, \eqref{Duhamel} becomes
\[
X(t)
=\delta %e^{\lambda_*t} 
e^{\mathbb L t} 
Y_{*} + \ds\int_0^t e^{\mathbb L(t-s)}F(X(s))\ud s.\]
We 
are going to
contradict the orbital stability by showing that  $\|X(m_\delta)\|_{\mathscr V}>\epsilon/4$:
the solution always exits the ball $B(0,\epsilon/4)$.

Let $$\tilde T_\delta
=\ds\sup\big\{t\in [0,m_\delta], \|X(s)-\delta e^{\mathbb L s} 
Y_{*}\|_{\mathscr V}
\leq \delta C_F e^{a_*s}, \textrm{ for }0\leq s\leq t\big\}
\in (0,m_\delta].$$ 
As a consequence of \eqref{Duhamel}, together with Lemma~\ref{EstNL} and~\ref{LShSt}-3, 
we get
\[
\|X(t)- \delta e^{\mathbb L t} 
Y_{*} \|_{\mathscr V}\leq \ds\int_0^t  K_1 e^{3a_*(t-s)/2}C_F\|X(s)\|_{\mathscr V}^2 \ud s.
\]
It follows that, for $0\leq t\leq \tilde T_\delta< m_\delta $,
\[
\begin{array}{lll}
\|X(t)-\delta e^{\mathbb L t} 
Y_{*}\|_{\mathscr V}&\leq&K_1C_F %(1+C_F)^2 \delta^2
 \ds\int_0^t e^{3a_*(t-s)/2} 
 \big|\delta \|  e^{\mathbb L s} 
Y_{*} \|_{\mathscr V}+ \|X(s)-\delta e^{\mathbb L s} 
Y_{*}\|_{\mathscr V}\big|^2
 \ud s
\\[.3cm]
&\leq&K_1C_F 
\ds\int_0^t e^{3a_*(t-s)/2} 
 \big|2\delta K_0 e^{a_* s}+\delta C_F e^{a_*s}\big|^2
 \ud s
 \\[.3cm]
 &&\hspace*{4cm}\textrm{ (by using Lemma~\ref{LShSt}-2)}
\\[.3cm]
&\leq&\delta^2 K_1C_F (2K_0 +C_F)^2 e^{3a_* t/2}
\ds\int_0^t 
 e^{a_* s/2}
 \ud s
\\[.3cm]
&\leq&
\ds\frac{2}{a_*} K_1C_F (2K_0 +C_F)^2 
\big(\delta e^{a_* t}\big)^2
\leq \epsilon \ds\frac{2 e^{a_*} }{a_*} K_1C_F (2K_0 +C_F)^2 \delta e^{a_* t}
%
%
%
%&\leq& \delta e^{a_*t}\left(1+\ds\frac{2K_1C_F(1+C_F)^2}{a_*}\delta e^{a_*T_\delta}\right)
%\leq  \delta e^{a_*t}\left(1+C_F\ds\frac{2K_1(1+C_F)^2}{a_*}\epsilon\right)
\end{array}\]
holds.
Hence, $\epsilon$ is chosen small enough so that 
this implies 
\[
\|X(t) - \delta e^{\mathbb L t} Y_*\|_{\mathscr V}
< \ds\frac{C_F}{2}\delta e^{a_*t},\]
which would contradict the definition of $\tilde T_\delta$ if $\tilde T_\delta< m_\delta$.
Accordingly,
\[\|X(t)-  \delta e^{\mathbb L t} Y_*\|_{\mathscr V}\leq   C_F\delta e^{a_*t}\]
holds for any $t\in [0,m_\delta]$.
Going back to the Duhamel formula thus 
yields, for $0\leq t\leq m_\delta$,
\[%\begin{array}{lll}
\|X(t)-\delta e^{\mathbb L t}Y_*\|_{\mathscr V}
%&\leq& 
%\ds\int_0^t \|e^{\mathbb L(t-s)}F(X(s))\|_{\mathscr V}\ud s
%\leq 
%\ds\int_0^t e^{3a_*(t-s)/2} K_1C_F \|X(s)\|_{\mathscr V}^2\ud s
%\\
%&\leq&
%K_1C_F (1+C_F)^2\delta^2 \ds\int_0^t e^{3a_*(t-s)/2} e^{2a_*s}\ud s
%=\ds\frac{2K_1C_F (1+C_F)^2}{a_*}\delta^2 e^{2a_*t}
%\\
%&\leq&
\leq  \ds\frac{2K_1C_F (2K_0+C_F)^2}{a_*}\delta^2 e^{2a_*m_\delta}.
%\end{array}
\]
Now, by using Lemma~\ref{LShSt}-1, we observe that
\[
\|e^{\mathbb L m_\delta} Y_*\|_{\mathscr V}
\geq \|e^{\lambda_* m_\delta} Y_*\|_{\mathscr V}-\ds\frac\epsilon {2\delta}
\geq e^{a_* m_\delta}-\ds\frac\epsilon {2\delta}\geq \ds\frac\epsilon {2\delta}
.\]
We deduce that
\[\begin{array}{lll}
\|X(m_\delta)\|_{\mathscr V}&\geq& 
\|\delta e^{\mathbb L m_\delta} Y_*\|_{\mathscr V}-
\|X(m_\delta)- \delta e^{\mathbb L m_\delta} Y_*\|_{\mathscr V}
\\[.3cm]
&\geq& \ds\frac\epsilon 2
-\ds\frac{2K_1C_F (2K_0+C_F)^2}{a_*}\delta^2 e^{2a_*m_\delta}
\\[.3cm]
&\geq& \epsilon \Big(
\ds\frac1 2 - \ds\frac{2K_1C_F (2K_0+C_F)^2e^{2a_*}}{a_*}\epsilon\Big)
> \ds\frac\epsilon 4
\end{array}
%\left(1
%-\ds\frac{2K_1C_F (1+C_F)^2}{a_*}\delta e^{a_*T_\delta}
%\right)
%=\epsilon\left(1-\ds\frac{2K_1C_F (1+C_F)^2}{a_*}\epsilon
%\right)\geq \ds\frac{\epsilon}2,
\]
as announced.
% {\color{black} (while $\|X(m_\delta)\|_{\mathscr V}\leq \|X(m_\delta)- \delta e^{\mathbb L m_\delta} Y_*\|_{\mathscr V}
%+\delta \| e^{\mathbb L m_\delta} Y_*\|_{\mathscr V}
%\leq \epsilon$).

 {\color{black}
That these estimates now imply the orbital instability of the plane wave solution, which  amounts to justify that  
\[
%\Theta_\epsilon=
\ds\inf_\theta\left\|X(m_\delta)+ \begin{pmatrix} 1 \\ 0 \\ -\gamma\avex{\sigma_1}\Gamma \\
0 \end{pmatrix} - \begin{pmatrix} \cos(\theta) \\ \sin(\theta) \\ -\gamma\avex{\sigma_1}\Gamma \\
0 \end{pmatrix} 
\right\|_{\mathscr V}\geq \kappa_*\epsilon\]
holds for a certain positive constant $\kappa_*$, 
follows by adapting the arguments of \cite[sp.~Theorem~6.2]{GSS2}, see also \cite{SRN3}.
%It amounts to justify that  
%\[
%\Theta_\epsilon=\ds\inf_\theta\left\|X(m_\delta)+ \begin{pmatrix} 1 \\ 0 \\ -\gamma\avex{\sigma_1}\Gamma \\
%0 \end{pmatrix} - \begin{pmatrix} \cos(\theta) \\ \sin(\theta) \\ -\gamma\avex{\sigma_1}\Gamma \\
%0 \end{pmatrix} 
%\right\|_{\mathscr V}\geq \kappa_*\epsilon\]
%holds for a certain positive constant $\kappa_*$.
%, where
%$R(\theta)$ stands for the (partial) rotation matrix 
%\[
%R(\theta)=\begin{pmatrix}
%\cos(\theta) & -\sin(\theta) & 0 & 0
%\\
%\sin(\theta) & \cos(\theta) & 0& 0 \\
%0 & 0 & 0 & 0\\
%0 & 0 & 0 & 0
%\end{pmatrix}\]
%and $X_*=(\mathbf 1,0,-\gamma\avex{\sigma_1}\Gamma,0)$.
%{\color{blue} The infimum is reached for an angle $\theta_\epsilon$ which tends to 0 as $\epsilon$ runs to 0, and
%we} get
%\[
%(R(\theta)-\mathbf I)X_*\underset{\theta\to 0}{\sim}(0,\theta,0,0)\]
%which lies in $\mathrm{Ker}(\mathbb L)$.
%In particular,  $Y_*$ has a non trivial part, $Y^\perp_*\neq 0$, orthogonal to  $\mathrm{Span}(0,\mathbf 1,0,0)$.
%It allows us to show that $\Theta_\epsilon \|Y^\perp_*\|_{\mathscr V}\geq \epsilon\frac{\|Y^\perp_*\|_{\mathscr V}^2}{2}$. }
}
\appendix

\section{Scaling of the model and physical interpretation}\label{adim}

It is worthwhile to discuss the meaning of the parameters 
that govern the equations and the asymptotic issues.
 Going back to physical units, the system reads
\begin{subequations}
    \begin{alignat}{11}
        \label{schro_niv0}&\ds \left(i\hbar \partial_{t}U+\frac{\hbar ^2}{2m}\Delta_{x}U\right)(t,x)=\left(\ds\int_{\mathbb T^d\times\mathbb R^n} 
        \sigma_{1}(x-y) \sigma_{2}(z)\Psi(t,y,z)\ud y\ud z\right)u(t,x),\\
        \label{wave_niv0}&(\partial_{tt}^{2}\Psi-\varkappa^{2}\Delta_{z}\Psi)(t,x,z)=-\sigma_{2}(z)\left(\ds\int_{\mathbb T^d} \sigma_{1}(x-y)|U(t,y)|^{2}
        \ud y\right).
    \end{alignat}
\end{subequations}
The quantum particle
is described by the wave function  $(t,x)\mapsto U(t,x)$: given $\Omega\subset \mathbb T^d$, the integral $\int_\Omega |U(t,x)|^2\ud x$ gives 
the probability of presence of the quantum  particle at time $t$ in the domain $\Omega$; this is a dimensionless quantity.
 In \eqref{schro_niv0}, $\hbar $ stands for the Planck constant; its  homogeneity is $\frac{\textrm{Mass}\times \textrm{Length}^2}{ \textrm{Time}}$ (and its value is $1.055 \times10^{-34}\  J s$) and $m$ is the inertial mass of the  particle.
Let us introduce mass, length and time units of observations: $\textrm M$, $\textrm L$ and $\textrm T$.
It helps the intuition to think of the $z$ directions as homogeneous to a length, but in fact this is not necessarily the case: we denote by $\Uppsi$ and $\mathrm Z$ the (unspecified)  units for $\Psi$ and the $z_j$'s. Hence, $\varkappa$ is homogeneous to the ratio $\frac{\mathrm Z}{\mathrm T}$.
The coupling between the vibrational field and the particle is driven by 
the product of the form functions $\sigma_1\sigma_2$, which   has the same homogeneity 
as 
$\frac{\hbar }{ \mathrm T  \Psi\mathrm L^d\mathrm Z^n}$
from \eqref{schro_niv0}
and 
as $\frac{ \Uppsi}{\mathrm L^d \mathrm T^2}$
from \eqref{wave_niv0}, both are thus measured with the same units. From now on, we denote by $\varsigma$ this coupling unit.
Therefore, we are led to the following dimensionless quantities 
\[\begin{array}{l}
U'(t',x')=U(t'\mathrm T, x'\mathrm L) \ \sqrt{\mathrm L^{d}\ds\frac{m}{\mathrm M}},
\\[.3cm]
\Psi'(t',x',z')=\ds\frac{1}{\Uppsi} \Psi(t'\mathrm T, x'\mathrm L,z'\mathrm Z), 
\\[.3cm]
\sigma_1'(x')\sigma_2(z')=\ds\frac1\varsigma\ \sigma_1(x'\mathrm L)\sigma_2(z'\mathrm Z).
\end{array}\]
Bearing in mind that $u$ is a probability density, we note that  
\[\ds\int_{\mathbb T^d} |U'(t',x')|^2\ud x'=\ds\frac{m}{\mathrm M}.\]
Dropping the primes, \eqref{schro_niv0}-\eqref{wave_niv0} becomes, in dimensionless form,
\begin{subequations}
    \begin{alignat}{11}
        \label{schro_res}&\ds \left(i\partial_{t}U+\frac{\hbar \mathrm T}{m\mathrm L^2}\frac12\Delta_{x}U\right)(t,x)=
        \ds\frac{ \varsigma\Uppsi \mathrm L^d\mathrm Z^n  \mathrm T}{\hbar}
        \left(\ds\int_{\mathbb T^d\times\mathbb R^n} 
        \sigma_{1}(x-y) \sigma_{2}(z)\Psi(t,y,z)\ud y\ud z\right)U(t,x),\\
        \label{wave_res}&\Big(\partial_{tt}^{2}\Psi-\frac{\varkappa^{2}\mathrm T^2}{\mathrm Z^2}
        \Delta_{z}\Psi\Big)(t,x,z)=-
       \ds\frac{ \varsigma\mathrm T^2}{\Uppsi}\ds\frac {\mathrm M}{m}
        \sigma_{2}(z)\left(\ds\int_{\mathbb T^d} \sigma_{1}(x-y)|U(t,y)|^{2}
        \ud y\right).
    \end{alignat}
\end{subequations}
Energy conservation plays a central role in the analysis of the system:
the total energy  is defined by using  the reference units  and
we obtain \[\begin{array}{l}
\mathscr E_0=
\Big(\ds\frac{\hbar \mathrm T}{m\mathrm L^2}\Big)^2\ds\frac12\ds\int_{\mathbb T^d}|\nabla_x U |^2 \ud x
+
\ds\frac{\Uppsi^2\mathrm L^d\mathrm Z^n}{\mathrm M\mathrm L^2}\ds\frac12 \ds\iint_{\mathbb T^d\times \mathbb R^n}
\Big(|\partial_t \Psi|^2 +\ds\frac{\varkappa^2\mathrm T^2}{\mathrm Z^2}  |\nabla_z\Psi|^2
\Big)
 \ud z\ud x
 \\
 \hspace*{3cm}
 +
 \varsigma \ds\frac{\Psi \mathrm L^d\mathrm Z^n\mathrm T^2} {m\mathrm L^2}
\ds \iint_{\mathbb T^d\times \mathbb R^n}
 |U|^2 \sigma_2\sigma_1\star \Psi \ud z\ud x
,
\end{array}\]
with $\mathscr E_0$ dimensionless (hence the total energy of the original system is $\mathscr E_0\frac{\mathrm{ML^2}}{\mathrm T^2}$). 
Therefore, we see that the dynamics is encoded by four independent parameters.
In what follows, we get rid of a parameter by assuming 
\[\ds\frac{\hbar \mathrm T}{m\mathrm L^2}=1,\]
and we work with the following three independent parameters
\[
\alpha=\ds\frac{\varsigma\Uppsi \mathrm L^d\mathrm Z^n   \mathrm T^2}{m\mathrm L^2}\ds\frac{m\mathrm L^2} {\hbar \mathrm T},\quad
\beta=\ds\frac{\varsigma \mathrm Z^2}{\varkappa^2\Uppsi}\ds\frac{\mathrm M}{m},
\quad
c=\ds\frac{\varkappa \mathrm T}{\mathrm Z}.
\]
It leads to 
\begin{subequations}
    \begin{alignat}{11}
        \label{schro_res1}&\ds \left(i\partial_{t}U+\frac12\Delta_{x}U\right)(t,x)=
       \alpha
        \left(\ds\int_{\mathbb T^d\times\mathbb R^n} 
        \sigma_{1}(x-y) \sigma_{2}(z)\Psi(t,y,z)\ud y\ud z\right)U(t,x),\\
        \label{wave_res1}&\Big(\ds\frac1{c^{2}}\partial_{tt}^{2}\Psi-
        \Delta_{z}\Psi\Big)(t,x,z)=-
       \beta
        \sigma_{2}(z)\left(\ds\int_{\mathbb T^d} \sigma_{1}(x-y)|U(t,y)|^{2}
        \ud y\right)
    \end{alignat}
\end{subequations}
together with
 \[\begin{array}{l}
\mathscr E_0=
\ds\frac12\ds\int_{\mathbb T^d}|\nabla_x U |^2 \ud x
+
\ds\frac12\ds\frac{\alpha }{\beta} \ds\iint_{\mathbb T^d\times \mathbb R^n}
\Big(\ds\frac1{c^2}|\partial_t \Psi|^2 +  |\nabla_z\Psi|^2
\Big)
 \ud z\ud x
 \\
 \hspace*{3cm}
 +
 \alpha 
\ds \iint_{\mathbb T^d\times \mathbb R^n}
 |U|^2 \sigma_2\sigma_1\star \Psi \ud z\ud x.
\end{array}\]
This relation allows us to interpret the scaling parameters as weights in the energy balance. 
Now, for notational convenience, we decide to work with $\sqrt{\frac {m}{\mathrm M}}\sqrt{\frac\alpha\beta}\Psi$ 
instead of $\Psi$ and $\sqrt{\frac {\mathrm M}{m}} U$ instead of $U$; it leads to
\eqref{Schro-s}-\eqref{Schro-p} and  \eqref{eq:hamiltonian} with $\gamma=\sqrt{\frac {M}{\mathrm m}}\sqrt{\alpha\beta}$. 
Accordingly, we shall implicitly work with solutions with amplitude of magnitude unity. 
The regime where $c\to \infty$, with $\alpha, \beta$ fixed leads, at least formally,   
to the Hartree system \eqref{Hartree}-\eqref{Poissonz_Hartree}; arguments are sketched in Appendix~\ref{AppB}.
The smallness condition \eqref{small} makes a threshold appear 
on the coefficients in order to guaranty the stability:
since it involves the product $\frac {\mathrm M}{ m}\alpha \beta$, it can be interpreted
as a condition on the strength of the coupling between the particle and the environment, and on the amplitude of the wave function. 
We shall see in the proof that a sharper condition can be derived, expressed by means of the Fourier coefficients of the form function $\sigma_1$.

\section {From Sch\"odinger-Wave to Hartree}
\label{AppB}

In this Section we wish to justify that solutions -- hereafter denoted $U_c$ -- of \eqref{Schro-s}-\eqref{Schro-p}
converge to the solution of \eqref{Hartree}-\eqref{Poissonz_Hartree} as $c\to\infty$.
We adapt the ideas in \cite{dBGV} where this question is investigated for Vlasov equations.
Throughout this section we consider a sequence of initial data $U^{\mathrm{Init}}_{c},\Psi^{\mathrm{Init}}_c,\Pi^{\mathrm{Init}}_c$ such that
\begin{subequations}
    \begin{alignat}{1}
\label{bdd0}
&\ds\sup_{c>0}
\ds\int_{\mathbb T^d}| U^{\mathrm{Init}}_{c}|^2\ud x=M_0<\infty,\\[.3cm]
\label{bdd1}
&\ds\sup_{c>0}
\ds\int_{\mathbb T^d}|\nabla_x U^{\mathrm{Init}}_{c}|^2\ud x=M_1<\infty,
\\[.3cm]
\label{bdd2}
&\ds\sup_{c>0}\left\{\ds\frac1{2c^2}\ds\iint_{\mathbb T^d\times\mathbb R^n}|\Pi^{\mathrm{Init}}_c|^2\ud z\ud x
+\ds\frac1{2}\ds\iint_{\mathbb T^d\times\mathbb R^n}
|\nabla_z\Psi^{\mathrm{Init}}_c|^2\ud z\ud x\right\}=M_2<\infty,
\\[.3cm]
\label{bdd3}
&\ds\sup_{c>0}\ds\iint |U^{\mathrm{Init}}_{c}|^2 \sigma_1\star \sigma_2|\Psi^{\mathrm{Init}}_c|\ud z\ud x=M_3<\infty.
\end{alignat}\end{subequations}
There are several direct consequences of these assumptions:
\begin{itemize}
\item The total energy is initially bounded uniformly with respect to $c>0$,
\item 
In fact, we shall see that the last assumption can be deduced from the previous ones.
\item 
Since the $L^2$ norm of $U_c$ is conserved by the equation, we already know that 
\[\text{$U_c$ is bounded in $L^\infty(0,\infty;L^2(\mathbb T^d))$}.\]
\end{itemize}

Next, we reformulate 
the expression of the potential, separating the contribution due to  the 
initial data of the wave equation and the self-consistent part. 
By using the linearity of the wave equation, we can split
\[
\Phi_c=\Phi_{\mathrm{Init},c}+\Phi_{\mathrm{Cou},c}\]
where $\Phi_{\mathrm{Init},c}$ is defined from the free-wave equation on $\mathbb R^n$
 and initial data $\Psi^{\mathrm{Init}}_c,\Pi^{\mathrm{Init}}_c$:
 \begin{equation}\label{freew}
 \begin{array}{l}
 \ds\frac1{c^2}\partial^2_{tt}\Upsilon_c-\Delta_z \Psi=0,
 \\
 (\Upsilon_c,\partial_t \Upsilon_c)\big|_{t=0}=(\Psi^{\mathrm{Init}}_c,\Pi^{\mathrm{Init}}_c).
 \end{array}
 \end{equation}
 Namely, we set
 \[\begin{array}{lll}
 \Phi_{\mathrm{Init},c}(t,x)
 &=&\ds\int_{\mathbb R^n}\sigma_2(z)\sigma_1\star\Upsilon_c(t,x,z)\ud z
 \\[.3cm]
 &=&\ds\int_{\mathbb R^n}
 \Big(
 \cos(c|\xi|t) \sigma_1\star \overline{\widehat \Psi^{\mathrm{Init}}_c(x,|\xi)}
 +\ds\frac{\sin(c|\xi|t}{c|\xi|}  \sigma_1\star\overline{ \widehat \Psi^{\mathrm{Init}}_c(x,|\xi)}\Big)
 \ds\frac{\widehat \sigma_2(\xi)\ud\xi}{(2\uppi)^n}.
\end{array} \]
Accordingly $\widetilde \Psi_c=\Psi_c-\Upsilon_c$ satisfies
 \begin{equation}\label{freew2}
 \begin{array}{l}
 \ds\frac1{c^2}\partial^2_{tt}\widetilde \Psi_c-\Delta_z \widetilde \Psi_c-=-\gamma\sigma_2\sigma_1\star |U_c|^2,
 \\
 (\widetilde \Psi_c,\partial_t \widetilde \Psi_c)\big|_{t=0}=(0,0).
 \end{array}\end{equation}
 and we get
 \[\begin{array}{lll}
 \Phi_{\mathrm{Cou},c}(t,x)
 &=&\gamma\ds\int_{\mathbb R^n}\sigma_2(z)\sigma_1\star\widetilde \Psi_c(t,x,z)\ud z
 \\[.3cm]
 &=&\gamma^2 c^2\ds\int_0^t \ds\int_{\mathbb R^n}
\ds\frac{\sin(c|\xi|s)}{c|\xi|}  \Sigma\star | U_{c}|^2(t-s,x)
|\widehat \sigma_2(\xi)|^2 \ds\frac{\ud\xi}{(2\uppi)^n}\ud s
\\
\\[.3cm]
 &=&\gamma^2\ds\int_0^{ct} \underbrace{\left(\ds\int_{\mathbb R^n}
\ds\frac{\sin(\tau |\xi|)}{|\xi|} 
|\widehat \sigma_2(\xi)|^2 \ds\frac{\ud\xi}{(2\uppi)^n}\right)}_{=p(\tau)}
 \Sigma\star | U_{c}|^2(t-\tau/c,x)\ud \tau
,
 \end{array}\]
 where it is known that the kernel $p$ is integrable on $[0,\infty)$ \cite[Lemma~14]{dBGV}.

\begin{lemma}
There exists a constant $M_w>0$ such that
\[\ds\sup_{c,t,x}|\Phi_{\mathrm{Init},c}(t,x)|\leq M_w,\qquad
\ds\sup_{c,t,x}|\Phi_{\mathrm{Cou},c}(t,x)|\leq M_w.\]
\end{lemma}

\begin{Proof}
Combining the Sobolev embedding theorem (mind the condition $n\geq 3$) and 
the standard  energy conservation for the free linear wave equation, we obtain 
$$
\|  \Upsilon_{c}\|_{L^\infty(0,\infty;L^2(\mathbb T^d ; L^{2n/(n-2)} (\R^n)))}
\leq C \| \nabla_z \Upsilon_{c} \|_{L^\infty(0,\infty;L^2(\mathbb T^d \times \mathbb R^n))} \leq C\sqrt{2 M_2} .$$
 Applying H\"older's 
 inequality, we are thus led to:
 \begin{equation}
\label{Phi0}
|\Phi_{\mathrm{Init},c}(t,x)|\leq C \|\sigma_2\|_{L^{2n/(n+2)}(\mathbb R^n)} \|\sigma_1\|_{L^2(\mathbb R^d)} \sqrt{2 M_2},
\end{equation}
which proves the first part of the claim.
Incidentally, it also shows that \eqref{bdd3} is a consequence of  \eqref{bdd0} and \eqref{bdd2}.
Next, we get
\[
|\Phi_{\mathrm{Cou},c}(t,x)|
 \leq
\gamma  \| \Sigma \|_{L^\infty(\mathbb T^d)} \| U_c \|_{L^\infty([0,\infty),L^2(\mathbb T^d))}\ds \int_0^\infty |p(\tau)| \ud \tau 
.\]
\end{Proof}

\begin{coro}
There exists a constant $M_S>0$ such that
\[\ds\sup_{c,t}\|\nabla U_c(t,\cdot)\|_{L^2(\mathbb T^d)}
\leq M_S.\]
\end{coro}

\begin{Proof}
This is a consequence of the energy conservation (the total energy being bounded by virtue of \eqref{bdd1}-\eqref{bdd3})  where the coupling term
\[
\ds\int_{\mathbb T^d}
(\Phi_{\mathrm{Init},c}+ \Phi_{\mathrm{Cou},c})|U_c|^2\ud x\]
can be dominated by $2M_wM_0$.
\end{Proof}
 
Coming back to 
\begin{equation}\label{eqsc}
\partial_t U_c=-\ds\frac{1}{2i}\Delta_x U_c +\frac{\gamma}i  (\Phi_{\mathrm{Init},c}+ \Phi_{\mathrm{Cou},c})U_c
\end{equation}
we see that $\partial_t U_c$ is bounded in $L^2(0,\infty;H^{-1}(\mathbb T^d))$.
Combining the obtained estimates with Aubin-Simon's lemma \cite[Corollary~4]{Simon}, we deduce that 
\[\text{$U_c$ is relatively compact in 
in $C^0([0,T];L^p(\mathbb T^d))$, $1\leq p<\ds\frac{2d}{d-2}$},\]
for any $0<T<\infty$. Therefore, possibly at the price of extracting a subsequence, we can suppose that $U_c$ converges strongly to $U$ in $C^0([0,T];L^2(\mathbb T^d))$.
It remains to pass to the limit  in \eqref{eqsc}. The difficulty consists in letting $c$ go to $\infty$ in the potential term and to justify the following claim.
\begin{lemma}
For any $\zeta\in C^\infty_c((0,\infty)\times\mathbb T^d)$, we have
\[\ds\lim_{c\to \infty}
\ds\int_0^\infty\ds\int_{\mathbb T^d} (\Phi_{\mathrm{Init},c}+ \Phi_{\mathrm{Cou},c})U_c\zeta\ud x\ud t
=\gamma\kappa \ds\int_0^\infty\ds\int_{\mathbb T^d} \Sigma\star |U_c|^2\ U_c\zeta\ud x\ud t
.\]
\end{lemma}

\begin{Proof}
We expect that $\Phi_{\mathrm{Cou},c}$ converges to $\gamma\kappa\Sigma\star |U|^2$:
\[\begin{array}{l}
\big|\Phi_{\mathrm{Cou},c}(t,x)- \gamma\kappa \Sigma \star |U|^2(t,x)\big|
\\[.3cm]
\qquad= \gamma
 \left| \ds\int_0^{ct} \Sigma \star |U_c|^2(t-\tau/c,x) p(\tau) \ud \tau - \kappa  \Sigma \star |U|^2(t,x) \right| \\[.4cm]
\qquad \leq
\gamma\ds 
 \ \ds\int_0^{ct} \Big|\Sigma \star |U_c|^2(t-\tau/c,x)-\Sigma \star |U|^2(t,x) \Big|\  |p(\tau)| \ud \tau  
+ 
 \gamma\ds \int_{ct} ^\infty |p(\tau)| \ud \tau \times \|\Sigma  \star |U|^2 \|_{L^\infty((0,\infty)\times\mathbb T^d)}
 \\[.3cm]
 \qquad\leq
\gamma\ds \int_0^{ct}  \Sigma \star \big|   |U_c|^2-  |U|^2\big|(t-\tau/c,x) \  |p(\tau)| \ud \tau 
\\[.3cm]\qquad\qquad
+\gamma \ds 
 \int_0^{ct}  \Sigma \star\big| |U|^2(t-\tau/c ,x) - |U|^2(t,x)\big| \  |p(\tau)| \ud \tau
 \\[.3cm]\qquad\qquad
+ 
\gamma\ds \int_{ct} ^\infty |p(\tau)| \ud \tau \  \| \Sigma \|_{L^\infty(\mathbb T^d)} \| U \|_{L^\infty((0,\infty);L^2(\mathbb T^d))}
 .
\end{array} \]
Let us denote by $\mathrm I_c(t,x)$, $\mathrm {II}_c(t,x)$, $\mathrm {III}_c(t)$, the three terms of the right hand side.
Since $p\in L^1([0,\infty))$, for any $t>0$, $\mathrm {III}_c(t)$ tends to 0 as $c\to \infty$, and it is dominated by $\|p\|_{L^1([0,\infty)}  \|\Sigma\|_{L^{\infty}(\mathbb T^d)}M_0$.
Next,  we
have $$\begin{array}{lll}
 | \mathrm {I}_c(t,x) | &\leq& \|p\|_{L^1([0,\infty)}\|\Sigma\|_{L^\infty(\mathbb T^d)} \ds\sup_{s\geq 0}
\ds\int_{\mathbb T^d}
\big| |U_c|^2 -|U|^2 \big|(s,y)\ud y 
\\[.3cm]
 &\leq& \|p\|_{L^1([0,\infty)}\|\Sigma\|_{L^\infty(\mathbb T^d)} \ds\sup_{s\geq 0}
\left(\ds\int_{\mathbb T^d}
 |U_c-U|^2(s,y) \ud y +2\mathrm{Re} \ds\int_{\mathbb T^d} (U_c-U)\overline U (s,y)\ud y \right)
 \end{array}$$
 which also goes to 0 as $c\to \infty$ and is dominated by $2M_0\|p\|_{L^1([0,\infty))}\|\Sigma\|_{L^\infty(\mathbb T^d)}$.
 Eventually, 
 we get
 \[\begin{array}{lll}
 |\mathrm{II}_c(t,x)| 
  \leq \|\Sigma\|_{L^\infty(\mathbb T^d)} \ds 
 \int_0^{ct} \left(\ds\int_{\mathbb T^d} \big| |U|^2(t-\tau/c,y)-|U|^2(t,y)\big|\ud y\right)\ |p(\tau)|\ud \tau.
\end{array} 
 \]
 Since $U\in C^0([0,\infty);L^2(\mathbb T^d))$, with $\|U(t,\cdot)\|_{L^2(\mathbb T^d)}\leq M_0$, 
 we can apply 
the Lebesgue theorem to show that $\mathrm{II}_c(t,x)$ tends to  0 for any $(t,x)$ fixed, and it is dominated by 
$2M_0\|p\|_{L^1([0,\infty))}\|\Sigma\|_{L^\infty(\mathbb T^d)}$.
This allows us to pass to the limit in 
\[\begin{array}{l}
\ds\int_0^\infty\int _ {\mathbb T^d}
\Phi_{\mathrm{Cou},c}U_c\zeta\ud x\ud t
-\kappa\ds\int_0^\infty\int _ {\mathbb T^d} \Sigma\star|U|^2 U\zeta\ud x\ud t
\\[.3cm]
\qquad
=\ds\int_0^\infty\int _ {\mathbb T^d}
\Phi_{\mathrm{Cou},c}(U_c-U)\zeta\ud x\ud t
+
\ds\int_0^\infty\int _ {\mathbb T^d} \Big(\Phi_{\mathrm{Cou},c}-\gamma\kappa \Sigma\star|U|^2 
\Big)U\zeta\ud x\ud t
.
\end{array}\]

It remains to justify that
\[\ds\lim_{c\to \infty}\ds\int_0^\infty\ds\int_{\mathbb T^d} \Phi_{\mathrm{init},c}U_c\zeta\ud x\ud t=0.\] 
 The space variable  $x$ is just a parameter for the free wave equation \eqref{freew}, which is equally satisfied
 by $\sigma_1\star\Upsilon_c$, with initial data $\sigma_1\star(\Psi^{\mathrm{Init}}_c,\Pi^{\mathrm{Init}}_c)$.
  We appeal to the Strichartz estimate for the wave equation, see  \cite[Corollary 1.3]{MT} or \cite[Theorem 4.2, for the case $n=3$]{Sog},which yields
\[\begin{array}{l}
c^{1/p}\left(\ds\int_0^\infty\left(\ds\int_{\mathbb R^n}|\sigma_1\star\Upsilon_c(t,x,y)|^q\ud y\right)^{p/q}\ud t\right)^{1/p}
\\[.3cm]
\qquad
\leq  C 
\left(\ds\frac{1}{c^2}\ds\int_{\mathbb R^n}|\sigma_1\star\Pi^{\mathrm{Init}}_c(x,z)|^2\ud z+ \ds\int_{\mathbb R^n}|\sigma_1\star\nabla_y\Psi^{\mathrm{Init}}_c(x,z)|^2\ud z\right)^{1/2},
\end{array}\]
 for any admissible pair:
\[2\leq p\leq q\leq\infty,\quad
\ds\frac1p+\ds\frac n q=\ds\frac n2-1,\quad
\ds\frac2p+\ds\frac{n-1}q\leq \ds\frac {n-1}2,\quad
(p,q,n)\neq (2,\infty,3).\]
The $L^2$ norm with respect to the space variable of the right hand side is dominated by 
$  \sqrt{\| \sigma_1\|_{L^1(\mathbb T^d)}\ 
M_{2}}$.
It follows that
\[\ds\int_{\mathbb T^d}
\left(\ds\int_0^\infty\left(\ds\int_{\mathbb R^n}|\sigma_1\star\Upsilon_c(t,x,z)|^q\ud z\right)^{p/q}\ud t
\right)^{2/p}\ud x
\leq 
 C^2\| \sigma_1\|_{L^1(\mathbb R^d)}\ M_2\ds\frac{1}{c^{2/p}}
\xrightarrow[c\to \infty]{}0.\]
Repeated use of the H\"older inequality (with $1/p+1/p'=1$) leads to
\[\begin{array}{l}
\left|\ds\int_0^\infty\ds\int_{\mathbb T^d} U_c\zeta  \Phi_{\mathrm{Init},c}\ud x\ud t\right|
\\[.4cm]
\quad\leq 
\left(\ds\int_{\mathbb T^d}
\left(
\ds\int_0^\infty |U_c\zeta(t,x)|^{p'}\ud t\right)^{2/p'}
\ud x\right)^{1/2}
\left(\ds\int_{\mathbb T^d}
\left(
\ds\int_0^\infty |\Phi_{\mathrm{Init},c}(t,x)|^{p}\ud t\right)^{2/p}
\ud x\right)^{1/2}.
\end{array}\]
On the one hand, assuming that $\zeta$ is supported in $[0,R]\times\mathbb T^d$ and $p>2$, we have
\[\begin{array}{lll}
\ds\int_{\mathbb T^d}
\left(
\ds\int_0^\infty |U_c\zeta|^{p'}\ud t\right)^{2/p'}
\ud x
&\leq&
\ds\int_{\mathbb T^d}
\left(
\ds\int_0^R |U_c|^{2}\ud t\right)
\left(
\ds\int_0^R |\zeta|^{2p'/(2-p')}\ud t\right)^{(2-p')/p'}
\ud x
 \\
&\leq &R^{1+(2-p')/p'}\| \zeta \|_{L^\infty((0,\infty)\times \mathbb T^d)}\|U_c\|_{L^\infty((0,\infty);L^2(\mathbb T^d))} 
\end{array}\]
which is thus 
bounded uniformly with respect to $c>0$.
On the other hand, we get
\[\begin{array}{l}
\ds\int_{\mathbb T^d}\left(\ds\int_0^\infty |\Phi_{\mathrm{Init},c}(t,x)|^{p}\ud t\right)^{2/p}
\ud x
=
\ds\int_{\mathbb T^d}\left(\ds\int_0^\infty \Big|
\ds\int_{\mathbb R^n}\sigma_2(z)
\sigma_1\star\Upsilon_{c}(t,x,z)\ud z
\Big|^{p}\ud t\right)^{2/p}\ud x
\\
\qquad
\leq
\|\sigma_2\|_{L^{q'}(\mathbb R^n)}
\ds\int_{\mathbb T^d}\left(\ds\int_0^\infty \Big|
\ds\int_{\mathbb R^n}
|\sigma_1\star\Upsilon_{c}(t,x,z)|^q\ud z
\Big|^{p/q}\ud t\right)^{2/p}\ud x
\end{array}
\]
which is of the order $\mathscr O(c^{-2/p})$.
\end{Proof}

\section{Well-posedness of the Schrödinger-Wave system}
\label{AppC}

The well-posedness of the Schr\"odinger-Wave system is justified by means of 
a fixed point argument.
The method described here works as well for the problem set on $\mathbb R^d$, and it is simpler than the 
approach in \cite{Vi3} since it avoids
the use of ``dual'' Strichartz estimates for the Schr\"odinger and the wave equations.

We define a mapping that associates to a function $(t,x)\in [0,T]\times\mathbb T^d \mapsto V(t,x)\in \mathbb C$:
\begin{itemize}
\item first, the solution $\Psi$ of the linear wave equation 
\[\ds\frac1{c^2}\partial^2_{tt}\Psi-\Delta_z\Psi=-\sigma_2\sigma_1\star |V|^2,\qquad
(\Psi,\partial_t\Psi)\big|_{t=0}=  (\Psi_{0},\Psi_1);\]
\item next, the potential $\Phi=\sigma_1\star\int_{\mathbb R^n}\sigma_2\Psi\ud z$;
\item and finally the solution of the linear Schr\"odinger equation 
\[i\partial_t U+\ds\frac12\Delta _x U=\gamma\Phi U,\qquad U\big|_{t=0}=U^{\mathrm{Init}}.\]
\end{itemize}
These successive steps 
define a mapping $\mathcal S: V\longmapsto U$
and we wish to show that this mapping admits a fixed point in $C^0([0,T];L^2(\mathbb T^d))$, which, in turn,
 provides a solution to the non linear system 
\eqref{Schro-s}-\eqref{Schro-p}.
In this discussion, the initial data $U^{\mathrm {Init}}, \Psi_{0},\Psi_1$  are fixed once for all in the space of finite energy:
\[U^{\mathrm {Init}}\in H^1(\mathbb T^d),\qquad \Psi_{0}\in L^2(\mathbb T^d;\overbigdot H^1(\mathbb R^n)),\qquad
\Psi_1\in L^2(\mathbb T^d\times \mathbb R^n).\]
We observe that $$\ds\frac{\ud}{\ud t}\ds\int_{\mathbb T^d}|U|^2\ud x=0.$$
Hence, the mapping $\mathcal S$ applies the ball $B(0,\|U^{\mathrm {Init}}\|_{L^2(\mathbb T^d)})$ of 
$C^0([0,T];L^2(\mathbb T^d))$ in itself; 
we thus consider $U=\mathcal S(V)$ for $V\in  C^0([0,T];L^2(\mathbb T^d))$ such that $\|V(t,\cdot)\|_{L^2(\mathbb T^d)}
\leq \|U^{\mathrm {Init}}\|_{L^2(\mathbb T^d)}$. 
Moreover, we can split
\[\Psi=\Upsilon+\widetilde\Psi\]
with $\Upsilon$ solution of the free wave equation
\[\ds\frac1{c^2}\partial^2_{tt}\Upsilon-\Delta_z\Upsilon=0,\qquad
(\Upsilon,\partial_t\Upsilon)\big|_{t=0}=  (\Psi_{0},\Psi_1),\]
and
\[\ds\frac1{c^2}\partial^2_{tt}\widetilde\Psi-\Delta_z\widetilde\Psi=0,\qquad
(\Upsilon,\partial_t\widetilde\Psi)\big|_{t=0}=  0.\] 
We write $\Phi=\Phi_{I}+\widetilde \Phi$ for the associated splitting of the potential.
In particular, the standard energy conservation for the wave equation tells us that 
\[\begin{array}{l}
\ds\frac{1}{2c^2}\ds\iint_{\mathbb T^d\times \mathbb R^n}|\partial_t\Upsilon|^2\ud z\ud x+
\ds\frac12\ds\iint_{\mathbb T^d\times \mathbb R^n}|\nabla_z\Upsilon|^2\ud z\ud x
\\[.3cm]
\qquad=\ds\frac{1}{2c^2}\ds\iint_{\mathbb T^d\times \mathbb R^n}|\Psi_1|^2\ud z\ud x+
\ds\frac12\ds\iint_{\mathbb T^d\times \mathbb R^n}|\nabla_z\Psi_0|^2\ud z\ud x=M_2\end{array}\]
holds. It follows that
\[
|\Phi_I(t,x)|\leq C\|\sigma_{2}\|_{L^{2n/(n+2}(\mathbb R^n)}\|\sigma_1\|_{L^2(\mathbb T^d)}\sqrt{2M_2}
\]
by using Sobolev's embedding.
Next, we obtain
\[\begin{array}{lll}
\widetilde \Phi(t,x)
 &=&\ds\int_{\mathbb R^n}\sigma_2(z)\sigma_1\star\widetilde \Psi(t,x,z)\ud z
 \\[.3cm]
 &=&\gamma\ds\int_0^{ct} \underbrace{\left(\ds\int_{\mathbb R^n}
\ds\frac{\sin(\tau |\xi|)}{|\xi|} 
|\widehat \sigma_2(\xi)|^2 \ds\frac{\ud\xi}{(2\uppi)^n}\right)}_{=p(\tau)}
 \Sigma\star | V|^2(t-\tau/c,x)\ud \tau
,
\end{array}\]
which  thus satisfies
$$\ds\sup_{x\in\mathbb T^d}|\widetilde \Phi(t,x)|\leq
\gamma\|\Sigma\|_{L^\infty(\mathbb T^d)}
\ds\int_0^{ct} |p(\tau)|
\left(\ds\int_{\mathbb T^d} |V|^2(t-\tau/c,y)\ud y\right)
\ud\tau.$$
In particular $$|\widetilde\Phi(t,x)|\leq \gamma\|\Sigma\|_{L^\infty(\mathbb T^d)}
\|p\|_{L^1((0,\infty))}\|V\|_{C^0([0,T];L^2(\mathbb T^d))}
\leq 
 \gamma\|\Sigma\|_{L^\infty(\mathbb T^d)}
\|p\|_{L^1((0,\infty))}\|U^{\mathrm{Init}}\|_{L^2(\mathbb T^d)}
$$ lies in $L^\infty((0,T)\times\mathbb T^d)$, and thus $\Phi\in L^\infty((0,T)\times\mathbb R^d)$.
This observation guarantees that $U=\mathcal; S(V)$ is well-defined.

Thus, let us pick $V_1,V_2$ in this ball of $C^0([0,T];L^2(\mathbb T^d))$ and consider $U_j=\mathcal S(V_j)$. We have
\[i\partial_t (U_2-U_1)+\ds\frac12\Delta _x (U_2-U_1)=\gamma\Phi_2 (U_2-U_1)+\gamma(\Phi_2-\Phi_1)U_1,\qquad
(U_2-U_1)\big|_{t=0}=0.\]
It follows that
\[\begin{array}{l}
\ds\frac{\ud}{\ud t}\ds\int_{\mathbb T^d} |U_2-U_1|^2\ud x
=2\gamma\mathrm{Im}\left(\ds\int_{\mathbb T^d}(\Phi_2-\Phi_1)\overline U_1(U_2-U_1)\ud x\right)
\\[.3cm]
\quad\leq 
2\gamma \|U_1\|_{L^2(\mathbb T^d)} \ \|U_2-U_1\|_{L^2(\mathbb T^d)} \ \|\Phi_2-\Phi_1\|_{L^\infty(\mathbb T^d)}
=2\gamma \|U_1\|_{L^2(\mathbb T^d)} \ \|U_2-U_1\|_{L^2(\mathbb T^d)} \ \|\widetilde\Phi_2-\widetilde\Phi_1\|_{L^\infty(\mathbb T^d)}
\\[.3cm]
\quad\leq
2\gamma^2 \|\Sigma\|_{L^\infty(\mathbb T^d)}
\|U^{\mathrm{Init}}\|_{L^2(\mathbb T^d)} \ \|U_2-U_1\|_{L^2(\mathbb T^d)} \ 
\ds\int_0^{ct} |p(\tau)|
\left(\ds\int_{\mathbb T^d}\big| |V_2|^2-|V_1|^2\big|(t-\tau/c,y)\ud y\right)
\ud\tau.
\end{array}\]
We use the elementary estimate $$
\ds\int_{\mathbb T^d}
\big||V_2|^2-|V_1|^2\big|\ud y=
\ds\int_{\mathbb T^d}\big||V_2-V_1|^2+2\mathrm{Re}(V_2-V_1)V_1\big|\ud y
\leq \|V_2-V_1\|^2_{L^2(\mathbb T^d)}+2\|V_2-V_1\|_{L^2(\mathbb T^d)} \ \|V_1\|_{L^2(\mathbb T^d)}.
$$ 
Combining this with Cauchy-Schwarz and Young inequalities, we arrive at
\[\begin{array}{l}\ds\frac{\ud}{\ud t}\ds\int_{\mathbb T^d} |U_2-U_1|^2\ud x
\\[.3cm]
\quad
\leq
2\gamma^2 \|\Sigma\|_{L^\infty(\mathbb T^d)}
\|U^{\mathrm{Init}}\|_{L^2(\mathbb T^d)}
\left(
2\|U^{\mathrm{Init}}\|_{L^2(\mathbb T^d)}\ds\int_0^{ct} |p(\tau)|  \|V_2-V_1\|^2(t-\tau/c)_{L^2(\mathbb T^d)}\ud\tau
\right.
\\[.3cm]
\qquad\qquad\qquad\qquad\left.+
\|U_2-U_1\|_{L^2(\mathbb T^d)} 
2\|U^{\mathrm{Init}}\|_{L^2(\mathbb T^d)}\ds\int_0^{ct} |p(\tau)|\|V_2-V_1\|(t-\tau/c)_{L^2(\mathbb T^d)}\ud \tau 
\right)
\\[.3cm]
\quad
\leq
2\gamma^2 \|\Sigma\|_{L^\infty(\mathbb T^d)}
\|U^{\mathrm{Init}}\|^2_{L^2(\mathbb T^d)}
\left(
\|U_2-U_1\|_{L^2(\mathbb T^d)}^2
\right.
\\[.3cm]
\qquad\qquad\qquad\qquad\left.+ 
(2+\|p\|_{L^1((0.\infty)})\ds\int_0^{ct} |p(\tau)|  \|V_2-V_1\|^2(t-\tau/c)_{L^2(\mathbb T^d)}\ud\tau
\right).
\end{array}\]
Set $L=2\gamma^2 \|\Sigma\|_{L^\infty(\mathbb T^d)}
\|U^{\mathrm{Init}}\|^2_{L^2(\mathbb T^d)}$.
We deduce that
\[
\|U_2-U_1\|(t)_{L^2(\mathbb T^d)}^2
\leq
(2+\|p\|_{L^1((0.\infty)})L \ds\int_0^t 
e^{L(t-s)} \ds\int_0^{cs} |p(\tau)|  \|V_2-V_1\|^2(s-\tau/c)_{L^2(\mathbb T^d)}\ud\tau\ud s.\]
We use this estimate for $0\leq t\leq T<\infty$ and we obtain
 \[
\|U_2-U_1\|(t)_{L^2(\mathbb T^d)}^2
\leq
(4+\|p\|_{L^1((0.\infty)})L Te^{LT} \|p\|_{L^1((0.\infty)}\ds\sup_{0\leq s\leq T} 
\|V_2-V_1\|^2(s)_{L^2(\mathbb T^d)}.
\]
Hence for $T$ small enough, $\mathcal S$ is a contraction in $C^0([0,T];L^2(\mathbb T^d))$, and consequently it admits a unique fixed point.
 Since the fixed point still has its $L^2$ norm equal to $\|U^{\mathrm{Init}}\|_{L^2(\mathbb T^d)}$, the solution can be extended 
 on the whole interval $[0,\infty)$. 
The argument can be adapted to handle the Hartree system.

\section*{Data availability statement}
Data sharing not applicable to this article as no datasets were generated or analysed during the current study.

%\bibliography{SchroPW}
%\bibliographystyle{plain}

\end{document}